\documentclass[11pt,reqno]{amsart}
\usepackage{amssymb,mathrsfs,graphicx,enumerate}
\usepackage{amsmath,amsfonts,amssymb,amscd,amsthm,bbm,mathtools}

\usepackage{tikz}
\usetikzlibrary{arrows.meta, positioning}
\usepackage[retainorgcmds]{IEEEtrantools}
\usepackage{colortbl}
\usepackage{cite}
\usepackage{comment}
\usepackage{subcaption}
\usepackage{nccmath}
\usepackage{bm}
\usepackage{upgreek}
\usepackage[hidelinks,linktoc=all,bookmarksnumbered=true,bookmarksopen=true]{hyperref}
\title[Orbital stability for the equivariant CS--LL equation]%
{Orbital stability for the equivariant Chern--Simons gauged Landau--Lifshitz equation in $\mathbb{R}^2$}

\author[Jin]{Guanghui Jin}
\address[Guanghui Jin]{\newline Department of Mathematics \newline Yanbian University, Yanji, 133002, China}
\email{jinguanghui@ybu.edu.cn}

\author[Seok]{Jinmyoung Seok} 
\address[Jinmyoung Seok]{\newline Department of Mathematics Education \newline Seoul National University, Seoul 08826, Republic of Korea}
\email{jmseok@snu.ac.kr}

\author[Yuan]{Jianjun Yuan}
\address[Jianjun Yuan]{\newline School of Artificial Intelligence and Information Technology
\newline 
Nanjing University of Chinese Medicine (NJUCM)
Nanjing, 210023, China}
\email{yjj@njucm.edu.cn}

\newcommand{\ip}[2]{\langle #1,#2\rangle}
\newcommand{\n}{\vec n}

\begin{document}

	\newtheorem{thm}{Theorem}[section]
	\newtheorem{lemma}[thm]{Lemma}
	\newtheorem{corollary}[thm]{Corollary}
	\newtheorem{prop}[thm]{Proposition}
	\newtheorem{remark}[thm]{Remark}
	\newtheorem{definition}[thm]{Definition}
	\newtheorem{example}[thm]{Example}
    \numberwithin{equation}{section}
    \renewcommand*{\theHequation}{\theHsection.\arabic{equation}}
	\newcommand{\R}{{\mathbb R}}
	
	\newcommand{\Z}{\mathbb Z}
	\newcommand{\N}{\mathbb N}
	\newcommand{\bbs}{\mathbb S}
	\newcommand{\bbp}{\mathbb P}
	\newcommand{\kp}{\kappa}
        \newcommand{\gm}{\gamma}
        \newcommand{\p}{\partial}
	\newcommand{\mi}{\mathrm i}
	\newcommand{\ddiv}{\textrm{div}}
	\newcommand{\bn}{\bf n}
	\newcommand{\rr}[1]{\rho_{{#1}}}
	\newcommand{\thh}{\theta}
	\def\charf {\mbox{{\text 1}\kern-.24em {\text l}}}
	\renewcommand{\arraystretch}{1.5}
    \newcommand{\D}{\mathcal D}
\newcommand{\pr}{\partial_r}

	

\begin{abstract}
We study the equivariant Chern--Simons gauged Landau--Lifshitz equation in
$\mathbb R^2$.  For equivariance index $m\ge2$ and a potential
$V\in C^{3,1}([0,2])$ satisfying $V(0)=V'(0)=0$, we first establish local well-posedness for the Cauchy problem at
the $H^2$ level, with the sphere constraint written in the natural affine
space $u-\vec n\in H^2$.  For nonnegative $V$ we then consider an
$m$--equivariant ground state
of $H_\omega=E+\omega M$, $\omega>0$.  Under a constrained
nondegeneracy condition for the polar Hessian, the conserved energy and mass
give orbital stability in $H^1$.  A Gustafson--Shatah type continuation
argument, in which the Chern--Simons feedback is treated as a tame nonlocal
zeroth-order perturbation, upgrades this conditional stability to global
$H^2$ existence for all $H^2$ initial data sufficiently close to the ground
state orbit in $H^1$.  For the quadratic potential we also verify the
nondegeneracy condition in a nonempty interior near-BPS regime approaching
the self-dual point $(\kappa,\omega)=(1,0)$.  This yields a concrete family
of globally defined orbitally stable ground states.
\end{abstract}

\maketitle

\tableofcontents
\bigskip

\section{Introduction}\label{sec:1}

\subsection{Background and motivation}
The Landau--Lifshitz equation governs the evolution of magnetization in
ferromagnetic media and is closely related to Schr\"odinger maps through its
geometric precession structure. In two spatial dimensions, physically
relevant models often include interactions with gauge fields, and
Chern--Simons theories are particularly significant. The Chern--Simons
action is topological and produces a flux attachment mechanism that is
central in planar phenomena such as anyonic statistics and effective vortex
dynamics. From an analytic viewpoint, the Chern--Simons coupling introduces
a nonlocal feedback between the spin field and the gauge potential through
elliptic constraints, and this creates an additional difficulty compared
with the ungauged Landau--Lifshitz or Schr\"odinger-map flows.

In this work we investigate the Chern--Simons gauged Landau--Lifshitz system
on $\mathbb R^2$ in an equivariant class and develop a stability theory near
a solitary wave.  The main analytic point is that, after the equivariant
Coulomb reduction, the Chern--Simons feedback does not change the principal
Schr\"odinger-map part of the equation.  It instead appears through radial
nonlocal coefficients which can be estimated without derivative loss.

The differentiated-field approach goes back to Chang--Shatah--Uhlenbeck
\cite{CSU}; its equivariant implementation near harmonic maps was developed
by Gustafson--Kang--Tsai \cite{GKT}.  Our continuation argument follows the
localized-wave strategy of Gustafson--Shatah \cite{GS}, with the additional
Chern--Simons source estimated explicitly below.  These references provide
the geometric and dispersive framework for the present equation.

\subsection{The Chern--Simons gauged Landau--Lifshitz system}
Let $u:\mathbb R^{2+1}\to\mathbb S^2$ be the spin field and
$A_\mu:\mathbb R^{2+1}\to\mathbb R$, $\mu=0,1,2$, an abelian gauge
potential. We fix the north pole $\vec n=(0,0,1)$. The gauge covariant derivative and curvature are
\[
D_\mu u=\partial_\mu u+A_\mu(\vec n\times u),\qquad
F_{\mu\nu}=\partial_\mu A_\nu-\partial_\nu A_\mu.
\]
Here $\partial_0=\partial_t$, and a repeated spatial index $j$ is summed over $j=1,2$.  We consider
\begin{align}
D_0u&=u\times D_jD_ju-V'(1-\vec n\cdot u)(\vec n\times u),
\qquad |u|=1,\label{A-1-1}\\
F_{01}&=\langle\vec n\times u,D_2u\rangle,\label{A-1-2}\\
F_{02}&=-\langle\vec n\times u,D_1u\rangle,\label{A-1-3}\\
F_{12}&=1-u_3.\label{A-1-4}
\end{align}
The boundary condition is
\begin{equation}\label{eq:boundary}
\lim_{|x|\to\infty}u(x,t)=\vec n.
\end{equation}
A principal example is
\[
V(s)=\kappa s^2,\qquad \kappa>0.
\]
For general potentials the flow results use $V\in C^{3,1}([0,2])$:
$V$ has three continuous derivatives and $V^{(3)}$ is Lipschitz.  We use
\[
 \|V\|_{C^{3,1}([0,2])}
 =\sum_{j=0}^3\|V^{(j)}\|_{L^\infty([0,2])}
  +\sup_{s\ne t}\frac{|V^{(3)}(s)-V^{(3)}(t)|}{|s-t|}.
\]
The one-map coefficient estimates need only the $C^3$ part of this norm;
the Lipschitz seminorm is used for differentiated differences.

Throughout we use the Coulomb gauge
\begin{equation}\label{eq:coulomb}
\nabla\cdot A=0,
\end{equation}
and fix its residual temporal freedom by
\begin{equation}\label{eq:A0-normalization}
A_0(r,t)\longrightarrow0\qquad(r\to\infty).
\end{equation}

\subsection{Equivariance, energy, and the ground-state orbit}
For an integer $m\ge1$, an $m$--equivariant map has the form
\[
u(r,\chi)=R_{-m\chi}v(r),
\]
where the target rotation about the $x_3$--axis is
\[
 R_\alpha:=
 \begin{pmatrix}
  \cos\alpha&-\sin\alpha&0\\
  \sin\alpha& \cos\alpha&0\\
  0&0&1
 \end{pmatrix}.
\]
We work in
the topological component
\[
v(0)=-\vec n,\qquad v(\infty)=\vec n.
\]
The variational energy space is
\[
\begin{split}
\mathscr E_m=\{u=R_{-m\chi}v(r):\;&|v|=1,\quad
 u-\vec n\in H^1(\mathbb R^2;\mathbb R^3),\quad
 v\in H^1_{\rm loc}(0,\infty;\mathbb R^3),\\
 &\lim_{r\downarrow0}v(r)=-\vec n,\quad
 \lim_{r\to\infty}v(r)=\vec n\}.
\end{split}
\]
The limits are imposed on the locally absolutely continuous radial
representative; they are not point traces of an arbitrary two-dimensional
$H^1$ function.  For $s\ge2$ the flow space is
\[
 \mathcal X_m^s=\{u\in\mathscr E_m:u-\vec n\in H^s(\mathbb R^2)\}.
\]
All Sobolev distances between sphere-valued maps are affine distances:
$u-u_*\in H^s$ means the ordinary vector-valued Sobolev norm of the
difference, not that $u$ itself is in $L^2$.
The energy $E$ and mass $M$ are defined in Section~\ref{sec:reduction};
for $\omega>0$ we set
\[
 H_\omega(u)=E(u)+\omega M(u).
\]
An $m$--equivariant ground state is a minimizer of $H_\omega$ over
$\mathscr E_m$ (with finite energy).  The variational space is therefore
fixed at the $H^1$ level, independently of the later flow regularity.
Lemma~\ref{lem:great-circle-reduction} identifies it with the scalar
minimization in \cite{JHS}; the ground-state regularity then places its
minimizers in $\mathcal X_m^2$.
If $\widetilde u$ is a ground state of $H_\omega$, its symmetry orbit is
\[
\mathcal O(\widetilde u)=\{R_\gamma\widetilde u:\gamma\in\mathbb R\},
\qquad
d_{H^s}(u,\mathcal O(\widetilde u))
=\inf_{\gamma\in\mathbb R}\|u-R_\gamma\widetilde u\|_{H^s}.
\]

\subsection{Main results}
Our first result is a local theory which is independent of any ground-state
or nondegeneracy assumption.

\begin{thm}[Equivariant $H^2$ local well-posedness]\label{thm:lwp}
Let $m\ge2$ and let $V\in C^{3,1}([0,2])$ satisfy $V(0)=V'(0)=0$.  For every
$u_0\in\mathcal X_m^2$ there exist maximal numbers
$T^{-},T^{+}\in(0,\infty]$ and a unique Coulomb-gauge solution of
\eqref{A-1-1}--\eqref{A-1-4}, normalized by
\eqref{eq:A0-normalization}, such that
\[
u-\vec n\in C((-T^{-},T^{+});H^2(\mathbb R^2))
\cap C^1((-T^{-},T^{+});L^2(\mathbb R^2)).
\]
The solution depends continuously on the initial data in $H^2$, and the
energy and mass are conserved.  Moreover,
\[
T^{+}<\infty\quad\Longrightarrow\quad
\limsup_{t\uparrow T^{+}}\|u(t)-\vec n\|_{H^2}=\infty,
\]
and the analogous statement holds at $-T^{-}$.  More precisely, for every
$R<\infty$ the local existence time can be chosen uniformly on the ball
$\|u_0-\vec n\|_{H^2}\le R$, with constants depending only on $m$,
$R$, and $\|V\|_{C^{3,1}([0,2])}$.
\end{thm}

The second theorem combines constrained orbital stability with a
Gustafson--Shatah type continuation argument.  Let
\[
\widetilde u=R_{-m\chi}(\sin\widetilde\theta,0,
\cos\widetilde\theta)
\]
be an $m$--equivariant ground state of $H_\omega$, $\omega>0$.  The polar
Hessian $L_{\widetilde\theta}^s$ and its form domain $\mathcal H$ are
introduced in Section~\ref{sec:stability}.  We assume the constrained
nondegeneracy condition
\begin{equation}\tag{ND}\label{ND-intro}
\ker L_{\widetilde\theta}^s\cap
\{f\in\mathcal H:\langle f,\sin\widetilde\theta\rangle_{L^2(2\pi r\,dr)}=0\}
=\{0\}.
\end{equation}

\begin{thm}[Global existence and orbital stability]\label{thm:global-stability}
Let $m\ge2$, let $V\in C^{3,1}([0,2])$ be nonnegative with
$V(0)=V'(0)=0$, and let $\omega>0$.  Suppose that
$\widetilde u\in\mathcal X_m^2$ is an $m$--equivariant ground state of
$H_\omega$ satisfying \emph{(ND)}.  Then for every $\varepsilon>0$ there
exists $\delta>0$ such that, whenever $u_0\in\mathcal X_m^2$ and
\[
d_{H^1}(u_0,\mathcal O(\widetilde u))<\delta,
\]
the maximal solution of Theorem~\ref{thm:lwp} is global and satisfies
\[
\sup_{t\in\mathbb R}d_{H^1}(u(t),\mathcal O(\widetilde u))
<\varepsilon.
\]
In addition, for every finite $T>0$,
\[
\sup_{|t|\le T}\|u(t)-\vec n\|_{H^2}<\infty.
\]
\end{thm}

The nondegeneracy hypothesis is not vacuous.  For the quadratic potential
we verify it in an entire interior sector approaching the self-dual point.

\begin{corollary}[Near-BPS ground states]\label{cor:near-BPS}
Let $m\ge2$ and $V(s)=\kappa s^2$.  There exists $\lambda_*>0$ with the
following property.  For every compact interval
$\mathcal K=[\lambda^{-},\lambda^{+}]\Subset(0,\lambda_*)$ there is $\delta_{\mathcal K}>0$ such
that, whenever
\[
0<\delta<\delta_{\mathcal K},\qquad
\kappa=1-\delta,\qquad
\omega=\delta\lambda,\qquad \lambda\in \mathcal K,
\]
the set of $m$--equivariant ground states of $H_\omega$ is nonempty and
every such ground state satisfies \emph{(ND)}.  Consequently all conclusions
of Theorem~\ref{thm:global-stability} hold near each of these ground states.
\end{corollary}

The restriction $\lambda\in \mathcal K\Subset(0,\lambda_*)$ is important.  The cusp
regime $\omega/(1-\kappa)\to0$ is not covered by the present compactness
argument; see Remark~\ref{remark-critical-cusp-not-covered}.

\subsection{Organization of the paper}
Section~\ref{sec:reduction} gives the equivariant Coulomb reduction and the
conservation laws. Section~\ref{sec:lwp-proof} proves
Theorem~\ref{thm:lwp}. Section~\ref{sec:stability} establishes constrained
Hessian coercivity and the conditional $H^1$ orbital stability estimate.
Section~\ref{sec:continuation} proves the $H^2$ continuation theorem in an
$H^1$ tube. Section~\ref{sec:bps} verifies the nondegeneracy hypothesis for
an interior near-BPS family and completes the proof of
Theorem~\ref{thm:global-stability} and Corollary~\ref{cor:near-BPS}.

\section{Equivariant Coulomb reduction and conservation laws}\label{sec:reduction}

Let
\[
u(r,\chi,t)=R_{-m\chi}v(r,t).
\]
For $x=(r\cos\chi,r\sin\chi)$, we use the standard orthonormal polar frame
\[
e_r=(\cos\chi,\sin\chi),\qquad
e_\chi=(-\sin\chi,\cos\chi),
\]
so that $x=re_r$ and $x^\perp=re_\chi$.
The Coulomb condition and \eqref{A-1-4} determine the spatial gauge field
uniquely among equivariant fields vanishing at the origin:
\begin{equation}\label{eq:A-h}
A=\frac{h(r,t)}r e_\chi,\qquad
h(r,t)=\int_0^r s(1-v_3(s,t))\,ds.
\end{equation}
Indeed $F_{12}=r^{-1}\partial_rh=1-v_3$.  With the normalization
\eqref{eq:A0-normalization}, equations \eqref{A-1-2}--\eqref{A-1-3}
give
\begin{equation}\label{eq:k-def-main}
A_0=k[v],\qquad
k[v](r,t)=\int_r^\infty
\frac{(v_1(s,t)^2+v_2(s,t)^2)}{s}\,(h[v](s,t)-m)\,ds,
\end{equation}
We use $k=A_0$ throughout.  The scalar notation in \cite{JHS}
uses the opposite sign for the tail integral in
\eqref{eq:k-def-main}; this accounts for the sign difference
between the two radial formulas.  Substitution into \eqref{A-1-1} yields the
single radial equation
\begin{equation}\label{eqn-single-m}
v_t=v\times\Delta_rv-
\left\{\frac{(h-m)^2}{r^2}v_3+V'(1-v_3)+k\right\}(\vec n\times v),
\qquad
\Delta_r=\partial_r^2+r^{-1}\partial_r.
\end{equation}
The remaining polar component of the time--space curvature equations is
automatic.  Indeed, differentiating \eqref{eq:A-h} in time and using
\eqref{eqn-single-m} gives
\[
h_t(r,t)=-\int_0^r s(v_t)_3\,ds
=-r\,(v\times v_r)_3,
\]
which is exactly the equivariant form of the remaining relation in
\eqref{A-1-2}--\eqref{A-1-3}.  Thus the gauge variables are determined
by $v$ and there is no independent gauge evolution in the equivariant
Coulomb class.

The energy and mass become
\begin{align}
E(u)&=2\pi\int_0^\infty
\left\{\frac12|v_r|^2+
\frac{(h-m)^2}{2r^2}(v_1^2+v_2^2)+V(1-v_3)\right\}r\,dr,
\label{eq:energy-radial}\\
M(u)&=2\pi\int_0^\infty(1-v_3)r\,dr.
\label{eq:mass-radial}
\end{align}
Since $|u-\vec n|^2=2(1-u_3)$,
\begin{equation}\label{eq:mass-L2}
\|u-\vec n\|_2^2=2M(u).
\end{equation}

\begin{lemma}[Conservation laws]\label{lem:conservation}
Let $u$ be an equivariant solution of \eqref{eqn-single-m} on an
interval $I\ni0$, with
$u-\vec n\in C(I;H^2(\mathbb R^2))\cap C^1(I;L^2(\mathbb R^2))$.
Then, for every $t\in I$,
\[
E(u(t))=E(u(0)),\qquad M(u(t))=M(u(0)).
\]
\end{lemma}
\begin{proof}
We first prove the identities for smooth solutions by radial
integration by parts and then pass to the $H^2$ solutions by approximation.
The calculation includes the contribution generated by eliminating the gauge field.  Since $(\vec n\times v)_3=0$ and
\[
 r(v\times\Delta_rv)_3
 =\partial_r\{r(v\times v_r)_3\},
\]
one has
\[
 \frac d{dt}M(u(t))
 =-2\pi\bigl[r(v\times v_r)_3\bigr]_{0}^{\infty}=0.
\]
For a smooth finite-energy profile the boundary term vanishes first after
cutting off at $\varepsilon<r<R$ and then sending
$\varepsilon\downarrow0$, $R\uparrow\infty$; the required limits follow
from the radial and angular parts of the Dirichlet energy.

For the energy, let $\xi\cdot v=0$ be a smooth compactly supported tangent
variation and write
\[
 h'_v[\xi](r)=-\int_0^r s\xi_3(s)\,ds.
\]
Differentiating the gauge term and using Fubini gives
\begin{align*}
 &2\pi\int_0^\infty
 \frac{(h-m)(1-v_3^2)}r h'_v[\xi] \,dr\\
 &\qquad=-2\pi\int_0^\infty
 \left\{\int_r^\infty
 \frac{(h(s)-m)(1-v_3(s)^2)}s\,ds\right\}\xi_3(r)r\,dr
 =-2\pi\int_0^\infty k\xi_3r\,dr.
\end{align*}
Consequently the full first variation is
\begin{equation}\label{eq:first-variation-energy}
 DE(v)[\xi]
 =-2\pi\int_0^\infty
 \left\langle
 \Delta_rv+\left\{k+\frac{(h-m)^2}{r^2}v_3
              +V'(1-v_3)\right\}\vec n,\xi
 \right\rangle r\,dr.
\end{equation}
Taking $\xi=v_t$ in \eqref{eq:first-variation-energy}, and using
\eqref{eqn-single-m}, yields
\[
 \begin{aligned}
 \frac d{dt}E(u(t))
 &=-2\pi\int_0^\infty
 \left\langle
 \Delta_rv+\left\{k+\frac{(h-m)^2}{r^2}v_3
                +V'(1-v_3)\right\}\vec n,\,v_t
 \right\rangle r\,dr=0.
 \end{aligned}
\]
because $v_t$ is the cross product of $v$ with the first vector
in the inner product.
The preceding calculation is justified for noncompact profiles by radial
cutoffs; all cutoff errors vanish by the finite-energy Hardy bounds and the
normalization $A_0(\infty)=0$.

For an $H^2$ solution, the local approximation and uniqueness
construction in Section~\ref{sec:lwp-proof} identifies it, on each
compact time interval, with the strong $C_tH_x^1$ limit of smooth
solutions.  The averaging and backward-Hardy estimates
\eqref{eq:h-over-r2-difference} and
\eqref{eq:backward-Hardy-difference} give convergence of the gauge-energy
terms.  Hence both identities pass to the limit.
\end{proof}

It is useful to compare \eqref{eqn-single-m} with the ordinary
$m$--equivariant Schr\"odinger-map equation
\[
v_t=v\times\Delta_rv-\frac{m^2}{r^2}v_3(\vec n\times v).
\]
Thus the CSLL equation can be written on $\mathbb R^2$ as
\begin{equation}\label{eq:forced-SM-main}
u_t=u\times\Delta u+\Gamma[u](\vec n\times u),
\end{equation}
where
\begin{equation}\label{eq:Gamma-main}
\Gamma[v]=\frac{h(2m-h)}{r^2}v_3-V'(1-v_3)-k.
\end{equation}
Whenever $u(r,\chi,t)=R_{-m\chi}v(r,t)$, the notation
$\Gamma[u](r,\chi,t)$ means the radial scalar coefficient
$\Gamma[v](r,t)$; it is independent of $\chi$.
When applied to a spatial scalar or vector field $f=f(x)$, the notation
$D^2f=(\partial_j\partial_k f)_{j,k=1}^2$ denotes the ordinary Cartesian
Hessian, not the gauge covariant derivative $D_\mu$ (or the second
variation of a functional).
This exact decomposition is the starting point for both the local theory
and the continuation argument.

\section{\texorpdfstring{$H^2$}{H2} local well-posedness}\label{sec:lwp-proof}

We prove Theorem~\ref{thm:lwp}.  The argument follows the modified-field
approach used for equivariant Schr\"odinger maps, with two additional
facts: the coefficient \eqref{eq:Gamma-main} is tame on $H^2$ bounded
sets, and the full Chern--Simons source enters the differentiated equation
without derivative loss.  The detailed coefficient estimates needed here
are proved again, in a form uniform on an $H^1$ tube, in
Lemmas~\ref{lem:tame} and \ref{lem:frame-source} below.  Those are static
coefficient estimates: their proofs use no local existence, stability,
or continuation theorem, so their use here does not presuppose the
dynamical conclusions of Section~\ref{sec:continuation}.

\subsection{Geometric notation and local recovery}
We distinguish the two connections from the outset.  The notation
$D_\mu$ denotes the gauge derivative defined in the introduction,
whereas $\nabla_\mu^u$ denotes the pullback Levi--Civita connection.
For sphere-valued $u$ and a tangent section $Y\cdot u=0$, set
\[
 P_uz=z-(z\cdot u)u,\qquad
 \nabla_\mu^uY=P_u\partial_\mu Y,\qquad J_uY=-u\times Y.
\]
Here $\mu$ can be $t$, a Cartesian spatial index, or $r,\chi$; in the
last case $\partial_\chi$ is the coordinate angular derivative, not
$r^{-1}\partial_\chi$.  We abbreviate $J_u$ to $J$ when the map is fixed.
Ordinary Cartesian derivatives are denoted by $\nabla_x$ or $D^j$;
these are not covariant derivatives.  Norms without a domain are spatial
norms on $\mathbb R^2$.  A radial amplitude uses $2\pi r\,dr$, unless a
harmless constant factor is explicitly suppressed.

Choose a radial parallel orthonormal frame $(e,Je)$,
$\nabla_r^ue=0$, and extend it equivariantly.  For a tangent field define
\[
 [Y]_e=\langle Y,e\rangle+i\langle Y,Je\rangle.
\]
In particular, at $\chi=0$ define the radial amplitudes
\begin{equation}\label{eq:W-local}
 W=u_r+\frac1rJu_\chi,\qquad q=[W]_e,
 \qquad W=(\Re q)e+(\Im q)Je.
\end{equation}
When planar norms are needed, the radial amplitude $q$ is lifted
explicitly as $e^{i(m-1)\chi}q(r)$; no separate lift operator is used.
We also use the canonical-frame coordinate
$\zeta=[P_u\vec n]_e$ for the tangential projection of $\vec n$.
For a smooth time-dependent frame its total temporal connection is the
real scalar $\alpha_{\rm tot}$ defined by
$\nabla_t^ue=\alpha_{\rm tot}Je$.  Its decomposition into the pure-flow
part $\alpha$ and the forcing part $\alpha_Z$ is made in
Lemma~\ref{lem:exact-q-equation}.
The restriction $m\ge2$ is used precisely through the Hardy estimate
\begin{equation}\label{eq:hardy-q-local}
\left\|\frac q r\right\|_{L^p(rdr)}
\le C_m\|\nabla(e^{i(m-1)\chi}q)\|_{L^p(\mathbb R^2)},
\qquad 1<p<\infty.
\end{equation}

\begin{lemma}[Short-time $H^2$ recovery]\label{lem:local-recovery}
For every $R<\infty$ there are $\rho_R>0$ and $C_R<\infty$ such that the
following holds.  If $u_*,u\in\mathcal X_m^2$ satisfy
\[
\|u_*-\vec n\|_{H^2}\le R,
\qquad \|u-u_*\|_{H^1}\le\rho_R,
\]
then
\begin{equation}\label{eq:local-recovery}
\|u-\vec n\|_{H^2}
\le C_R\left(1+\|u-\vec n\|_2+
\|e^{i(m-1)\chi}q\|_{H^1}\right).
\end{equation}
\end{lemma}
\begin{proof}
The Bogomolny derivative satisfies the exact identity
\begin{equation}\label{eq:tension-local}
P_u\Delta u=\nabla_r^u W+\frac1r(1-J\nabla_\chi^u)W.
\end{equation}
By equivariance and \eqref{eq:hardy-q-local},
\[
\|P_u\Delta u\|_2\le C_m
\|e^{i(m-1)\chi}q\|_{\dot H^1}.
\]
Since $|u|=1$,
\[
\Delta u=P_u\Delta u-|\nabla u|^2u,
\qquad \|D^2u\|_2=\|\Delta u\|_2.
\]
Put $w=u-u_*$.  The two-dimensional Gagliardo--Nirenberg inequality gives
\[
\|\nabla u\|_4^2
\le C_R+C\|\nabla w\|_2\|D^2w\|_2
\le C_R+C\rho_R(\|D^2u\|_2+R).
\]
Choose $\rho_R$ so that the last multiple of $\|D^2u\|_2$ is absorbed.
Adding the $L^2$ part proves \eqref{eq:local-recovery}.
\end{proof}

\subsection{The exact modified-field equation}
For the local theory we use the following exact equation for the derivative field.
The calculation is first made for smooth solutions and then used for the
regularized approximants.  As in the preceding subsection let
\[
 J=-u\times,\qquad
 W=u_r+\frac1rJu_\chi=(\Re q)e+(\Im q)Je.
\]
where $\nabla_r^u e=0$.  We now remove the remaining phase freedom once and for
all.  Fix $\mathbf e_1=(1,0,0)\in T_{\vec n}\mathbb S^2$.  At $\chi=0$
we choose the unique radial parallel frame satisfying
\begin{equation}\label{eq:canonical-frame-infinity}
 \nabla_r^u e=0,\qquad \lim_{r\to\infty}e(r,0,t)=\mathbf e_1,
\end{equation}
and extend it equivariantly by
$e(r,\chi,t)=R_{-m\chi}e(r,0,t)$.  We call this the \emph{canonical
frame}.  Since the limiting vector is fixed in time, its temporal
connection satisfies $\alpha_{\rm tot}(\infty,t)=0$.  In particular, no
 time-dependent rephasing is made in the local difference estimates.

The following comparison lemma is the reason that the normalization at
infinity is compatible with the $H^1$ difference theory.

\begin{lemma}[Canonical radial-frame comparison]
\label{lem:canonical-frame-comparison}
For every $R<\infty$ there are $r_R>0$, $\rho_R>0$, and $C_R<\infty$
with the following property.  Let $u,\bar u\in\mathcal X_m^2$,
\[
 \|u-\vec n\|_{H^2}+\|\bar u-\vec n\|_{H^2}\le R,
 \qquad \|u-\bar u\|_{H^1}\le\rho_R,
\]
and let $e,\bar e$ be their canonical frames.  Then
\begin{equation}\label{eq:canonical-frame-comparison}
 \|u-\bar u\|_\infty+\sup_{r>0}|e(r,0)-\bar e(r,0)|
 \le C_R\|u-\bar u\|_{H^1}.
\end{equation}
\end{lemma}

\begin{proof}
We give the argument because it is precisely what prevents a hidden
relative temporal phase from entering the difference equation.  The
radial lemma and the $H^2$ bound give a radius $R_0=R_0(R)$ such that both
maps lie in a fixed north-pole chart for $r\ge R_0$.  In that chart use the
smooth unit tangent frame
\[
 E_1^{+}(u)=\mathbf e_1-\frac{u_1}{1+u_3}(u+\vec n),
 \qquad E_2^{+}(u)=J_uE_1^{+}(u).
\]
A direct calculation gives for its radial connection coefficient
$a_u^{+}=\langle \nabla_r^u E_1^{+}(u),E_2^{+}(u)\rangle$
\begin{align}
 |a_u^{+}|&\le C\sqrt{u_1^2+u_2^2}\,|u_r|,                                      \label{eq:north-connection-bound}\\
 |a_u^{+}-a_{\bar u}^{+}|
 &\le C\bigl[(\sqrt{u_1^2+u_2^2}+\sqrt{\bar u_1^2+\bar u_2^2})|u_r-\bar u_r|
       +\sqrt{(u_1-\bar u_1)^2+(u_2-\bar u_2)^2}(|u_r|+|\bar u_r|)\bigr].             \label{eq:north-connection-difference}
\end{align}
The right sides are integrable in $dr$.  Indeed, equivariance gives
\[
\begin{aligned}
 &\sum_{j=1}^2\left(\|u_j/r\|_{L^2(rdr)}
           +\|\bar u_j/r\|_{L^2(rdr)}\right)\le C_mR,\\
 &\sum_{j=1}^2\|(u_j-\bar u_j)/r\|_{L^2(rdr)}
        \le C_m\|u-\bar u\|_{H^1}.
\end{aligned}
\]
so weighted Cauchy--Schwarz applies to
\eqref{eq:north-connection-bound}--\eqref{eq:north-connection-difference}.
Writing
\[
 e=\cos\psi_u\,E_1^{+}(u)+\sin\psi_u\,E_2^{+}(u),
 \qquad
 \psi_u(r)=\int_r^\infty a_u^{+}(s)\,ds,
\]
and similarly for $\bar u$, therefore defines exactly the normalization
\eqref{eq:canonical-frame-infinity} and yields
\[
 \sup_{r\ge R_0}|e(r,0)-\bar e(r,0)|
 \le C_R\|u-\bar u\|_{H^1}.
\]
Choose $r_R>0$ so that the uniform $H^2$ bound and the common
value $u(0)=\bar u(0)=-\vec n$ put both maps in a fixed south-pole
chart on $(0,r_R)$.  The transverse difference obeys
\[
\begin{aligned}
 \sup_{r>0}\sum_{j=1}^2|u_j(r)-\bar u_j(r)|^2
 &\le2\sum_{j=1}^2\|(u_j-\bar u_j)/r\|_{L^2(rdr)}
           \|u_{j,r}-\bar u_{j,r}\|_{L^2(rdr)}\\
 &\le C_m\|u-\bar u\|_{H^1}^2.
\end{aligned}
\]
On either pole chart the third component is a smooth function of the
transverse components, with uniformly bounded derivative.  On the fixed
annulus $[r_R,R_0]$ use one-dimensional Sobolev embedding instead.
Thus $\|u-\bar u\|_\infty\le C_R\|u-\bar u\|_{H^1}$ globally.
On the annulus, subtract the extrinsic parallel-transport equations
\[
 e_r=-(u_r\cdot e)u,\qquad
 \bar e_r=-(\bar u_r\cdot\bar e)\bar u.
\]
Their $L^1(dr)$ coefficients are bounded uniformly, so Gronwall extends
the frame bound down to $r_R$.  On the inner interval use
\[
 E_1^{-}(u)=\mathbf e_1-\frac{u_1}{1-u_3}(u-\vec n),\qquad
 E_2^{-}(u)=J_uE_1^{-}(u).
\]
Its connection satisfies the analogues of
\eqref{eq:north-connection-bound}--\eqref{eq:north-connection-difference}.
The phase difference at $r_R$ is already controlled, and its integral
from $r_R$ to any $r>0$ is bounded by
\[
 C_R\sum_{j=1}^2\Bigl[
 (\|u_j/r\|_2+\|\bar u_j/r\|_2)\|u_r-\bar u_r\|_2
 +\|(u_j-\bar u_j)/r\|_2(\|u_r\|_2+\|\bar u_r\|_2)
 \Bigr].
\]
This proves the asserted global comparison.  The construction uses
absolutely continuous radial representatives and integrable connection
coefficients, so it applies directly to $H^2$ maps.  For smooth maps the
frame depends smoothly on time on compact annuli; the normalized tail
is obtained from the same convergent connection integral.
\end{proof}

\begin{lemma}[Fixed-profile frame comparison in an $H^1$ tube]
\label{lem:canonical-frame-H1-tube}
Fix $u_*\in\mathcal X_m^2$, let $e_*$ be its canonical frame,
and let $q_*$ be the amplitude in \eqref{eq:W-local} for $u_*$.  There are constants $\rho_*>0$ and
$C_*<\infty$, depending only on $u_*$ and $m$, such that every
$u\in\mathcal X_m^2$ with $\|u-u_*\|_{H^1}\le\rho_*$ satisfies
\begin{equation}\label{eq:canonical-frame-H1-tube}
 \|u-u_*\|_{L^\infty}
 +\sup_{r>0}|e(r,0)-e_*(r,0)|
 +\|e^{i(m-1)\chi}(q-q_*)\|_{L^2}
 \le C_*\|u-u_*\|_{H^1}.
\end{equation}
\end{lemma}

\begin{proof}
Write $u=R_{-m\chi}v$, $u_*=R_{-m\chi}v_*$, and
$w=v-v_*$.  All one-dimensional norms in this proof use the measure
$r\,dr$.  The angular part of the $H^1$ distance controls
$\sum_{j=1}^2\|w_j/r\|_2$.  Consequently
\[
 \sup_{r>0}\sum_{j=1}^2|w_j(r)|^2
 \le 2\sum_{j=1}^2\int_0^\infty |w_j||w_j'|\,dr
 \le 2\sum_{j=1}^2\|w_j/r\|_2\|w_j'\|_2
 \le C_m\|u-u_*\|_{H^1}^2.
\]
Choose fixed radii $0<a<A$ so that the Dirichlet energies of $u_*$ on
$(0,a)$ and $(A,\infty)$ are small.  The elementary estimates
\[
 \begin{aligned}
 |v(r)+\vec n|^2&\le \frac1m\int_0^r
 \left(|v_s|^2+\frac{m^2(v_1^2+v_2^2)}{s^2}\right)s\,ds,\\
 |v(r)-\vec n|^2&\le \frac1m\int_r^\infty
 \left(|v_s|^2+\frac{m^2(v_1^2+v_2^2)}{s^2}\right)s\,ds
 \end{aligned}
\]
follow from $|v_3'|\le \sqrt{v_1^2+v_2^2}|v'|$ and the endpoint values.
After decreasing $\rho_*$, both maps therefore lie in a fixed south-pole
chart on $(0,a)$ and a fixed north-pole chart on $(A,\infty)$.
On either chart $v_3$ is a smooth function of $(v_1,v_2)$, with uniformly
bounded derivative.  The preceding transverse estimate controls $w_3$
there.  On $[a,A]$, one-dimensional Sobolev embedding controls the whole
profile difference.  Thus
\[
 \sup_{r>0}|v(r)-v_*(r)|\le C_*\|u-u_*\|_{H^1}.
\]

On the north-pole tail use $E_1^{+}$ and its connection from
\eqref{eq:north-connection-bound}--\eqref{eq:north-connection-difference}.
The connection difference has the integrable bound
\[
\begin{aligned}
 \int_A^\infty |a_u^{+}-a_{u_*}^{+}|\,dr
 &\le C_*\sum_{j=1}^2\bigl[\|v_j/r\|_2\|w'\|_2\\
 &\hspace{28mm}+\|w_j/r\|_2(\|v'\|_2+\|v_*'\|_2)\bigr]\\
 &\le C_*\|u-u_*\|_{H^1}.
\end{aligned}
\]
Together with the common frame value at infinity, this bounds
$e-e_*$ uniformly on $[A,\infty)$.  On $[a,A]$, subtract the two
extrinsic parallel-transport equations and use Gronwall's inequality;
their $L^1(dr)$ coefficients and the one-dimensional $H^1$ difference
are controlled by the $H^1$ tube alone.

To reach the origin, use the south-pole frame
\[
 E_1^{-}(v)=\mathbf e_1-\frac{v_1}{1-v_3}(v-\vec n),
 \qquad E_2^{-}(v)=J_vE_1^{-}(v).
\]
Its connection satisfies the same estimates as the north-pole frame,
with $1-v_3$ in the denominator.  Choose the phase difference at $r=a$
using the already established frame bound.  Integrating the connection
difference from $a$ to $r$ costs at most
\[
 C_*\sum_{j=1}^2\Bigl[
 (\|v_j/r\|_2+\|v_{*,j}/r\|_2)\|w'\|_2
 +\|w_j/r\|_2(\|v'\|_2+\|v_*'\|_2)\Bigr].
\]
The estimates are uniform for $0<r<a$.  This proves the global frame
bound without using any higher derivative.

Finally, since $J_u-J_{u_*}$ is bounded pointwise by $|u-u_*|$,
\[
 \|W[u]-W[u_*]\|_2
 \le C\|u-u_*\|_{H^1}
    +\|u-u_*\|_\infty\|r^{-1}(u_*)_\chi\|_2
 \le C_*\|u-u_*\|_{H^1}.
\]
Taking coordinates in the two frames adds at most
$C\|W[u_*]\|_2(\|e-e_*\|_\infty+\|u-u_*\|_\infty)$.
The angular lift preserves the $L^2$ norm, which proves
\eqref{eq:canonical-frame-H1-tube}.
\end{proof}

\begin{lemma}[Forward continuity of the canonical modified field]
\label{lem:forward-modified-field}
For every $R<\infty$ there are $\rho_R>0$ and $C_R<\infty$ such that
$u,\bar u\in\mathcal X_m^2$ with
$\|u-\vec n\|_{H^2}+\|\bar u-\vec n\|_{H^2}\le R$ and
$\|u-\bar u\|_{H^1}\le\rho_R$ satisfy
\begin{equation}\label{eq:forward-modified-field}
 \|e^{i(m-1)\chi}(q-\bar q)\|_{H^1}\le C_R\|u-\bar u\|_{H^2}.
\end{equation}
The canonical field of every $H^2$ map belongs to $H^1$, with a bound
on each $H^2$ ball.  Smooth equivariant sphere-valued maps are dense in
$\mathcal X_m^2$ in the affine $H^2$ topology.
\end{lemma}
\begin{proof}
We first record the polar estimates that control the apparent singular
terms.  If $F=e^{ik\chi}f(r)$, its Euclidean Hessian satisfies
\[
 \frac1{2\pi}\|D^2F\|_2^2
 =\int_0^\infty\left\{|f''|^2
  +2k^2\left|\frac{f'}r-\frac f{r^2}\right|^2
  +\left|\frac{f'}r-\frac{k^2f}{r^2}\right|^2\right\}r\,dr.
\]
For $|k|\ge2$ the last two quadratic expressions are positive definite
in $f'/r,f/r^2$, and hence
\[
 \|f'/r\|_2+\|f/r^2\|_2\le C_k\|D^2F\|_2.
\]
For $k=0$ the same identity controls $f'/r$; no bound for $f/r^2$ is
claimed.  These formulas apply to the transverse and third components,
respectively, of an equivariant map or a difference of two such maps.
They hold by density at the $H^2$ level.

We will also use a useful consequence of the nonzero angular modes.
For $F=e^{ik\chi}f$ with $k\ne0$,
\[
 \|f\|_\infty^2\le2\|f'\|_{L^2(rdr)}\|f/r\|_{L^2(rdr)}
 \le C_k\|\nabla_xF\|_2^2.
\]
The complex first derivatives of $u_1+iu_2$ have modes $-m\pm1$,
and those of $u_3$ have modes $\pm1$.  None is zero for $m\ge2$.
Consequently
\begin{equation}\label{eq:equivariant-gradient-infty}
 \|du\|_\infty\le C_m\|D^2u\|_2.
\end{equation}

Write $w=u-\bar u$ and use radial profiles at $\chi=0$.  In components
\[
 P_u\vec n=(-u_3 v_1,-u_3 v_2,(v_1^2+v_2^2)),\qquad
 \partial_r(P_u\vec n)_j=-(\partial_r u_3)v_j-u_3\partial_r v_j\quad(j=1,2).
\]
The polar Hessian estimates, the bounds $|v|=1$, and
$\|w\|_\infty\lesssim\|w\|_{H^2}$ give
\[
 \left\|\frac{P_u\vec n-P_{\bar u}\vec n}{r^2}\right\|_2
 +\left\|\frac{(P_u\vec n-P_{\bar u}\vec n)_r}{r}\right\|_2
 \le C_R\|w\|_{H^2}.
\]
For instance, the first two components are the sums of
$-u_3(v_j-\bar v_j)/r^2$ and
$-(u_3-\bar u_3)\bar v_j/r^2$, $j=1,2$; the third component uses
$((v_1^2+v_2^2)-(\bar v_1^2+\bar v_2^2))/r^2$.  For the derivative quotient,
subtract $-(\partial_r u_3)v_j-u_3\partial_r v_j$ componentwise and use the radial $\partial_r u_3/r$ and $v_r/r$ bounds.
The analogous first-order bound, together with the definition of $W$,
therefore yields
\[
 \|W[u]-W[\bar u]\|_2
 +\|\partial_r(W[u]-W[\bar u])\|_2
 +\|(W[u]-W[\bar u])/r\|_2
 \le C_R\|w\|_{H^2}.
\]
The corresponding one-map quantities are bounded on the same ball.
Since the frame is parallel and $W$ is tangent,
\[
 q_r=\langle W_r,e\rangle+i\langle W_r,Je\rangle.
\]
Subtract this identity and $q=[W]_e$, and apply the global comparison
of Lemma~\ref{lem:canonical-frame-comparison}.  Finally
\[
 \|e^{i(m-1)\chi}q\|_{H^1}^2
 =2\pi\int_0^\infty
       (|q|^2+|q_r|^2+(m-1)^2|q/r|^2)r\,dr.
\]
This proves \eqref{eq:forward-modified-field} and the one-map bound.

For density, radially mollify $u-\vec n$ in Cartesian variables and add
$\vec n$ back.  Uniform convergence follows from $H^2$ convergence;
for sufficiently small mollification radii the resulting vector has
length bounded away from zero.  Normalize it by $z\mapsto z/|z|$.
Composition in $H^2$ gives convergence, radial convolution preserves
equivariance, and normalization preserves it as well.  The transverse
components at zero vanish and the third component remains negative,
so the normalized value there is $-\vec n$; the limit at infinity is
$\vec n$.  The approximants belong to every finite Sobolev class.
The forward estimate and the connection-integral construction show that
the limiting field is exactly the canonical field of the limit map.
In particular no inverse recovery estimate is used in this argument.
\end{proof}

With the present orientation convention, a constant target rotation obeys
\[
 q[R_\sigma u]=e^{-i\sigma}q[u],\qquad
 \zeta[R_\sigma u]=e^{-i\sigma}\zeta[u].
\]
Indeed $R_\sigma e$ has limiting vector $R_\sigma\mathbf e_1$; the
constant frame phase that restores the boundary value $\mathbf e_1$
multiplies both complex coordinates by $e^{-i\sigma}$.
A genuinely time-dependent factor will be introduced only in
Section~\ref{sec:continuation}, where its derivative is kept explicitly
in the equation.

The following elementary identities will be used repeatedly.  Since
$u(r,\chi)=R_{-m\chi}v(r)$, every equivariant tangent field $Y$ satisfies
\[
 J\nabla_\chi^u Y=-mu_3Y.
\]
Moreover $Ju_\chi=mP_u\vec n$, and hence
\begin{equation}\label{eq:local-reconstruction-identities}
 [u_r]_e=q-\frac m r\zeta,\qquad
 \zeta_r=-u_3\left(q-\frac m r\zeta\right),\qquad
 \frac{\partial_r u_3}{r}=\Re\left(\frac\zeta r\,\overline q\right)
                 -m\left|\frac\zeta r\right|^2.
\end{equation}

Set
\[
  (TF)(r)=\int_r^\infty F(s)\frac{ds}{s}.
\]
The operator $T$ is the backward radial Hardy operator; it is bounded on
$L^p((0,\infty);r\,dr)$ for $1<p<\infty$.

\begin{lemma}[Exact modified equation]\label{lem:exact-q-equation}
Let $u$ be a smooth equivariant solution of
\[
 u_t=u\times\Delta u+Z,\qquad Z\cdot u=0,
\]
and let $q$ be defined above.  Decompose the total temporal frame connection as
$\alpha_{\rm tot}=\alpha+\alpha_Z$, where $\alpha$ is produced by the pure
Schr\"odinger-map vector field and $\alpha_Z$ is produced by $Z$.
Then
\begin{equation}\label{eq:alphaZ-local}
 \alpha_Z(r)=-\int_r^\infty\langle JZ(s),u_s(s)\rangle\,ds
\end{equation}
and the forcing contributes to the $q$-equation through
\begin{equation}\label{eq:G-local-exact}
 \mathfrak G_u[Z]
 =\left[\nabla_r^u Z+\frac1rJ\nabla_\chi^u Z\right]_e-i\alpha_Zq.
\end{equation}
For the pure Schr\"odinger-map part one has
\begin{equation}\label{eq:exact-q-radial}
 iq_t-q_{rr}-\frac1rq_r+\frac{(m-1)^2}{r^2}q=\mathcal V[u,q]q,
\end{equation}
where
\begin{equation}\label{eq:Vq-exact}
 \mathcal V[u,q]
 =\alpha+\frac{m\partial_r u_3}{r}
 +\frac{(m-1)^2-(1+m u_3)^2}{r^2},
\end{equation}
and the temporal connection is given exactly by
\begin{equation}\label{eq:alpha-exact}
 \begin{aligned}
 \alpha&=\frac12|q|^2-\frac m r\Re(q\overline\zeta)
          +\alpha_{\mathrm{nl}},\\
 \alpha_{\mathrm{nl}}(r)&:=-T\left(|q|^2-
              \frac{2m}{r}\Re(q\overline\zeta)\right)(r).
 \end{aligned}
\end{equation}
Consequently the full CSLL equation satisfies, with $Q(r,\chi)=e^{i(m-1)\chi}q(r)$,
\begin{equation}\label{eq:exact-Q-CSLL}
 (i\partial_t-\Delta)Q
 =e^{i(m-1)\chi}\mathcal V[u,q]q
 +ie^{i(m-1)\chi}\mathfrak G_u[\Gamma[u](\vec n\times u)].
\end{equation}
\end{lemma}

\begin{proof}
Put $p=[u_t^{\rm SM}]_e$, where
$u_t^{\rm SM}=u\times\Delta u$.  The identity
\eqref{eq:tension-local} and $J\nabla_\chi^u W=-mu_3W$ give
\[
 p=-i\left(\partial_r+\frac{1+m u_3}{r}\right)q.
\]
If $\nabla_t^u e=\alpha Je$ for the pure flow, compatibility of the pullback
connection gives
\[
 q_t+i\alpha q
 =\left(\partial_r-\frac{m u_3}{r}\right)p.
\]
Combining the last two formulas yields
\[
 i q_t
 =q_{rr}+\frac1rq_r+\frac{m\partial_r u_3}{r}q
  -\frac{(1+m u_3)^2}{r^2}q+\alpha q,
\]
which is exactly \eqref{eq:exact-q-radial}--\eqref{eq:Vq-exact}.

We next compute $\alpha$.  Since $\nabla_r^u e=0$, the curvature identity on
$\mathbb S^2$ gives
\[
 \alpha_r=\langle Ju_t^{\rm SM},u_r\rangle
 =\langle P_u\Delta u,u_r\rangle.
\]
From \eqref{eq:local-reconstruction-identities} and
$[P_u\Delta u]_e=q_r+(1+m u_3)q/r$ we obtain
\[
 \alpha_r
 =\frac12(|q|^2)_r-\frac m r\partial_r\Re(q\overline\zeta)
  +\frac1r|q|^2-\frac m{r^2}\Re(q\overline\zeta).
\]
Therefore
\[
 \partial_r\left(\alpha-\frac12|q|^2+\frac mr\Re(q\overline\zeta)\right)
 =\frac1r\left(|q|^2-\frac{2m}{r}\Re(q\overline\zeta)\right).
\]
The chosen temporal normalization and the $H^1$ decay of the lifted
field give zero boundary value at infinity, and integration proves
\eqref{eq:alpha-exact}.

Finally, if the tangent vector field $Z$ is added to the evolution, then
$p$ is replaced by $p+[Z]_e$.  The curvature identity gives
$(\alpha_Z)_r=\langle JZ,u_r\rangle$ and the normalization at infinity
implies \eqref{eq:alphaZ-local}.  Subtracting the pure equation from the
full compatibility relation gives precisely
\eqref{eq:G-local-exact}.  Multiplication by $i$ produces the second term
on the right-hand side of \eqref{eq:exact-Q-CSLL}.
\end{proof}

The advantage of \eqref{eq:exact-Q-CSLL} is that every coefficient is now
an explicit algebraic or Hardy expression in $q$ and the zeroth-order
geometric variables.  No derivative is hidden in an unspecified
``standard modified equation''.  We record the estimates in the exact
form used below.

\begin{lemma}[Tame and difference estimates for the pure modified
nonlinearity]\label{lem:pure-q-tame}
For every $R<\infty$ there is $C_R<\infty$ with the following property.
Let $u\in\mathcal X_m^2$ be smooth and
$\|u-\vec n\|_{H^2}\le R$, and let $Q=e^{i(m-1)\chi} q$.  Then
\begin{equation}\label{eq:pure-q-tame}
 \left\|e^{i(m-1)\chi}(\mathcal V[u,q]q)\right\|_{W^{1,4/3}}
 \le C_R\left(1+\|Q\|_{W^{1,4}}\right).
\end{equation}
Suppose in addition that $\bar u$ is another such map, both derivative
fields are written in the canonical frames of
Lemma~\ref{lem:canonical-frame-comparison}, and
$\|u-\bar u\|_{H^1}\le\rho_R$ with $\rho_R$ sufficiently small.  With
barred notation for the second map,
\begin{align}
 &\left\|e^{i(m-1)\chi}(\mathcal V[u,q]q)
       -e^{i(m-1)\chi}(\mathcal V[\bar u,\bar q]\bar q)
 \right\|_{W^{1,4/3}}                                      \notag\\
 &\qquad\le C_R
 \left(1+\|Q\|_{W^{1,4}}+\|\bar Q\|_{W^{1,4}}\right)
 \left(\|Q-\bar Q\|_{H^1}
      +\|Q-\bar Q\|_{W^{1,4}}
      +\|u-\bar u\|_{H^1}\right).                         \label{eq:pure-q-difference}
\end{align}
\end{lemma}

\begin{proof}
We give the coefficient calculation because this is the modified-field estimate
which is needed for Theorem~\ref{thm:lwp}.  The first and third identities in
\eqref{eq:local-reconstruction-identities} become
\begin{equation}\label{eq:b-c-identities}
 [u_r]_e=q-\frac{m\zeta}{r},
 \qquad
 \frac{\partial_r u_3}{r}=\Re(\frac{\zeta}{r}\overline q)-m\frac{|\zeta|^2}{r^2}.
\end{equation}
Also
\begin{equation}\label{eq:geom-potential-factor}
 \frac{(m-1)^2-(1+m u_3)^2}{r^2}
 =m\frac{1+u_3}{r^2}\bigl(2(m-1)-m(1+u_3)\bigr).
\end{equation}
On the present $H^2$ ball, the two geometric coefficient bounds
are best kept separate:
\begin{equation}\label{eq:b-hardy-local}
 \begin{aligned}
  \|\frac{\zeta}{r}\|_{L^p}&\le C_{R,p}(1+\|Q\|_{L^p}),&&2\le p<\infty,\\
  \|\frac{1+u_3}{r^2}\|_{L^p}&\le C_{R,p}(1+\|Q\|_{L^{2p}}^2),&&1<p<\infty.
 \end{aligned}
\end{equation}
Here the $p=2$ bound for $\frac{\zeta}{r}$ follows directly from the angular energy.
For $p>2$, choose $r_R>0$ so that $1+u_3\le\varepsilon_R$ on $(0,r_R)$.
This radius is uniform on the $H^2$ ball by the common value $u_3(0)=-1$
and the embedding $H^2\hookrightarrow C^{0,\gamma}$, $\gamma<1$.
The regular reconstruction equation has the integral form
\[
 \frac{\zeta(r)}r=\int_0^1t^m
 \left\{-u_3(tr)q(tr)+\frac{m(1+u_3(tr))}{tr}\zeta(tr)\right\}\,dt.
\]
Its $L^p((0,r_R);r\,dr)$ operator norm is at most
$(m+1-2/p)^{-1}$.  The last term is absorbed by choosing
$\varepsilon_R$ small.  On the exterior $\frac{|\zeta|}{r}\le r^{-1}$, which is in
$L^p$ for $p>2$.  Finally,
\[
 \frac{1+u_3}{r^2}=\frac{|\zeta|^2}{r^2(1-u_3)}
 \le \frac{|\zeta|^2}{r^2}\quad(0<r<r_R),\qquad
 0\le \frac{1+u_3}{r^2}\le2r^{-2}\quad(r\ge r_R).
\]
Applying the bound for $\frac{\zeta}{r}$ at exponent $2p>2$ proves the second line.
This comparison is used only on the inner pole chart;
no global $L^1$ estimate for $\frac{1+u_3}{r^2}$ is asserted.
To control one derivative, on an interval
$(0,r_R)$ on which $\|1+u_3\|_{L^\infty}\le\varepsilon_R$, the second
identity in \eqref{eq:local-reconstruction-identities} is equivalent to
\begin{equation}\label{eq:b-first-order}
 \partial_r\bigl(\frac{\zeta}{r}\bigr)+\frac{m+1}{r^2}\zeta
 =-u_3\frac q r+m(1+u_3)\frac{\zeta}{r^2}.
\end{equation}
The regular inverse of $\partial_r+(m+1)/r$ sends $f$ to
$r^{-m-1}\int_0^r s^{m+1}f(s)\,ds$.  The weighted Hardy inequality
for this inverse reads
\[
\begin{aligned}
&\left\|r^{-m-2}\int_0^r s^{m+1}f(s)\,ds\right\|_{L^4(rdr)}\\
&\quad+\left\|\partial_r\left(r^{-m-1}\int_0^r s^{m+1}f(s)\,ds\right)\right\|_{L^4(rdr)}
 \le C_m\|f\|_{L^4(rdr)}.
\end{aligned}
\]
Taking $\varepsilon_R$ small, the last term in
\eqref{eq:b-first-order} is absorbed, while
$\|q/r\|_4\lesssim_m\|\nabla Q\|_4$ by
\eqref{eq:hardy-q-local}.  On $(r_R,\infty)$ the factors $r^{-1}$ are
bounded, so \eqref{eq:b-first-order} and \eqref{eq:b-hardy-local} give
the same estimate.  Hence
\begin{equation}\label{eq:b-derivative-local}
 \|D\bigl(\frac{\zeta}{r}\bigr)\|_{L^4}
 \le C_R\left(1+\|Q\|_{W^{1,4}}\right).
\end{equation}
On $(0,r_R)$ we use
$\frac{1+u_3}{r^2}=\frac{|\zeta|^2}{r^2(1-u_3)}$.  Thus
\[
 |\partial_r\bigl(\frac{1+u_3}{r^2}\bigr)|\le C_R\bigl(\frac{|\zeta|}{r}|\partial_r\bigl(\frac{\zeta}{r}\bigr)|+\frac{|\zeta|^2}{r^2}|\partial_r u_3|\bigr),
\]
and \eqref{eq:b-c-identities}, \eqref{eq:b-hardy-local}, and
\eqref{eq:b-derivative-local} imply
\begin{equation}\label{eq:d-derivative-local}
 \|D\bigl(\frac{1+u_3}{r^2}\bigr)\|_{L^2}
 \le C_R\left(1+\|Q\|_{W^{1,4}}\right).
\end{equation}
On $(r_R,\infty)$ the same conclusion follows directly from
$\partial_r\bigl(\frac{1+u_3}{r^2}\bigr)=\partial_r u_3/r^2-2(1+u_3)/r^3$.  Here $D$ denotes the Cartesian derivative
after the natural angular lift.  Thus no unestimated polar derivative
occurs.

The local terms in \eqref{eq:Vq-exact} are now controlled by
\eqref{eq:b-c-identities}, \eqref{eq:geom-potential-factor},
\eqref{eq:b-derivative-local}, and \eqref{eq:d-derivative-local}.  For
the nonlocal part of the temporal connection in
\eqref{eq:alpha-exact}, the $L^2(r\,dr)$ boundedness of $T$ gives
\[
 \|\alpha_{\mathrm{nl}}\|_2
 \le C\left\||q|^2-\frac{2m}{r}\Re(q\overline\zeta)\right\|_2\le C_R,
 \qquad (\alpha_{\mathrm{nl}})_r
 =\frac{|q|^2-\frac{2m}{r}\Re(q\overline\zeta)}{r}.
\]
Consequently
\[
 \|\alpha_{\mathrm{nl}}q\|_{4/3}
   \le\|\alpha_{\mathrm{nl}}\|_2\|q\|_4,
 \qquad
 \|\alpha_{\mathrm{nl}}q_r\|_{4/3}
   \le\|\alpha_{\mathrm{nl}}\|_2\|q_r\|_4.
\]
The remaining differentiated term is estimated without differentiating a
Hardy integral:
\begin{align*}
 \|(\alpha_{\mathrm{nl}})_rq\|_{4/3}
 &\le C\left(\|q/r\|_4\|q\|_4^2
       +\|\zeta/r^2\|_4\|q\|_4^2\right)\\
 &\le C_R\left(1+\|Q\|_{W^{1,4}}\right),
\end{align*}
where the $\zeta/r^2$ bound is contained in the first-order Hardy estimate used
in \eqref{eq:b-first-order}.  The angular derivative of the lift is
estimated in the same way, with $q/r$ controlled by
\eqref{eq:hardy-q-local}.  The local cubic terms
$|q|^2q$ and $\frac1r\Re(q\overline\zeta)q$ are handled by the same
$L^4\cdot L^4\cdot L^4$ H\"older estimate.  Thus
\begin{equation}\label{eq:alpha-q-product}
 \|e^{i(m-1)\chi}(\alpha q)\|_{W^{1,4/3}}
 \le C_R\left(1+\|Q\|_{W^{1,4}}\right).
\end{equation}
Together with the local terms this proves \eqref{eq:pure-q-tame}.

For the fixed-time difference estimate we keep the canonical frames.
On $[r_R,\infty)$ the frame comparison lemma, together with the smooth
dependence of $P_u\vec n$ on $u$, controls the coordinate difference.
On $(0,r_R)$, subtract the two equations
\eqref{eq:b-first-order}, use $\lim_{r\downarrow0}\zeta(r)/r
=\lim_{r\downarrow0}\bar\zeta(r)/r=0$, and apply the same
weighted Hardy inverse used above; the small coefficients
$1+u_3$ and $1+\bar u_3$ are absorbed.  These two regions give
\begin{equation}\label{eq:frame-difference-local}
 \|\frac{\zeta}{r}-\frac{\bar\zeta}{r}\|_{L^p}
 +\left\|\frac{(1+u_3)-(1+\bar u_3)}{r^2}\right\|_{L^p}
 \le C_{R,p}\left(\|u-\bar u\|_{H^1}
                  +\|Q-\bar Q\|_{L^{2p}}\right).
\end{equation}
Subtracting the first-order equations \eqref{eq:b-first-order} for
the two maps and applying the same inverse Hardy estimate gives
\begin{align*}
 \|D(\frac{\zeta}{r}-\frac{\bar\zeta}{r})\|_4
 +\|D(\frac{1+u_3}{r^2}-\frac{1+\bar u_3}{r^2})\|_2
 \le C_R\bigl(&\|u-\bar u\|_{H^1}
 +\|Q-\bar Q\|_{H^1}
 +\|Q-\bar Q\|_{W^{1,4}}\bigr).
\end{align*}
Subtracting \eqref{eq:alpha-exact} for the two maps and using linearity
of $T$, together with the explicit identity $(Tf)_r=-f/r$, gives the same
Lipschitz estimate for the products
$\alpha q-\bar\alpha\bar q$ in $W^{1,4/3}$.  Applying the product
rule to the remaining two terms in \eqref{eq:Vq-exact} yields
\eqref{eq:pure-q-difference}.  Notice that the constants depend only on
the $H^2$ ball, not on higher derivatives.
\end{proof}

\begin{lemma}[Short-time difference recovery]\label{lem:local-difference-recovery}
For every $R<\infty$ there are $\rho_R>0$ and $C_R<\infty$ such that, for
two $H^2$ maps on the ball of Lemma~\ref{lem:forward-modified-field} with
$\|u-\bar u\|_{H^1}\le\rho_R$ and both radial frames chosen canonically,
\begin{equation}\label{eq:local-difference-recovery}
 \|u-\bar u\|_{H^2}
 \le C_R\left(\|u-\bar u\|_{H^1}
              +\|Q-\bar Q\|_{H^1}\right).
\end{equation}
\end{lemma}

\begin{proof}
The exact tension identity gives in the canonical frame
\[
 P_u\Delta u=[q_r+(1+m u_3)q/r]_e^{-1},
\]
where $[z]_e^{-1}=(\Re z)e+(\Im z)Je$.  Subtract the two expressions.
The global frame and map comparison from
Lemma~\ref{lem:canonical-frame-comparison}, the one-map bound in
Lemma~\ref{lem:forward-modified-field}, and $m-1\ne0$ imply
\[
 \|P_u\Delta u-P_{\bar u}\Delta\bar u\|_2
 \le C_R(\|Q-\bar Q\|_{H^1}+\|u-\bar u\|_{H^1}).
\]
For the normal part, put $w=u-\bar u$.  By
\eqref{eq:equivariant-gradient-infty},
\[
 \bigl\||du|^2u-|d\bar u|^2\bar u\bigr\|_2
 \le C_R(\|w\|_2+\|dw\|_2).
\]
Subtracting $\Delta u=P_u\Delta u-|du|^2u$ and using
$\|D^2w\|_2=\|\Delta w\|_2$ proves
\eqref{eq:local-difference-recovery}.  All identities extend from the
smooth approximants of Lemma~\ref{lem:forward-modified-field}, so the
same estimate holds for $H^2$ maps with the same constants.
\end{proof}

\subsection{Strichartz estimates for smooth solutions}
For smooth solutions the preceding calculation is classical.  At low
regularity the moving-frame identities are justified by the same
regularization argument as in \cite[Proposition~3.2]{ShimizuLocal}; the
important point here is that the equation to be estimated is now the
explicit equation \eqref{eq:exact-Q-CSLL}.
For reference, the free two-dimensional Strichartz estimate used here is
\begin{equation}\label{eq:inhomogeneous-Strichartz-displayed}
 \|Q\|_{L_t^\infty H_x^1(I)}
 +\|Q\|_{L_t^4W_x^{1,4}(I)}
 \le C_{\rm Str}\bigl(\|Q(t_0)\|_{H^1}
       +\|F\|_{L_t^{4/3}W_x^{1,4/3}(I)}
       +\|G\|_{L_t^1H_x^1(I)}\bigr)
\end{equation}
for a Duhamel solution of $(i\partial_t-\Delta)Q=F+G$ on
$I=[t_0,t_1]$.  The constant is independent of the length of $I$;
spatial derivatives commute with the free propagator.

Set
\[
X_T=\|Q\|_{L_t^\infty H_x^1([-T,T])}
+\|Q\|_{L_t^4W_x^{1,4}([-T,T])}.
\]

We fix the bootstrap constants in an order that does not use a
conclusion as an assumption.  Given the initial radius $R$, let
$\rho_R$ and $L_R=C_{\rm rec}(R)$ be the constants in
Lemma~\ref{lem:local-recovery}, and set
\[
 \begin{split}
 M_R&=\sup_{\|u_0-\vec n\|_{H^2}\le R}\|e^{i(m-1)\chi}q[u_0]\|_{H^1}<\infty,\qquad
 A_R=C_{\rm Str}M_R,\\
 K_R&=2(A_R+1),\qquad
 B_R=2\max\{R,L_R(1+R+K_R)\}.
 \end{split}
\]
The supremum is over $\mathcal X_m^2$ (and is taken to be zero for an
empty ball); it is finite by
Lemma~\ref{lem:forward-modified-field}; $C_{\rm Str}$ is a fixed
Strichartz constant.  In particular $A_R,K_R,B_R$ are selected before
any time length is selected.

\begin{lemma}[Uniform short-time derivative estimate]\label{lem:q-local}
There exists $T_R>0$, depending only on $m$, $R$, and
$\|V\|_{C^{3,1}([0,2])}$, such that the following holds.
Let $u$ be a smooth unregularized solution on $[-T,T]$, $T\le T_R$, with
$\|u(0)-\vec n\|_{H^2}\le R$.  If
\[
 \sup_{|t|\le T}\|u(t)-\vec n\|_{H^2}\le B_R,\qquad
 \sup_{|t|\le T}\|u(t)-u(0)\|_{H^1}\le\rho_R,
\]
then
\begin{equation}\label{eq:q-local-bound}
 X_T\le K_R.
\end{equation}
\end{lemma}
\begin{proof}
The $(4,4)$ Strichartz estimate, together with its $L_t^\infty H_x^1$
component, gives
\begin{align}\label{eq:q-local-strichartz-exact}
 X_T&\le C_{\rm Str}\|Q(0)\|_{H^1}
 +C\|e^{i(m-1)\chi}(\mathcal Vq)\|_{L_t^{4/3}W_x^{1,4/3}}
 +C\|e^{i(m-1)\chi}\mathfrak G_u[\Gamma[u](\vec n\times u)]\|_{L_t^1H_x^1}.
\end{align}
Apply Lemma~\ref{lem:pure-q-tame} on the already fixed ball of radius
$B_R$.  With a constant $C_{B_R}$ that may depend on this radius but
not on $T$,
\begin{equation}\label{eq:pure-q-time}
 \|e^{i(m-1)\chi}(\mathcal Vq)\|_{L_t^{4/3}W_x^{1,4/3}}
 \le C_{B_R}(T^{3/4}+T^{1/2}X_T).
\end{equation}
Lemma~\ref{lem:frame-source} and the bootstrap bound similarly give
\begin{equation}\label{eq:CS-q-time-local}
 \|e^{i(m-1)\chi}\mathfrak G_u[\Gamma[u](\vec n\times u)]\|_{L_t^1H_x^1}
 \le C_{B_R}T(1+X_T).
\end{equation}
Enlarging $C_{B_R}$ to include the Strichartz constants, we obtain
\[
 X_T\le A_R+C_{B_R}(T^{1/2}+T)X_T+C_{B_R}(T^{3/4}+T).
\]
Now choose $T_R\le1$ so that
$C_{B_R}(T_R^{1/2}+T_R)\le1/2$ and
$C_{B_R}(T_R^{3/4}+T_R)\le1$.  Then
$X_T\le2(A_R+1)=K_R$, as claimed.  Dependence on $B_R$ causes no
circularity: that radius was fixed using only the initial ball and the
one-map recovery constant, before these smallness conditions on $T_R$.
\end{proof}

For later difference estimates we fix the geometric distance now.
When $\|u-\bar u\|_\infty<1$, let
$\Pi(x):T_{\bar u(x)}\mathbb S^2\to T_{u(x)}\mathbb S^2$ be
parallel transport along the unique shorter geodesic.  For $a,b\in\mathbb S^2$
with $a\cdot b>-1$, the transport of $Y\in T_a\mathbb S^2$ is
\begin{equation}\label{eq:geodesic-transport-first-use}
 \Pi_{a\to b}Y=Y-\frac{Y\cdot b}{1+a\cdot b}(a+b),
 \qquad Y\cdot a=0.
\end{equation}
Here $a=\bar u(x)$ and $b=u(x)$.  In particular, the denominator
is bounded below by $3/2$ under $\|u-\bar u\|_\infty<1$.
Define the difference energy by
\begin{equation}\label{eq:geometric-difference-energy}
 \mathcal Z_j=\partial_ju-\Pi\partial_j\bar u,\qquad
 G(t)=\frac12\int_{\mathbb R^2}
       \left(|u-\bar u|^2+\sum_{j=1}^2|\mathcal Z_j|^2\right)dx.
\end{equation}
On a fixed equivariant $H^2$ ball, this is locally equivalent to the
squared affine $H^1$ distance.  Indeed
$\partial_j(u-\bar u)=\mathcal Z_j+(\Pi-I)\partial_j\bar u$,
$|\Pi-I|\lesssim|u-\bar u|$, and
\eqref{eq:equivariant-gradient-infty} controls
$\|d\bar u\|_\infty$.
We also record the uniform equivariant comparison
\[
 \|u-\bar u\|_\infty\le C_R\|u-\bar u\|_{H^1}.
\]
For $\|u-\bar u\|_{H^1}\le\rho_R$ this is the map part of
Lemma~\ref{lem:canonical-frame-comparison}.  For larger distances it
follows from $\|u-\bar u\|_\infty\le2\le
2\rho_R^{-1}\|u-\bar u\|_{H^1}$, after increasing $C_R$.
Thus this comparison is valid for every pair in the fixed $H^2$ ball;
it does not assert the false general embedding $H^1(\mathbb R^2)
\subset L^\infty$.

\begin{lemma}[The instantaneous pure geometric difference estimate]
\label{lem:pure-geometric-difference}
Let $m\ge2$ and $R<\infty$.  Suppose that
$u,\bar u\in\mathcal X_m^2$ satisfy
$\max\{\|u-\vec n\|_{H^2},\|\bar u-\vec n\|_{H^2}\}\le R$
and $\|u-\bar u\|_\infty<1$.  Let $\Pi$, $\mathcal Z_i$, and $G$
be defined by \eqref{eq:geometric-difference-energy}.  For smooth
profiles define
\begin{equation}\label{eq:pure-energy-directional-derivative}
 G'_{\rm SM}:=DG(u,\bar u)
   [u\times\Delta u,\,\bar u\times\Delta\bar u].
\end{equation}
Then
\begin{equation}\label{eq:pure-geometric-bound}
 |G'_{\rm SM}|\le C_RG.
\end{equation}
If $u,\bar u$ are equivariant sphere-valued paths on a compact time
interval, with the same uniform $H^2$ and short-geodesic bounds and
$u-\vec n,\bar u-\vec n\in C_tH_x^2\cap C_t^1L_x^2$, then
$G$ is absolutely continuous and
\begin{equation}\label{eq:main-geometric-low-regularity}
 \begin{gathered}
 \mathcal Z_i\in C_tH_x^1,\qquad
 \partial_t\mathcal Z_i\in L_t^\infty H_x^{-1},\\
 \frac12\frac d{dt}\|\mathcal Z_i(t)\|_2^2
 =\langle\nabla_t^u\mathcal Z_i,\mathcal Z_i\rangle_{H^{-1},H^1}
 \quad\text{for a.e. }t.
 \end{gathered}
\end{equation}
Along such paths solving the forced CSLL equation, the pure
Schr\"odinger-map contribution to $G'$ satisfies
\eqref{eq:pure-geometric-bound} for almost every time.
\end{lemma}
\begin{proof}
The argument uses the geometric parallel-transport method of
\cite{M,ShimizuUnique}.  The requisite fixed-profile pairing bounds
and the low-regularity energy identity are proved in
Appendix~\ref{app:instantaneous-geometric}.  The directional derivative
in \eqref{eq:pure-energy-directional-derivative} can equivalently be
evaluated along any pair of sphere-valued curves with the indicated
initial velocities; its value is independent of the curves.

Put $w=u-\bar u$, $\psi_i=\partial_i u$, and
$\bar\psi_i=\partial_i\bar u$.  Let $\Phi_s(a,b)$ be the shorter
constant-speed geodesic from $a$ to $b$, parametrized by $s\in[0,1]$,
and write $\gamma(s)=\Phi_s(\bar u,u)$ and $\Pi_{t,s}$ for transport
along $\gamma$ from parameter $s$ to parameter $t$.
For a section $Y$ along $\bar u$, set
\[
 \mathcal C_jY=\nabla_j^u(\Pi Y)-\Pi\nabla_j^{\bar u}Y.
\]
Differentiation of the parallel-transport equation gives
\begin{equation}\label{eq:transport-commutator-explicit}
 \mathcal C_jY
 =\int_0^1\Pi_{1,s}R^{\mathbb S^2}(\gamma_s,\partial_j\gamma)
                         \Pi_{s,0}Y\,ds,
 \qquad
 R^{\mathbb S^2}(a,b)c=(c\cdot b)a-(c\cdot a)b.
\end{equation}

We make the velocity dependence explicit.  Given tangent endpoint
velocities $\bar V,V$, the variation field of the geodesic is
\begin{equation}\label{eq:geodesic-velocity-linear}
 \mathcal J[\bar V,V](s)
 =\mathrm d_1\Phi_s(\bar u,u)[\bar V]
  +\mathrm d_2\Phi_s(\bar u,u)[V].
\end{equation}
It is the Jacobi field with boundary values $\bar V,V$ and is linear
in this pair.  Its short-geodesic boundary-value estimate, also proved
in Appendix~\ref{app:instantaneous-geometric}, is
\begin{equation}\label{eq:main-velocity-Jacobi-bound}
 |\mathcal J[\bar V,V](s)|\le C(|\bar V|+|V|),\qquad 0\le s\le1.
\end{equation}
Define the temporal transport contribution at these velocities by
\begin{equation}\label{eq:temporal-transport-linear}
 \mathcal C[\bar V,V]Y
 =\int_0^1\Pi_{1,s}
 R^{\mathbb S^2}(\gamma_s,\mathcal J[\bar V,V](s))
                   \Pi_{s,0}Y\,ds.
\end{equation}
This too is linear in the endpoint velocities.  If
$F=u\times\Delta u$, $\bar F=\bar u\times\Delta\bar u$,
$Z=\Gamma[u](\vec n\times u)$, and
$\bar Z=\Gamma[\bar u](\vec n\times\bar u)$, then
\begin{equation}\label{eq:transport-velocity-splitting}
 \mathcal C_t=\mathcal C[\bar u_t,u_t]
 =\mathcal C[\bar F,F]+\mathcal C[\bar Z,Z]
 =:\mathcal C_t^{\rm SM}+\mathcal C_t^{\rm CS}.
\end{equation}
The same linearity holds for the full directional derivative of $G$.
In particular, the variation of the transport is not frozen when the
forcing is split off.

Since transport intertwines the complex structures $J=-u\times$,
the exact difference equation along a forced flow is
\begin{equation}\label{eq:geometric-forced-difference}
 \nabla_t^u\mathcal Z_i
 =-J\sum_j(\nabla_j^u)^2\mathcal Z_i
   +\sum_{\alpha=0}^3\mathcal R_{\alpha,i}+\mathcal S_i,
\end{equation}
where repeated spatial indices are summed and
\begin{align*}
 \mathcal R_{0,i}&=-\mathcal C_t^{\rm SM}\bar\psi_i,\\
 \mathcal R_{1,i}&=-J\mathcal C_j\nabla_j^{\bar u}\bar\psi_i,\\
 \mathcal R_{2,i}&=-J\nabla_j^u(\mathcal C_j\bar\psi_i),\\
 \mathcal R_{3,i}&=-J\left\{
 R^{\mathbb S^2}(\psi_i,\psi_j)\psi_j
 -\Pi R^{\mathbb S^2}(\bar\psi_i,\bar\psi_j)\bar\psi_j\right\},\\
 \mathcal S_i&=\nabla_i^uZ-\Pi\nabla_i^{\bar u}\bar Z
                   -\mathcal C_t^{\rm CS}\bar\psi_i.
\end{align*}
Thus the source-dependent terms are exactly $\mathcal S_i$ and the
zeroth-order pairing $\int w\cdot(Z-\bar Z)\,dx$.
All four $\mathcal R_{\alpha,i}$ depend only on the spatial profiles and
on the prescribed pure velocities $\bar F,F$.

We state explicitly the estimates proved in
Appendix~\ref{app:instantaneous-geometric}.  At a fixed spatial profile
pair define
\begin{equation}\label{eq:geometric-size-functions}
 \begin{gathered}
 L=d_{\mathbb S^2}(\bar u,u),\qquad
 \mathsf a=|du|+|d\bar u|,\\
 \mathsf b=|D^2u|+|D^2\bar u|+\mathsf a^2,\qquad
 \mathsf z=\Bigl(\sum_{i=1}^2|\mathcal Z_i|^2\Bigr)^{1/2}.
 \end{gathered}
\end{equation}
Here $L$ is the short-geodesic distance, not the extrinsic distance
$|w|$.  For every $s\in[0,1]$ and $j=1,2$ the Jacobi boundary-value
calculation gives
\begin{equation}\label{eq:main-Jacobi-bound}
 |\nabla_j^\gamma\gamma_s|
 \le C\bigl(|\mathcal Z_j|+L^2\mathsf a\bigr).
\end{equation}
In particular, the differentiated transport commutator can be estimated
without any $L^2$ bound on $\nabla\mathcal Z_j$.  More precisely, for
$i=1,2$ and almost every $x\in\mathbb R^2$ the four remainders above obey
\begin{equation}\label{eq:main-geometric-remainder-bounds}
 \begin{aligned}
 |\mathcal R_{0,i}|+|\mathcal R_{1,i}|&\le CL\mathsf a\mathsf b,\\
 |\mathcal R_{2,i}|&\le C\bigl(\mathsf a^2\mathsf z
              +L\mathsf a\mathsf b+L^2\mathsf a^3\bigr),\\
 |\mathcal R_{3,i}|&\le C\mathsf a^2\mathsf z.
 \end{aligned}
\end{equation}
Here $C$ is a universal geometric constant (with the fixed number of
spatial components absorbed).  The $H^2$ bound, the uniform equivariant
comparison and H\"older's inequality then give the 
\emph{instantaneous pairing estimate}
\begin{equation}\label{eq:four-pairings-instantaneous}
 \begin{aligned}
 \sum_{i=1}^2\sum_{\alpha=0}^3
  \left|\langle\mathcal R_{\alpha,i},\mathcal Z_i\rangle\right|
 &\le C\int_{\mathbb R^2}\bigl(\mathsf a^2\mathsf z^2
          +L\mathsf a\mathsf b\mathsf z+L^2\mathsf a^3\mathsf z\bigr)\,dx\\
 &\le C_RG.
 \end{aligned}
\end{equation}
In particular, \eqref{eq:four-pairings-instantaneous} is a bound at
\emph{each fixed time}, not an integrated bound.  Appendix
\ref{app:instantaneous-geometric} proves
\eqref{eq:main-Jacobi-bound},
\eqref{eq:main-geometric-remainder-bounds} and
\eqref{eq:four-pairings-instantaneous} in this order.  The Jacobi bound
\eqref{eq:main-Jacobi-bound}, rather than $\|du\|_\infty$ alone, is what
removes the uncontrolled derivative from $\mathcal R_{2,i}$.

The principal term in \eqref{eq:geometric-forced-difference} has zero
energy pairing.  The zeroth-order pure contribution is
\[
 \frac12\frac d{dt}\|w\|_2^2\bigg|_{\rm SM}
 =-\sum_j\int w\cdot
      (\partial_j\bar u\times\partial_jw)\,dx,
\]
whose absolute value is at most $C_RG$, by
\eqref{eq:equivariant-gradient-infty} and $G\simeq_R\|w\|_{H^1}^2$.
Together with \eqref{eq:four-pairings-instantaneous} this proves
\eqref{eq:pure-geometric-bound}.  The last part of
Appendix~\ref{app:instantaneous-geometric} establishes
\eqref{eq:main-geometric-low-regularity} and justifies the transport
identities and energy chain rule at $H^2$ regularity.
Consequently
\[
 G'=G'_{\rm SM}+G'_{\rm CS}
\]
holds almost everywhere along forced paths.  The estimate of
$G'_{\rm CS}$ is supplied separately by
Lemma~\ref{lem:forcing-uniqueness}.  Agreement of the two initial data
is used only when the resulting energy inequality is applied to
uniqueness, not in any of the geometric estimates above.
\end{proof}

The same calculation supplies the difference estimate which will be used
for continuous dependence.  If $u,\bar u$ are two smooth solutions on the
same $H^2$ ball, use the canonical-frame comparison of Lemma~\ref{lem:canonical-frame-comparison} and put
\[
 D_T=\|Q-\bar Q\|_{L_t^\infty H_x^1}
    +\|Q-\bar Q\|_{L_t^4W_x^{1,4}}.
\]
The difference version of the Chern--Simons estimates
\eqref{eq:P-difference} and \eqref{eq:G-L2}--\eqref{eq:G-H1}, together
with Lemma~\ref{lem:local-difference-recovery}, gives
\begin{equation}\label{eq:CS-q-difference-local}
 \|e^{i(m-1)\chi}(\mathfrak G_u[\Gamma[u](\vec n\times u)]
 -\mathfrak G_{\bar u}[\Gamma[\bar u](\vec n\times\bar u)])\|_{L_t^1H_x^1}
 \le C_RT\left(D_T+\|u-\bar u\|_{C_tH^1}\right).
\end{equation}
Here the potential-dependent constants are controlled by
$\|V\|_{C^{3,1}}$.  In particular the term involving the difference of
third derivatives is bounded by
\[
 |V^{(3)}(1-u_3)-V^{(3)}(1-\bar u_3)|
 \le [V^{(3)}]_{\mathrm{Lip}}\,|u_3-\bar u_3|.
\]
This bound is used when differentiating the potential contribution to
$\mathfrak G$; continuity of $V^{(3)}$ alone would not give this
Lipschitz estimate.
Using \eqref{eq:pure-q-difference} for the pure part and H\"older in time,
Strichartz gives
\begin{equation}\label{eq:q-difference-Strichartz}
 D_T\le C_R\|Q(0)-\bar Q(0)\|_{H^1}
 +\varepsilon_R(T)D_T
 +\varepsilon_R(T)\|u-\bar u\|_{C_tH^1},
 \qquad \varepsilon_R(T)\to0.
\end{equation}
Thus, after decreasing $T_R$, the $D_T$ term is absorbed.  To make the
remaining low-order dependence explicit, use the energy $G(t)$ in
\eqref{eq:geometric-difference-energy}.  Lemma~\ref{lem:pure-geometric-difference}
and the forcing estimate of Lemma~\ref{lem:forcing-uniqueness} give
\begin{equation}\label{eq:Gronwall-local-estimate}
 |G'(t)|\le C_RG(t)
\end{equation}
for almost every time on any common comparison interval where the maps
lie in the fixed $H^2$ ball and $\|u-\bar u\|_\infty<1$.
Here and below the constants for the full flow may also depend on
$m$ and $\|V\|_{C^{3,1}}$.  Absolute continuity of $G$ is justified in
Appendix~\ref{app:instantaneous-geometric} without using uniqueness.
Gronwall's inequality, in either time direction, yields
\[
 G(t)\le e^{C_R|t-t_0|}G(t_0).
\]
The equivalence $G\simeq_R\|u-\bar u\|_{H^1}^2$ therefore gives
\[
 \|u(t)-\bar u(t)\|_{H^1}
 \le C_R e^{C_R|t-t_0|}\|u(t_0)-\bar u(t_0)\|_{H^1}.
\]
The multiplicative constant includes the norm-equivalence constants;
this is a Lipschitz estimate on the fixed $H^2$ ball, not a claim about
arbitrary $H^1$ data.

For clarity, the comparison interval does not impose an extra
hypothesis on nearby initial data.  Suppose both $H^2$ norms remain
bounded by $R$ on a common interval $[-T,T]$.  The uniform equivariant
$L^\infty$ comparison preceding Lemma~\ref{lem:pure-geometric-difference}
shows that, if $\|u(0)-\bar u(0)\|_{H^1}$ is sufficiently small in
terms of $R$ and $T$, the preceding bound implies
$\|u(t)-\bar u(t)\|_\infty\le1/2$ wherever the weaker condition
$\|u(t)-\bar u(t)\|_\infty<1$ holds.  Continuity in $L^\infty$ then
extends the comparison over the whole common interval.  Decreasing the
initial distance further also keeps the $H^1$ distance below the
thresholds in the fixed-time difference lemmas.  In particular
\[
 \sup_{|t|\le T}\|u(t)-\bar u(t)\|_{H^1}
 \le C_R e^{C_RT}\|u(0)-\bar u(0)\|_{H^1}.
\]
Equal initial data give $G=0$ by the same argument.

By Lemma~\ref{lem:forward-modified-field}, $H^2$ convergence of the
initial maps gives $H^1$ convergence of their initial canonical fields.
Combining this forward implication, \eqref{eq:q-difference-Strichartz},
the preceding $H^1$ bound, and the reverse estimate in
Lemma~\ref{lem:local-difference-recovery} yields continuous dependence
in $C_tH^2$ for smooth solutions, with constants controlled by the
$H^2$ ball.  The forward and reverse maps are thus used in different,
explicitly justified steps.  No higher Sobolev norm occurs.

\begin{lemma}[Smooth construction and persistence of regularity]
\label{lem:smooth-construction}
Let $V\in C^\infty([0,2])$ satisfy $V(0)=V'(0)=0$, and let
$u_0\in\mathcal X_m^\infty:=\bigcap_{s\ge2}\mathcal X_m^s$.
Then the CSLL equation has a unique smooth equivariant solution on a
nontrivial interval $(-T_{\mathrm{loc}},T_{\mathrm{loc}})$.  If
$(-T_*^{-},T_*^{+})$ is its maximal
smooth interval and
\[
 \sup_{0\le t<T}\|u(t)-\vec n\|_{H^2}\le B
 \qquad(0<T<T_*^{+}),
\]
then, for every integer $s\ge3$,
\begin{equation}\label{eq:high-regularity-persistence}
 \sup_{0\le t<T}\|u(t)-\vec n\|_{H^s}
 \le C_{s,B,T,V}\bigl(1+\|u_0-\vec n\|_{H^s}\bigr).
\end{equation}
Thus a finite maximal smooth time can occur only if the $H^2$ norm is
unbounded.  The analogous assertion holds backward in time.
\end{lemma}

\begin{proof}
We supply the covariant energy estimate behind persistence.  In this proof
$\nabla^u_jY=P_u\partial_jY$ denotes the pullback Levi--Civita connection;
it is not the Chern--Simons covariant derivative of the introduction.

\smallskip
\noindent\emph{The equivariant gradient bound.}
Equation~\eqref{eq:equivariant-gradient-infty}, proved in
Lemma~\ref{lem:forward-modified-field} from the nonzero angular modes,
gives $\|du\|_\infty\le C_mB$.  This is not a general two-dimensional
Sobolev embedding.  In the following source calculation write
$U=u_1+iu_2=e^{-im\chi}(v_1+iv_2)$; the third component is $u_3=v_3$.

\smallskip
\noindent\emph{High-order estimates for the source.}
On an $H^2$ ball the tangent source satisfies
\begin{equation}\label{eq:P-high-tame}
 \|\Gamma[u](\vec n\times u)\|_{H^s}
 \le C_{s,B,V}\bigl(1+\|u-\vec n\|_{H^s}\bigr),
 \qquad s\ge3.
\end{equation}
We verify explicitly that no radial singularity costs a derivative.
Put $f=1-u_3$, $g=|U|^2$, and $a=h/r^2$.
The radial vector $(h/r)e_r$ equals $\nabla\Delta^{-1}f$.  The averaging identity gives
\[
 \|D^ja\|_2\le (j+1)^{-1}\|D^jf\|_2,
 \qquad j\ge0,
\]
whereas
\[
 \|(\tfrac{h}{r}e_r)\|_\infty\le C_B,
 \qquad
 \|D^j(\tfrac{h}{r}e_r)\|_2\le C_j\|f\|_{H^{j-1}},\quad j\ge1.
\]
No $L^2$ bound for $(\tfrac{h}{r}e_r)$ itself is used.  Define the Cartesian
Jacobian
\[
 \operatorname{Jac}(U)
 =\partial_1u_1\partial_2u_2-\partial_2u_1\partial_1u_2.
\]
Since $U$ has index $-m$,
$\operatorname{Jac}(U)=-mg_r/(2r)$.  Hence the identity for $k$ is
\begin{equation}\label{eq:k-Cartesian-high}
 \Delta k=-fg-(h-m)\frac{g_r}{r}
 =-fg-(\tfrac{h}{r}e_r)\cdot\nabla g-2\operatorname{Jac}(U).
\end{equation}
For $s\ge3$, the product rule and Moser estimates, using
\eqref{eq:equivariant-gradient-infty}, give
\[
 \begin{aligned}
 \|\operatorname{Jac}(U)\|_{H^{s-2}}
 &\le C_s\|du\|_\infty\|du\|_{H^{s-2}},\\
 \|(\tfrac{h}{r}e_r)\cdot\nabla g\|_{H^{s-2}}
 &\le C_s\bigl(
 \|(\tfrac{h}{r}e_r)\|_\infty\|\nabla g\|_{H^{s-2}}
 +\|D(\tfrac{h}{r}e_r)\|_{H^{s-3}}\|\nabla g\|_\infty\bigr).
 \end{aligned}
\]
The estimate for $fg$ is the ordinary product estimate.  Since
$\|k\|_2\le C_B$, elliptic regularity applied to
\eqref{eq:k-Cartesian-high} therefore yields
\[
 \|k\|_{H^s}\le C_{s,B}\bigl(1+\|u-\vec n\|_{H^{s-1}}\bigr).
\]
Together with $h(2m-h)/r^2=2ma-|(\tfrac{h}{r}e_r)|^2$, these bounds prove
\eqref{eq:P-high-tame} by the product and composition rules.  For a fixed
smooth $V$, the constants may depend on derivatives of $V$ up to order
$s+1$.  The same calculation on bounded Sobolev sets gives the local
Lipschitz estimates needed in a smooth approximation scheme.

\smallskip
\noindent\emph{The covariant energy calculation.}
Set $\psi_j=\partial_ju$, $J=-u\times$, and let $\nabla^u$ act on all
pullback tangent indices.  We sum over all ordered spatial multi-indices
and use
\[
 \mathcal E_s(u)
 =1+\|u-\vec n\|_2^2
   +\sum_{k=0}^{s-1}\|(\nabla^u)^k\psi\|_2^2.
\]
On the present $H^2$ ball this is equivalent to
$1+\|u-\vec n\|_{H^s}^2$, with constants depending on $s,m,B$.
Indeed, $\partial_jY=\nabla^u_jY-(Y\cdot\psi_j)u$ for a tangent field
$Y$.  Repeated use of this identity converts the two kinds of derivatives
into each other.  All correction terms are products of lower derivatives
of $\psi$; interpolation and \eqref{eq:equivariant-gradient-infty} bound
them by $\varepsilon\|D^su\|_2+C_{s,B,\varepsilon}$, or by the
corresponding covariant top derivative, so they may be absorbed.

For clarity, the interpolation used for the curvature terms is
\[
 \|(\nabla^u)^a\psi\|_{L^{2k/a}}
 \le C_k\|\psi\|_\infty^{1-a/k}
             \|(\nabla^u)^k\psi\|_2^{a/k},
 \qquad 1\le a\le k,
\]
with $a=0$ interpreted as the $L^\infty$ norm.  This is the covariant
Gagliardo--Nirenberg inequality, obtained by integration by parts for a
metric connection and H\"older's inequality.  Thus, if $a+b+c=k$,
\begin{equation}\label{eq:covariant-cubic-interpolation}
 \|\,| (\nabla^u)^a\psi|\,
          |(\nabla^u)^b\psi|\,
          |(\nabla^u)^c\psi|\,\|_2
 \le C_k\|\psi\|_\infty^2\|(\nabla^u)^k\psi\|_2.
\end{equation}
For $k=0$ the assertion is just the $L^\infty\cdot L^\infty\cdot L^2$
product bound.

The connection identities are
\[
 \begin{aligned}
 \nabla_t^u\psi_i&=-J\nabla_i^u\nabla_j^u\psi_j
                  +\nabla_i^u[\Gamma[u](\vec n\times u)],\\
 [\nabla_i^u,\nabla_j^u]Y
   &=R^{\mathbb S^2}(\psi_i,\psi_j)Y,\\
 [\nabla_t^u,\nabla_i^u]Y
   &=R^{\mathbb S^2}(u_t,\psi_i)Y.
 \end{aligned}
\]
Both $J$ and the curvature tensor of the unit sphere are parallel.
Commuting $k$ spatial derivatives gives
\[
 \nabla_t^u((\nabla^u)^k\psi)
 =-J\nabla_j^u\nabla_j^u((\nabla^u)^k\psi)
   +\mathcal R_k+\mathcal S_k.
\]
In schematic tensor notation (with contractions of tensor and spatial
indices denoted by $\ast$), the two remainders have the form
\begin{equation}\label{eq:covariant-remainder-schematic}
 \begin{aligned}
 \mathcal R_k&=\sum_{a+b+c=k}C_{abc}
   (\nabla^u)^a\psi\ast(\nabla^u)^b\psi\ast(\nabla^u)^c\psi,\\
 \mathcal S_k&=(\nabla^u)^{k+1}[\Gamma[u](\vec n\times u)]
   +\sum_{a+b+c=k-1}C'_{abc}
    (\nabla^u)^a[\Gamma[u](\vec n\times u)]\ast(\nabla^u)^b\psi
                         \ast(\nabla^u)^c\psi.
 \end{aligned}
\end{equation}
The sums are over nonnegative integers; the second sum is empty for
$k=0$.  The coefficients $C_{abc},C'_{abc}$ denote finite
combinations of metric, complex-structure, and curvature contractions,
with constants depending on $k$, not on $u$.
These descriptions follow directly from the two commutators above: the
pure part of $u_t$ contributes one derivative of $\psi$, while its source
part contributes $\Gamma[u](\vec n\times u)$.

The principal term has zero real energy pairing, since
\[
 \int\langle Y,-J\nabla_j^u\nabla_j^u Y\rangle
 =\int\langle\nabla_j^u Y,J\nabla_j^u Y\rangle=0.
\]
Equation \eqref{eq:covariant-cubic-interpolation} bounds the pure
remainder by $C_s\|du\|_\infty^2\mathcal E_s$.  The source product
estimates and \eqref{eq:P-high-tame} give
$\|\mathcal S_k\|_2\le C_{s,B,V}\mathcal E_s^{1/2}$ for $k\le s-1$.
Mass conservation controls the zeroth-order term.  Hence
\begin{equation}\label{eq:covariant-persistence-energy}
 \frac d{dt}\mathcal E_s(u(t))\le C_{s,B,V}\mathcal E_s(u(t)).
\end{equation}
Gronwall and the norm equivalence prove
\eqref{eq:high-regularity-persistence} for smooth solutions on their
interval of existence.

\smallskip
\noindent\emph{Construction before the $H^2$ local theorem.}
To avoid using that theorem in the construction of its own approximants,
one may use the following viscous regularization.  This is an alternative
to the intrinsic wave-map approximation in \cite{M}:
\[
 u_t^\varepsilon=u^\varepsilon\times\Delta u^\varepsilon
 +\varepsilon\bigl(\Delta u^\varepsilon+|du^\varepsilon|^2u^\varepsilon\bigr)
 +\Gamma[u^\varepsilon](\vec n\times u^\varepsilon),\qquad 0<\varepsilon\le1.
\]
For the off-sphere construction use the projection $\pi(z)=z/|z|$ on
$|z|>1/2$ and extend the source by
$z\mapsto\Gamma[\pi(z)](\vec n\times z)$.  This is a locally smooth
extension in the affine high-Sobolev space, is orthogonal to $z$, and
agrees with $\Gamma[u](\vec n\times u)$ on the sphere.  The coefficients are computed
from the sphere-valued projection, so the mass integral defining $h$ is
still finite.  Extend $V$ smoothly off $[0,2]$ if needed.
For fixed $\varepsilon$ the principal symmetric part is $\varepsilon\Delta$,
so the equation is a smooth quasilinear parabolic problem with a locally
smooth lower-order source.  The scalar function
$z=|u^\varepsilon|^2-1$ solves
\[
 z_t=\varepsilon\Delta z+2\varepsilon|du^\varepsilon|^2z,
 \qquad z(0)=0;
\]
thus the sphere constraint is preserved.  Uniqueness at this smooth
parabolic level also preserves equivariance.

The same covariant calculation has the additional nonnegative
left-hand term
\[
 2\varepsilon\sum_{k=0}^{s-1}\|\nabla^u((\nabla^u)^k\psi)\|_2^2.
\]
The added curvature terms have the same three-factor structure and
coefficients bounded for $\varepsilon\le1$.  The zeroth-order identity is
\[
 \frac12\frac d{dt}\|u^\varepsilon-\vec n\|_2^2
 =-\varepsilon\int u_3^\varepsilon|du^\varepsilon|^2,
\]
which is bounded in absolute value by $\varepsilon\|du^\varepsilon\|_2^2$.
Start with a fixed sufficiently high Sobolev index, for instance $s_0=4$.
Under a bootstrap on that norm, these estimates give a time depending on
the initial $H^{s_0}$ norm and the fixed smooth potential, but not on
$\varepsilon$.  The estimate at every higher index is then linear, with
coefficients controlled by this base bound.  All Sobolev norms of the
smooth initial datum consequently persist on this common short interval.

Local compactness and a diagonal extraction give a smooth limit solving
the unregularized equation.  The gauge limit is also local: $h$ is an
integral from zero, and the omitted tail of $k$ is bounded uniformly by
$C_B R^{-2}$ on $s\ge R$, using $|h-m|\le C_B$ and
$\left(\|u_1\|_2^2+\|u_2\|_2^2\right)^{1/2}\le B$.  Thus it passes to the limit before
$R\to\infty$.  Bounds two Sobolev orders higher also bound the time
derivative in each fixed $H^s$ and give the initial value and strong time
continuity there.  Smooth uniqueness follows from the geometric
$H^1$ difference energy (all coefficients are bounded at this regularity)
and the local Lipschitz source bound.  Applying the same construction to
the sign-reversed vector field gives the backward solution.
Finally, \eqref{eq:covariant-persistence-energy} allows restarting near
any finite time at which the $H^2$ norm stays bounded; the smooth local
time can be chosen from a fixed high-order ball.  This proves the stated
persistence and continuation assertion without invoking
Theorem~\ref{thm:lwp}.
\end{proof}

\begin{prop}[Uniform $H^2$ lifespan for smooth approximants]
\label{prop:smooth-H2-lifespan}
For every $R<\infty$ there is $T_R>0$, depending only on $m$,
$R$, and $\|V\|_{C^{3,1}([0,2])}$, such that every smooth initial map
with $\|u_0-\vec n\|_{H^2}\le R$ gives a smooth solution throughout
$[-T_R,T_R]$, and
\begin{equation}\label{eq:uniform-H2-lifespan}
 \sup_{|t|\le T_R}\|u(t)-\vec n\|_{H^2}\le B_R/2.
\end{equation}
Here $B_R$ is the constant fixed before Lemma~\ref{lem:q-local}.
\end{prop}
\begin{proof}
Use the maximal smooth interval from
Lemma~\ref{lem:smooth-construction}.  On its compact subintervals impose
\[
 \sup_{|t|\le T}\|u(t)-\vec n\|_{H^2}\le B_R,\qquad
 \sup_{|t|\le T}\|u(t)-u(0)\|_{H^1}\le\rho_R.
\]
These bounds initially hold with room to spare.  The equation and the
one-map source bounds give $\|u_t\|_2\le C(B_R)$, so interpolation yields
\[
 \|u(t)-u(0)\|_2\le C(B_R)|t|,\qquad
 \|u(t)-u(0)\|_{H^1}\le C(B_R)|t|^{1/2}.
\]
Decrease the previously chosen $T_R$ so that
$C(B_R)T_R^{1/2}\le\rho_R/2$.  This improves the second bootstrap.
Smooth mass conservation gives $\|u(t)-\vec n\|_2\le R$.
Lemma~\ref{lem:q-local} gives $X_T\le K_R$, and recovery, applied with
$u_*=u(0)$, gives
\[
 \sup_{|t|\le T}\|u(t)-\vec n\|_{H^2}
 \le L_R(1+R+K_R)\le B_R/2.
\]
Thus both bootstraps improve.  Continuity proves the bounds throughout
the intersection with $[-T_R,T_R]$.  If a smooth maximal time occurred
before $T_R$, the resulting $H^2$ bound and
Lemma~\ref{lem:smooth-construction} would keep each higher Sobolev norm
finite up to that time.  Restarting in a fixed sufficiently high
Sobolev ball would contradict maximality.  The backward argument is
identical for the sign-reversed vector field.
\end{proof}

\subsection{Approximation, uniqueness, and continuous dependence}
The passage from smooth data to $H^2$ data requires a genuine tightness
argument at spatial infinity.  A family $f_j\in L^2(\mathbb R^2)$ is
\emph{uniformly $L^2$-tight} if
\begin{equation}\label{eq:uniform-L2-tightness-definition}
 \lim_{R\to\infty}\sup_j\int_{|x|>R}|f_j(x)|^2\,dx=0.
\end{equation}
Conservation of the \emph{total} mass alone does not provide such a
statement.  What is used below is the local mass flux identity together
with the uniform $H^2$ bound.

\begin{lemma}[Finite-time $H^1$ tightness]\label{lem:finite-time-tightness}
Let $u$ be a smooth equivariant solution of \eqref{eq:forced-SM-main} on
$[-T,T]$ and suppose
\[
 \sup_{|t|\le T}\|u(t)-\vec n\|_{H^2}\le B.
\]
There is a constant $C_B$ such that, for every $R\ge1$,
\begin{equation}\label{eq:local-mass-tail}
 \sup_{|t|\le T}\int_{|x|\ge2R}(1-u_3(x,t))\,dx
 \le
 \int_{|x|\ge R}(1-u_3(x,0))\,dx+\frac{C_BT}{R}.
\end{equation}
Consequently, if a family $u^{(j)}$ satisfies the same bound with one
constant $B$ and if $u^{(j)}(0)-\vec n$ is uniformly tight in $L^2$, then
\begin{equation}\label{eq:H1-tail-tightness}
 \lim_{R\to\infty}\sup_j\sup_{|t|\le T}
 \|u^{(j)}(t)-\vec n\|_{H^1(\{|x|\ge R\})}=0.
\end{equation}
\end{lemma}

\begin{proof}
Put $\rho=1-u_3$.  The Chern--Simons forcing
$\Gamma[u](\vec n\times u)$ has vanishing third component.  Moreover
\[
(u\times\Delta u)_3
 =\partial_j(u\times\partial_j u)_3.
\]
Hence every smooth solution satisfies the local conservation law
\begin{equation}\label{eq:local-mass-flux}
 \partial_t\rho+\partial_j\mathcal J_j=0,
 \qquad
 \mathcal J_j=(u\times\partial_j u)_3.
\end{equation}
Let $\chi_R$ be radial, $0\le\chi_R\le1$, with $\chi_R=0$ on
$\{|x|\le R\}$, $\chi_R=1$ on $\{|x|\ge2R\}$, and
$|\nabla\chi_R|\le C/R$.  Since
\[
 |(u\times\partial_j u)_3|
 \le \sqrt{u_1^2+u_2^2}\,|\partial_j u|
 \le |u-\vec n|\,|\partial_j u|,
\]
we obtain from \eqref{eq:local-mass-flux}
\[
 \left|\frac d{dt}\int_{\mathbb R^2}\chi_R\rho\,dx\right|
 \le \frac{C}{R}\|u-\vec n\|_2\|\nabla u\|_2
 \le \frac{C_B}{R}.
\]
Integrating in time gives \eqref{eq:local-mass-tail}.  Notice that this is
the point at which the local flux identity, rather than conservation of the
total mass, is essential.

We next upgrade the $L^2$ tail to an $H^1$ tail.  Write $f=u-\vec n$ and
choose a radial cutoff $\eta_R$ which vanishes on $\{|x|\le2R\}$, equals
one on $\{|x|\ge4R\}$, and satisfies $|\nabla\eta_R|\le C/R$.  Integration
by parts gives
\[
 \begin{aligned}
 \int \eta_R^2|\nabla f|^2
 &=-\int \eta_R^2 f\cdot\Delta f
   -2\int \eta_R f\,\nabla\eta_R\cdot\nabla f\\
 &\le \frac12\int \eta_R^2|\nabla f|^2
   +C\|f\|_{L^2(|x|\ge2R)}\|D^2f\|_2
   +\frac{C}{R^2}\|f\|_{L^2(|x|\ge2R)}^2.
 \end{aligned}
\]
Therefore
\begin{equation}\label{eq:localized-interpolation-tail}
 \|\nabla f\|_{L^2(|x|\ge4R)}^2
 \le C_B\|f\|_{L^2(|x|\ge2R)}
 +\frac{C}{R^2}\|f\|_{L^2(|x|\ge2R)}^2.
\end{equation}
Since $|f|^2=2\rho$, equations \eqref{eq:local-mass-tail} and
\eqref{eq:localized-interpolation-tail} prove
\eqref{eq:H1-tail-tightness}.
\end{proof}

Take smooth equivariant $u_0^{(j)}\to u_0$ in $H^2$.  Extend $V$
as a $C^{3,1}$ function to a neighborhood of $[0,2]$, mollify, and
subtract the value and first derivative at zero.  This gives smooth
potentials satisfying
\[
 V_j\to V\text{ in }C^3([0,2]),\qquad V_j(0)=V_j'(0)=0,
 \qquad \sup_j\|V_j\|_{C^{3,1}([0,2])}<\infty.
\]
Convergence in the $C^{3,1}$ norm is not required.  For each fixed $j$,
Lemma~\ref{lem:smooth-construction} gives a maximal smooth solution of the
unregularized CSLL equation with potential $V_j$.  The one-map $H^2$
estimates use the uniform $C^3$ bound; differentiated difference estimates
also use the uniform Lipschitz bound for $V_j^{(3)}$.  Higher derivatives
of $V_j$ enter only the persistence constant for that fixed smooth
solution.  Proposition~\ref{prop:smooth-H2-lifespan}
therefore extends every approximant to the common interval
$[-T_R,T_R]$ and gives the uniform bounds
\[
 u^{(j)}-\vec n\stackrel{*}{\rightharpoonup}u-\vec n
 \quad\hbox{in }L^\infty_tH^2_x,
 \qquad
 u_t^{(j)}\stackrel{*}{\rightharpoonup}u_t
 \quad\hbox{in }L^\infty_tL^2_x.
\]

We now justify the strong compactness used in the nonlinear terms.  On
every fixed ball $B_A$, the uniform bounds imply
\[
 \|u^{(j)}(t)-u^{(j)}(s)\|_{H^1(B_A)}
 \le C_R
 \|u^{(j)}(t)-u^{(j)}(s)\|_{L^2}^{1/2}
 \le C_R|t-s|^{1/2}.
\]
Here we used interpolation with the uniform $H^2$ bound and
$\|u_t^{(j)}\|_{L_t^\infty L_x^2}\le C_R$.  Since
$H^2(B_A)\Subset H^1(B_A)$, the Arzel\`a--Ascoli theorem gives, after a
diagonal extraction,
\begin{equation}\label{eq:local-strong-H1}
 u^{(j)}\longrightarrow u
 \qquad\hbox{strongly in }C([-T_R,T_R];H^1(B_A))
 \quad\hbox{for every }A<\infty.
\end{equation}
The initial convergence in $L^2$ implies uniform $L^2$ tightness of
$u_0^{(j)}-\vec n$.  Lemma~\ref{lem:finite-time-tightness} therefore
upgrades \eqref{eq:local-strong-H1} to the global convergence
\begin{equation}\label{eq:global-strong-H1}
 u^{(j)}\longrightarrow u
 \qquad\hbox{strongly in }C([-T_R,T_R];H^1(\mathbb R^2)).
\end{equation}

\begin{lemma}[Identification of the limiting modified field]
\label{lem:identify-Q-limit}
Let $u$ be the limit map in \eqref{eq:global-strong-H1}, and let
$e^{(j)}$ and $Q^{(j)}=e^{i(m-1)\chi}q^{(j)}$ be the canonical frames
and modified fields of $u^{(j)}$.  Then $u$ has a canonical radial frame
$e$, and, after passing to a subsequence, with
$Q(r,\chi,t)=e^{i(m-1)\chi}q[u](r,t)$, one has
\begin{align}
 e^{(j)}(\cdot,t)&\longrightarrow e(\cdot,t)
   &&\hbox{uniformly on every compact annulus},
                                                        \label{eq:frame-limit-annulus}\\
 Q^{(j)}&\longrightarrow Q
   &&\hbox{strongly in }C_tL_x^2,
                                                        \label{eq:Q-strong-L2-identification}\\
 Q^{(j)}&\stackrel{*}{\rightharpoonup}Q
   &&\hbox{in }L_t^\infty H_x^1,
                                                        \label{eq:Q-weak-H1-identification}
\end{align}
and weakly in $L_t^4W_x^{1,4}$.
\end{lemma}

\begin{proof}
The global $C_tH_x^1$ convergence gives convergence of the radial
profiles in the one-dimensional $H^1$ topology on every compact annulus.
The coefficients of the extrinsic parallel-transport ODE therefore
converge in $L^1$ there.  At infinity, the north-pole construction used in
Lemma~\ref{lem:canonical-frame-comparison} depends continuously on the
$H^1$ tails; Lemma~\ref{lem:finite-time-tightness} makes those tails
uniform in time.  Stability of the ODE then gives
\eqref{eq:frame-limit-annulus}, uniformly in $t$.  Since
\[
 W^{(j)}=u_r^{(j)}+r^{-1}J_{u^{(j)}}u_\chi^{(j)}
 \longrightarrow
 W=u_r+r^{-1}J_uu_\chi
\]
in $L^2$ on each annulus, one obtains
$q^{(j)}=[W^{(j)}]_{e^{(j)}}\to[W]_e=q[u]$ there.

The two ends are uniform.  At the origin, the $H^2$ bound and the
Gagliardo--Nirenberg inequality give
\[
 \sup_j\int_{|x|<a}|W^{(j)}|^2\,dx
 \le C\sup_j\int_{|x|<a}|\nabla u^{(j)}|^2\,dx
 \le C_Ra,
\]
and the same estimate holds for $u$.  At infinity,
$|W|\lesssim|u_r|+\sqrt{u_1^2+u_2^2}/r$, so the uniform $H^1$ tightness of
Lemma~\ref{lem:finite-time-tightness} makes the $L^2$ tails uniformly
small.  Annular convergence plus these two tail estimates proves
\eqref{eq:Q-strong-L2-identification}, uniformly in time.  Finally the
uniform bounds from Lemma~\ref{lem:q-local} provide weak compactness in
$L_t^\infty H_x^1\cap L_t^4W_x^{1,4}$; the strong $L^2$ identification
forces every such weak limit to equal $Q$, which proves
\eqref{eq:Q-weak-H1-identification}.  Interpolation gives strong
convergence in every spacetime $L^p$ with $2\le p<4$.  Together with the
strong convergence of the geometric coefficients and the boundedness of
the Hardy operator, this is sufficient to pass all polynomial and nonlocal
terms in the modified equation distributionally; weak lower
semicontinuity retains the stated $L_t^{4/3}W_x^{1,4/3}+L_t^1H_x^1$
forcing bounds.
\end{proof}

By interpolation with the uniform $H^2$ bound, the convergence is also
strong in $C_tH_x^{1+\sigma}$ for every $0<\sigma<1$.  We spell out the
passage through the nonlocal gauge coefficients, since this is where the
backward Hardy estimate is used.  Put
\[
 U_j=u^{(j)}-\vec n,\qquad U=u-\vec n,\qquad
 f_j=1-u_3^{(j)}=\frac12|U_j|^2,\qquad
 f=1-u_3=\frac12|U|^2.
\]
For every finite $1\le p<\infty$, the two-dimensional Sobolev embedding
and \eqref{eq:global-strong-H1} give
\begin{equation}\label{eq:fj-strong-Lp}
 \|f_j-f\|_{C_tL_x^p}
 \le C_{R,p}
 \|U_j-U\|_{C_tL_x^{2p}}
 \longrightarrow0.
\end{equation}
Here and in what follows the time interval is $[-T_R,T_R]$.  Since
\[
 h_j(r)-h(r)=\int_0^r s\{f_j(s)-f(s)\}\,ds,
\]
we first obtain
\begin{equation}\label{eq:h-uniform-convergence}
 \|h_j-h\|_{C_tL_r^\infty}
 \le \frac1{2\pi}\|f_j-f\|_{C_tL_x^1}
 \longrightarrow0.
\end{equation}
The combinations containing $h/r^2$ converge as well.  Indeed, with
$a_j=h_j/r^2$ and $a=h/r^2$, the exact averaging identity gives
\[
 a_j(r)-a(r)=\int_0^1 t\{f_j(tr)-f(tr)\}\,dt,
\]
and hence, for $1<p<\infty$,
\begin{equation}\label{eq:h-over-r2-difference}
 \|a_j-a\|_{C_tL_x^p}
 \le C_p\|f_j-f\|_{C_tL_x^p}\longrightarrow0.
\end{equation}
Moreover, if
\[
 A_j=\frac{h_j(r)}r e_\chi,
 \qquad A=\frac{h(r)}r e_\chi,
\]
then $A_j-A=\nabla^\perp\Delta^{-1}(f_j-f)$, so the
Hardy--Littlewood--Sobolev estimate and \eqref{eq:fj-strong-Lp} imply
\[
 \|A_j-A\|_{C_tL_x^4}
 \lesssim \|f_j-f\|_{C_tL_x^{4/3}}\longrightarrow0.
\]
Consequently
\begin{equation}\label{eq:h-coefficient-convergence}
 \frac{h_j(2m-h_j)}{r^2}
 =2ma_j-|A_j|^2
 \longrightarrow
 2ma-|A|^2
 =\frac{h(2m-h)}{r^2}
 \quad\hbox{in }C_tL_x^2.
\end{equation}

For the coefficient $k$, the exact source in \eqref{eq:k-def-main} is
$(h-m)(1-u_3^2)$.  The corresponding approximating sources converge
strongly in $C_tL_x^2$: since $|(u_3^{(j)})^2-u_3^2|
\le2|u^{(j)}-u|$ and $h_j$ is uniformly bounded, one has
\begin{equation}\label{eq:b-difference-convergence}
\begin{aligned}
&\|(h_j-m)(1-(u_3^{(j)})^2)-(h-m)(1-u_3^2)\|_{C_tL_x^2}\\
&\qquad\le C_R\bigl(\|h_j-h\|_{C_tL_r^\infty}
                +\|u^{(j)}-u\|_{C_tL_x^2}\bigr)\longrightarrow0.
\end{aligned}
\end{equation}
Recall the backward radial Hardy operator
$(Tf)(r)=\int_r^\infty f(s)\,ds/s$.  To verify its $L^2$ bound,
suppress the harmless factor $2\pi$, set $r=e^t$, and write
\begin{equation}\label{eq:backward-Hardy-log-coordinate}
 e^t(Tf)(e^t)
 =\int_0^\infty e^{-y}\,e^{t+y}f(e^{t+y})\,dy.
\end{equation}
The right-hand side is convolution, up to the direction of translation,
with the $L^1$ kernel $e^{-y}\mathbf1_{\{y>0\}}$.  The map
$f(r)\mapsto e^tf(e^t)$ is an isometry from $L^2((0,\infty),r\,dr)$
onto $L^2(\mathbb R,dt)$.  Young's inequality thus gives
\begin{equation}\label{eq:backward-Hardy-difference}
\begin{aligned}
 &\|Tf\|_{L^2(rdr)}\le\|f\|_{L^2(rdr)},\\
 &\|k_j-k\|_{C_tL_x^2}
 \le\|(h_j-m)(1-(u_3^{(j)})^2)
         -(h-m)(1-u_3^2)\|_{C_tL_x^2}\longrightarrow0.
\end{aligned}
\end{equation}
If $\Gamma_j$ denotes the coefficient corresponding to $V_j$, then the
potential term also converges in $C_tL^2$.  Indeed, using
$V_j'(0)=V'(0)=0$,
\[
 \|V_j'(f_j)-V'(f)\|_2
 \le C\|f_j-f\|_2
 +\|V_j''-V''\|_\infty\|f_j\|_2\longrightarrow0.
\]
Combining this with \eqref{eq:h-coefficient-convergence} and
\eqref{eq:backward-Hardy-difference} yields
$\Gamma_j\to\Gamma$ in $C_tL_x^2$.  The uniform $L^\infty$ coefficient
bounds from Lemma~\ref{lem:tame} and \eqref{eq:global-strong-H1} therefore
show that
\begin{equation}\label{eq:P-strong-limit-detailed}
 \Gamma_j[u^{(j)}](\vec n\times u^{(j)})\longrightarrow\Gamma[u](\vec n\times u)
 \qquad\hbox{strongly in }C_tL_x^2,
\end{equation}
where $\Gamma_j$ denotes the scalar coefficient in
\eqref{eq:Gamma-main} for the potential $V_j$.

Thus all nonlinear terms pass to the limit and $u$ satisfies CSLL
distributionally.  By Lemma~\ref{lem:identify-Q-limit}, the field
$Q=e^{i(m-1)\chi} q[u]$ is the canonical modified field of this limit
map.  We next make precise the time continuity recovered from it.  The
bounds used in Lemma~\ref{lem:q-local} pass to the limit and give
\[
 F_{\rm SM}:=e^{i(m-1)\chi}(\mathcal V[u,q]q)
 \in L_t^{4/3}W_x^{1,4/3},
 \qquad
 F_{\rm CS}:=ie^{i(m-1)\chi}\mathfrak G_u[\Gamma[u](\vec n\times u)]
 \in L_t^1H_x^1.
\]
Hence for $s,t\in[-T_R,T_R]$ the Duhamel formula for
\eqref{eq:exact-Q-CSLL} is
\begin{equation}\label{eq:Duhamel-Q-limit}
 Q(t)=e^{-i(t-s)\Delta}Q(s)
 -i\int_s^t e^{-i(t-\tau)\Delta}
 \{F_{\rm SM}(\tau)+F_{\rm CS}(\tau)\}\,d\tau.
\end{equation}
The dual Strichartz estimate for the admissible pair $(4,4)$ and
the unitarity of the free group
on $H^1$ give
\begin{align*}
 &\left\|Q(t)-e^{-i(t-s)\Delta}Q(s)\right\|_{H^1}\notag\\
 &\qquad\le C
 \|F_{\rm SM}\|_{L_t^{4/3}W_x^{1,4/3}([s,t])}
 +\|F_{\rm CS}\|_{L_t^1H_x^1([s,t])}
 \longrightarrow0\qquad(t\to s).
\end{align*}
Together with the strong continuity of $e^{-it\Delta}$ on $H^1$, this
proves
\begin{equation}\label{eq:Q-CtH1-limit}
 Q\in C([-T_R,T_R];H^1(\mathbb R^2)).
\end{equation}

We finally explain the recovery of \emph{time continuity} in $H^2$.
Lemma~\ref{lem:local-recovery} gives the one-map $H^2$ control, whereas
continuity is obtained from the difference estimate
Lemma~\ref{lem:local-difference-recovery}.  The latter estimate extends
from smooth maps to $H^2$ maps by density, with the same constants on a
fixed $H^2$ ball.  Let $B$ be the uniform weak-$*$ $L_t^\infty H_x^2$
bound above and let $E\subset[-T_R,T_R]$ be a full-measure set on which
$u(t)-\vec n\in H^2$ and $\|u(t)-\vec n\|_{H^2}\le B$.  The canonical frame is used at every time, so the fields
$Q(t)$ and $Q(s)$ already belong to the same fixed complex gauge.
Therefore, for $s,t\in E$ close
enough that $\|u(t)-u(s)\|_{H^1}\le\rho_B$,
\begin{equation}\label{eq:H2-time-difference-limit}
 \|u(t)-u(s)\|_{H^2}
 \le C_B\left(
 \|u(t)-u(s)\|_{H^1}
 +\|Q(t)-Q(s)\|_{H^1}\right).
\end{equation}
The right-hand side tends to zero as $t\to s$ by
\eqref{eq:global-strong-H1} and \eqref{eq:Q-CtH1-limit}.  Thus $u:E\to
H^2$ is locally uniformly continuous.  Given an arbitrary
$t_0\in[-T_R,T_R]$, choose $t_j\in E$ with $t_j\to t_0$.  Then
$\{u(t_j)\}$ is Cauchy in $H^2$ by
\eqref{eq:H2-time-difference-limit}; its $H^2$ limit agrees with $u(t_0)$
because $u(t_j)\to u(t_0)$ in $H^1$.  This defines the unique continuous
$H^2$ representative and proves
\[
 u-\vec n\in C([-T_R,T_R];H^2(\mathbb R^2)).
\]
The right-hand side of \eqref{eq:forced-SM-main} is then continuous
in $L^2$: use $H^2\hookrightarrow L^\infty$, continuity of $\Delta u$
in $L^2$, and the coefficient difference estimates.  The distributional
equation therefore also gives $u-\vec n\in C^1_tL^2_x$.

We next record explicitly how the Chern--Simons source enters the
parallel-transport uniqueness argument.

\begin{lemma}[Forcing contribution to the geometric difference energy]
\label{lem:forcing-uniqueness}
Let $u,\bar u$ be two $H^2$ solutions on a time interval $I$ with
\[
 \sup_{t\in I}\bigl(
 \|u(t)-\vec n\|_{H^2}+
 \|\bar u(t)-\vec n\|_{H^2}\bigr)\le R.
\]
On a subinterval on which $\|u-\bar u\|_\infty<1$, use
$\Pi$, $\mathcal Z_j$, and $G(t)$ as defined in
\eqref{eq:geometric-difference-energy}.
The part of $G'(t)$ generated by the forcing terms
$\Gamma[u](\vec n\times u)$ and $\Gamma[\bar u](\vec n\times\bar u)$ satisfies
\begin{equation}\label{eq:forcing-G-bound}
 |G'_{\rm CS}(t)|\le C_R G(t).
\end{equation}
\end{lemma}

\begin{proof}
The coefficient estimates used in Lemma~\ref{lem:tame}, applied to
differences, give both
\begin{align}
 \|\Gamma[u](\vec n\times u)-\Gamma[\bar u](\vec n\times\bar u)\|_2
 &\le C_R\|u-\bar u\|_{H^1},                              \label{eq:P-L2-Lipschitz}\\
 \|\nabla^u[\Gamma[u](\vec n\times u)]
       -\Pi\nabla^{\bar u}[\Gamma[\bar u](\vec n\times\bar u)]\|_2
 &\le C_R\left(
 \|u-\bar u\|_2+ \sum_j\|\mathcal Z_j\|_2\right).                                 \label{eq:P-covariant-Lipschitz}
\end{align}
Here $G\simeq_R\|u-\bar u\|_{H^1}^2$, and the uniform comparison
preceding Lemma~\ref{lem:pure-geometric-difference} gives
$\|u-\bar u\|_\infty\le C_RG^{1/2}$.
This bound uses the common equivariant pole charts (and the trivial
large-distance bound), not a general $H^1$ embedding into $L^\infty$.
Apply the averaging and Hardy operators to differences of the coefficients;
then apply the product rule.  Products with a coefficient difference use
this $L^\infty$ comparison or the $L^4\cdot L^4$ pairing, and terms with
a derivative difference use $\|du\|_\infty+\|d\bar u\|_\infty\le C_R$.
These give the ordinary $H^1$ source-difference estimate.  If the
$H^1$ distance is not small, the separate uniform source bounds and
$G\simeq_R\|u-\bar u\|_{H^1}^2$ give the same estimate after increasing
$C_R$.  Replacing
ordinary subtraction by transport adds
$O(|u-\bar u|)\,|\nabla^{\bar u}[\Gamma[\bar u](\vec n\times\bar u)]|$, whose $L^2$ norm
is at most $C_RG^{1/2}$.  The normal projections in the covariant
derivatives have the same bound.  This proves the two displayed estimates.

The zeroth-order part of the difference energy contributes
\[
 \left|\int (u-\bar u)\cdot
       \{
       \Gamma[u](\vec n\times u)-\Gamma[\bar u](\vec n\times\bar u)
       \}\,dx\right|
 \le C_RG(t).
\]
For the derivative part use the exact source term $\mathcal S_j$ in
\eqref{eq:geometric-forced-difference}.  Its direct term is
$\nabla_j^u[\Gamma[u](\vec n\times u)]
 -\Pi\nabla_j^{\bar u}[\Gamma[\bar u](\vec n\times\bar u)]$; the
transport term is exactly $-\mathcal C_t^{\rm CS}\bar\psi_j$.
Equations \eqref{eq:temporal-transport-linear} and
\eqref{eq:main-velocity-Jacobi-bound}, together with the uniform $L^\infty$
source bound, give the pointwise estimate
\[
 |\mathcal C_t^{\rm CS}\bar\psi_j|
 \le C_R|u-\bar u|\,|\partial_j\bar u|.
\]
Thus every time-dependent change of transport caused by the forcing is
included, not discarded.  Using
\eqref{eq:P-covariant-Lipschitz}, the equivariant bound \eqref{eq:equivariant-gradient-infty},
and Cauchy--Schwarz, every source contribution is therefore bounded
by
\[
 C_R\left(
 \|u-\bar u\|_2^2+
 \sum_j\|\mathcal Z_j\|_2^2\right)=2C_RG(t).
\]
This proves \eqref{eq:forcing-G-bound}.  Since equal initial data give
zero $L^\infty$ distance at $t=0$, the short-interval argument applies;
once uniqueness is obtained there, it is iterated over $I$.
\end{proof}

For uniqueness we use the parallel-transport difference energy of
\cite{ShimizuUnique}.  On an $H^2$ ball the Chern--Simons vector field
satisfies the difference estimate
\begin{equation}\label{eq:P-difference}
\|\Gamma[u](\vec n\times u)-\Gamma[\bar u](\vec n\times\bar u)\|_{H^1}
\le C_R\|u-\bar u\|_{H^1}.
\end{equation}
Indeed $h[u]-h[\bar u]$ is linear in $u_3-\bar u_3$, the averaging operator
$h/r^2$ is bounded on $L^p$ for $1<p<\infty$, and the backward Hardy operator
in $k$ is bounded on the radial weighted spaces; the product estimates are
the difference versions of Lemma~\ref{lem:tame}, using the equivariant
comparison just recorded.  Lemma~\ref{lem:pure-geometric-difference}
provides the pure contribution, while
Lemma~\ref{lem:forcing-uniqueness} contributes only $C_RG(t)$.
Therefore $|G'(t)|\le C_RG(t)$ almost everywhere on the comparison
interval.  Gronwall's inequality and the comparison-interval argument
after \eqref{eq:Gronwall-local-estimate} give uniqueness when $G(0)=0$.
For nearby data in the same bounded $H^2$ class, the same calculation
gives the local $H^1$ Lipschitz estimate recorded there.  Combining it with
\eqref{eq:q-difference-Strichartz} and
Lemma~\ref{lem:local-difference-recovery} gives the $H^2$ continuous
dependence estimate.  To pass this estimate through smooth potential
approximation, split any source difference into a map difference at the
same potential and a potential difference at the same map.  On a fixed
$H^2$ ball the latter is bounded in the differentiated source norm by
$C_R\|V_j-V_k\|_{C^3}$: the explicit source formula uses derivatives
of $V$ only through order three.  Thus the Strichartz estimate acquires
only $C_RT\|V_j-V_k\|_{C^3}$.  The lower-order energy estimate,
by Cauchy--Schwarz and Young's inequality, becomes
\[
 |G'|\le C_RG+C_R\|V_j-V_k\|_{C^3}^2,
\]
so Gronwall gives the corresponding vanishing error in the $H^1$
difference bound.  The uniform $C^{3,1}$ bounds control
the map-difference constants.  Passing to the smooth limit therefore
preserves the asserted $H^2$ continuous dependence.

Conservation of $E$ and $M$ follows from Lemma~\ref{lem:conservation} by
passage to the smooth limit.  Finally, if a maximal solution had
$T^{+}<\infty$ while $\sup_{t<T^{+}}\|u(t)-\vec n\|_{H^2}<\infty$, the uniform
local time just proved would restart the solution near $T^{+}$ and extend it
past $T^{+}$.  This is the blow-up alternative.  Reversibility gives the
backward statement and completes the proof of Theorem~\ref{thm:lwp}.

\begin{corollary}[Local spacetime regularity of the canonical field]
\label{cor:local-spacetime}
Let $u$ be the maximal solution of Theorem~\ref{thm:lwp}, and let $Q(r,\chi,t)=e^{i(m-1)\chi}q[u](r,t)$ be its
canonical planar modified field.  Then
\begin{equation}\label{eq:local-spacetime-class}
 Q\in C((-T^{-},T^{+});H^1(\mathbb R^2))
       \cap L^4_{\rm loc}((-T^{-},T^{+});W^{1,4}(\mathbb R^2)).
\end{equation}
\end{corollary}
\begin{proof}
The assertion $Q\in C_tH^1$ follows directly from $u-\vec n\in C_tH^2$
and the forward continuity in Lemma~\ref{lem:forward-modified-field}.
On the initial local interval the $L_t^4W_x^{1,4}$ assertion is the bound
in Lemma~\ref{lem:q-local} and its identification in
Lemma~\ref{lem:identify-Q-limit}.
On any compact subinterval the already proved $C_tH^2$ regularity gives
a finite $H^2$ bound.  Restart the local construction and cover the
compact subinterval by finitely many of its uniform local intervals.
Uniqueness identifies each restarted solution with $u$.  Their canonical
frames have the same normalization, so the resulting fields also agree.
This proves \eqref{eq:local-spacetime-class}.  No continuation estimate
from Section~\ref{sec:continuation} is used in this argument.
\end{proof}

\section{Constrained Hessian coercivity and orbital stability}\label{sec:stability}

\subsection{Orbital stability and linearized analysis of a ground state}
Throughout this subsection, we assume that $m\geq 2$, $V\in C^{3,1}([0,2])$ is
nonnegative, $V(0)=V'(0)=0$, and $\omega>0$.  Define
$H(u) = E(u) +\omega M(u)$.

For the scalar great-circle reduction, we first fix the comparison
class and the topology used in the variational problem.
Here $H^1_{\rm rad}(\mathbb R^2)$ denotes the radial subspace of $H^1$,
with
\begin{equation}\label{eq:H1rad-def}
 \|\theta\|_{H^1_{\rm rad}}^2
 =2\pi\int_0^\infty\bigl(|\theta_r|^2+|\theta|^2\bigr)r\,dr.
\end{equation}
The scalar comparison class for our fixed topological component is
\begin{equation}\label{eq:scalar-great-circle-class}
 \mathscr A_m^{\rm gc}=\{\theta\in H^1_{\rm rad}(\mathbb R^2):
 0\le\theta\le\pi,\quad\theta(0)=\pi, \quad \int_0^\infty\frac{\sin^2\theta(r)}r\,dr<\infty\}.
\end{equation}
The value at zero is the radial representative's limit; the limit at
infinity is zero by the radial $H^1$ property.  The larger admissible
class in \cite{JHS} allows any odd multiple of $\pi$ at the origin.
The energy-preserving truncation of \cite[Lemma~3.1]{JHS}
reduces its minimization to \eqref{eq:scalar-great-circle-class}.

\begin{lemma}[Great-circle reduction of equivariant minimizers]
\label{lem:great-circle-reduction}
Let
\[
 u(r,\chi)=R_{-m\chi}v(r)
\]
be an admissible $m$--equivariant map in the topological component
$v(0)=-\vec n$, $v(\infty)=\vec n$, with $u-\vec n\in H^1(\mathbb R^2)$
and $H(u)<\infty$.  Define
\[
 \theta(r):=\arccos v_3(r)\in[0,\pi],\qquad
 v^\sharp(r):=(\sin\theta(r),0,\cos\theta(r)),
 \qquad u^\sharp=R_{-m\chi}v^\sharp.
\]
Then $u^\sharp$ belongs to the same topological class,
\begin{equation}\label{eq:great-circle-invariants}
 h[v^\sharp]=h[v],\qquad M(u^\sharp)=M(u),
\end{equation}
and
\begin{equation}\label{eq:great-circle-energy-drop}
 H(u^\sharp)\leq H(u).
\end{equation}
If in addition $0<\theta(r)<\pi$ for every $r>0$, then the pair $(v_1,v_2)$ admits a phase lift
\[
 v=(\sin\theta\cos\varphi,\sin\theta\sin\varphi,\cos\theta),
 \qquad \varphi\in H^1_{\rm loc}(0,\infty),
\]
and
\begin{equation}\label{eq:great-circle-energy-identity}
 H(u)-H(u^\sharp)
 =\pi\int_0^\infty \sin^2\theta\,|\varphi_r|^2 r\,dr.
\end{equation}
Consequently the infimum of $H$ over the full vector-valued equivariant
class equals the infimum of the scalar functional
\begin{equation}\label{eq:JHS-scalar-functional}
\begin{aligned}
 H_\omega\bigl(R_{-m\chi}(\sin\theta,0,\cos\theta)\bigr)
 &=2\pi\int_0^\infty\Bigl\{\frac12\theta_r^2
 +\frac{(h_\theta-m)^2}{2r^2}\sin^2\theta\\
 &\hspace{7em}+V(1-\cos\theta)+\omega(1-\cos\theta)\Bigr\}r\,dr,
\end{aligned}
\end{equation}
where
\[
 h_\theta(r)=\int_0^r s(1-\cos\theta(s))\,ds.
\]
The scalar infimum is taken over $\mathscr A_m^{\rm gc}$ in
\eqref{eq:scalar-great-circle-class}.  In particular, every
vector-valued ground state is, after one constant target rotation about
$\vec n$, of the form
\[
 R_{-m\chi}(\sin\widetilde\theta,0,\cos\widetilde\theta),
\]
where $\widetilde\theta$ is a scalar ground state in the sense of
\cite{JHS}; conversely every such scalar ground state gives a ground state
in the full vector-valued equivariant class.
\end{lemma}

\begin{proof}
The function $\theta$ is the spherical distance from $v(r)$ to $\vec n$.
Since $v$ is locally absolutely continuous as a curve in $\mathbb S^2$,
so is $\theta$, and the metric derivative of the distance function gives
\begin{equation}\label{eq:theta-speed}
 |\theta_r|\leq |v_r|\qquad\hbox{for a.e. }r>0.
\end{equation}
Equivalently, on the set $0<\theta<\pi$ one may differentiate
$v_3=\cos\theta$ and use $v\cdot v_r=0$ to obtain
\[
 |(v_3)_r|^2\leq(1-v_3^2)|v_r|^2
 =\sin^2\theta\,|v_r|^2,
\]
which is exactly \eqref{eq:theta-speed}.  Moreover
\[
 |u^\sharp-\vec n|^2=2(1-\cos\theta)=2(1-v_3)=|u-\vec n|^2,
\]
and
\[
 |\partial_\chi u^\sharp|^2=m^2\sin^2\theta
 =m^2(v_1^2+v_2^2)=|\partial_\chi u|^2.
\]
Thus $u^\sharp-\vec n\in H^1$, and it has the same endpoint values as
$u$.  Moreover, because $0\leq\theta\leq\pi$,
\[
 \theta^2\leq C(1-\cos\theta)=C(1-v_3),
\]
so $\theta\in L^2(\mathbb R^2)$; together with
\eqref{eq:theta-speed} and the angular derivative estimate above, this
shows that $\theta$ belongs to the scalar admissible class of \cite{JHS}.
Since $v_3=\cos\theta$, both $h$ and the mass are unchanged, which
proves \eqref{eq:great-circle-invariants}.  The angular gauge term and the
potential term in \eqref{eq:energy-radial} depend only on $v_3$ and
$(v_1^2+v_2^2)=1-v_3^2$, so \eqref{eq:theta-speed} gives
\eqref{eq:great-circle-energy-drop}.

For the sharper identity, assume $0<\theta<\pi$ on $(0,\infty)$.
Then $v_1+iv_2$ never vanishes, and since $(0,\infty)$ is simply
connected it admits a global $H^1_{\rm loc}$ phase lift $\varphi$.  The
standard polar-coordinate identity gives
\[
 |v_r|^2=\theta_r^2+\sin^2\theta\,\varphi_r^2
 \qquad\hbox{a.e. on }(0,\infty).
\]
Substitution into \eqref{eq:energy-radial} gives
\eqref{eq:great-circle-energy-identity}.

It remains to compare the variational classes.  The preceding construction
maps every vector-valued competitor to a scalar one.  Conversely, any
scalar competitor $\theta$ in the admissible class of \cite{JHS}, after the
energy-preserving truncation of \cite[Lemma~3.1]{JHS}, belongs to
$\mathscr A_m^{\rm gc}$ in \eqref{eq:scalar-great-circle-class}; then
$R_{-m\chi}(\sin\theta,0,\cos\theta)$ is an admissible vector-valued
competitor with exactly the energy \eqref{eq:JHS-scalar-functional}.
Hence the two infima coincide.

Finally, let $u$ be a vector-valued ground state.  Then $u^\sharp$ is a
scalar ground state.  By \cite[Theorem~1.2]{JHS}, its profile satisfies
$0<\theta(r)<\pi$ for every $r>0$.  Equality of the two ground-state
energies in \eqref{eq:great-circle-energy-identity} therefore gives
$\varphi_r=0$ a.e. on $(0,\infty)$, so $\varphi$ is constant.  Thus $u$ is
a constant target rotation of $u^\sharp$.  The converse follows from the
equality of the two infima.
\end{proof}

By Lemma~\ref{lem:great-circle-reduction}, after a fixed target rotation we
may and shall write any ground state as
\[
\widetilde u = R_{-m\chi}\widetilde v(r)
= R_{-m\chi}(\sin\widetilde\theta(r), 0, \cos\widetilde\theta(r)).
\]
In particular $\widetilde\theta$ is a scalar ground state in the sense of
\cite{JHS}; hence $0<\widetilde\theta<\pi$ on $(0,\infty)$ and, because
$\omega>0$, it decays exponentially by \cite[Theorem~1.2]{JHS}.
We note that $\widetilde v$ satisfies the elliptic equation
\begin{equation}\label{EL-v}
\widetilde v \times \left(\Delta_r\widetilde v+(k[\widetilde v]+\frac{(h[\widetilde v]-m)^2}{r^2}\widetilde v_3+V'(1-\widetilde v_3)+\omega)\vec{n} \right) = 0.
\end{equation}
It follows that the standing wave $R_{\omega t}\widetilde u$ is a solution of \eqref{eqn-single-m}.

In this section, we prove the orbital stability of the ground state
$\widetilde u$ under $m$-equivariant perturbations, subject to the
nondegeneracy condition stated below.
\begin{thm}[Conditional orbital stability]\label{thm:conditional-stability}
Assume that the ground state $\widetilde u$ satisfies the constrained
nondegeneracy condition \eqref{theta-ND} below.  Then, for any $\varepsilon > 0$,
there exists $\delta > 0$ such that if $u_0$ is a sphere-valued
$m$--equivariant map in the same topological component as $\widetilde u$,
$u_0-\widetilde u\in H^1(\R^2)$, and
$\|u_0-\widetilde u\|_{H^1} < \delta$, then
\[
\sup_{t \in [0, T]}d_{H^1}(u(t),\mathcal O(\widetilde u)) < \varepsilon,
\]
where $u-\widetilde u\in C([0,T];H^1(\R^2))$ is an $m$--equivariant
sphere-valued solution of \eqref{eqn-single-m} with initial condition
$u_0$ for which the energy and mass are conserved on $[0,T]$.
\end{thm}

We distinguish radial coefficients from spatial tangent sections.  Let
$\xi(r)$ be a smooth compactly supported radial vector field satisfying
$\xi(r)\cdot\widetilde v(r)=0$; its corresponding spatial tangent section is
$R_{-m\chi}\xi(r)$.  Choose a twice
differentiable sphere-valued curve $u_s=R_{-m\chi}v_s$ with
$v_0=\widetilde v$, $\partial_sv_s|_0=\xi$.  Write
$\upsilon=\partial_s^2v_s|_0$ and define, before differentiating,
\[
 \widetilde h=h[\widetilde v],\qquad
 \widetilde k=k[\widetilde v],\qquad
 \delta h[Y](r)=-\int_0^r \sigma Y_3(\sigma)\,d\sigma.
\]
Since $h$ is affine in the third component,
$\partial_sh[v_s]|_0=\delta h[\xi]$ and
$\partial_s^2h[v_s]|_0=\delta h[\upsilon]$; the second expression is a
curve derivative, not a nonzero second Frechet derivative of $h$.
The extrinsic calculation gives
\[
\begin{aligned}
\frac{d^2}{ds^2}\bigg|_{s=0}H(u_s)
={}&\frac{d^2}{ds^2}\bigg|_{s=0}
\int_0^\infty\Bigl[\frac12|\partial_rv_s|^2
 +\frac{(h[v_s]-m)^2}{2r^2}\bigl(1-(v_s)_3^2\bigr)\\
&\hspace{38mm}+V\bigl(1-(v_s)_3\bigr)
 +\omega\bigl(1-(v_s)_3\bigr)\Bigr]2\pi r\,dr\\
={}&\int_0^\infty\Bigl[|\partial_r\xi|^2
 +\frac{\delta h[\xi]^2(1-\widetilde v_3^2)}{r^2}
 -\frac{4(\widetilde h-m)\widetilde v_3\delta h[\xi]\xi_3}{r^2}\\
&\hspace{21mm}-\frac{(\widetilde h-m)^2\xi_3^2}{r^2}
 +V''(1-\widetilde v_3)\xi_3^2\Bigr]2\pi r\,dr+(\star),
\end{aligned}
\]
where
\[
\begin{aligned}
(\star)=\int_0^\infty\Bigl[&\langle\partial_r\widetilde v,\partial_r\upsilon\rangle
 +\frac{\widetilde h-m}{r^2}(1-\widetilde v_3^2)\delta h[\upsilon]
 -\frac{(\widetilde h-m)^2}{r^2}\widetilde v_3\upsilon_3\\
&-V'(1-\widetilde v_3)\upsilon_3-\omega\upsilon_3\Bigr]2\pi r\,dr.
\end{aligned}
\]
Applying Fubini's theorem and integration by parts to $(\star)$, we obtain
\[
(\star) = \int_0^\infty \left\langle -\Delta_r\widetilde v-\left(\widetilde k+\frac{(\widetilde h-m)^2}{r^2}\widetilde v_3+V'(1-\widetilde v_3)+\omega\right)\vec{n}, \upsilon \right\rangle \, 2\pi rdr. 
\]
The equation \eqref{EL-v} shows
\[
\Delta_r\widetilde v+\left(\widetilde k+\frac{(\widetilde h-m)^2}{r^2}\widetilde v_3+V'(1-\widetilde v_3)+\omega\right)\vec{n} = \widetilde\lambda \widetilde v,
\]
where
\[
\widetilde\lambda = \langle\Delta_r\widetilde v, \widetilde v\rangle +(\widetilde k+\frac{(\widetilde h-m)^2}{r^2}\widetilde v_3+V'(1-\widetilde v_3)+\omega)\widetilde v_3.
\]
Differentiating the equation $\langle v_s, v_s \rangle = 1$ with respect to $s$ twice, we get
\[
\langle \upsilon, \widetilde v \rangle = -|\xi|^2.
\]
Then we see that $(\star) = \int_0^\infty -\widetilde\lambda \langle\widetilde v, \upsilon\rangle\,2\pi rdr  
=  \int_0^\infty \widetilde\lambda|\xi|^2 \,2\pi rdr $.

The radial and spatial bases are, respectively,
\[
 \begin{aligned}
 e_\theta^0(r)&=(-\cos\widetilde\theta,0,\sin\widetilde\theta),
 &e_\phi^0&=(0,1,0),\\
 e_\theta(x)&=R_{-m\chi}e_\theta^0(r),
 &e_\phi(x)&=R_{-m\chi}e_\phi^0.
 \end{aligned}
\]
Thus $(e_\theta(x),e_\phi(x))$ is a basis of
$T_{\widetilde u(x)}\mathbb S^2$, and
\[
 \xi=\xi_\theta e_\theta^0+\xi_\phi e_\phi^0.
\]
For any scalar angle $\theta$, define the integral operator
\[
 (\mathcal I_\theta f)(r)
 :=\int_0^r \sigma\sin\theta(\sigma)f(\sigma)\,d\sigma.
\]
The sign distinction between an angle variation and a polar coefficient
is then explicit:
\[
 D_\theta h_\theta[g]=\mathcal I_\theta g,\qquad
 \delta h[f e_\theta^0]=-\mathcal I_{\widetilde\theta}f,
\]
because $e_\theta^0=-\partial_\theta(\sin\theta,0,\cos\theta)$.

For real radial functions write
\begin{equation}\label{eq:radial-L2-inner-product}
 L^2_{\rm rad}:=L^2((0,\infty),2\pi r\,dr),\qquad
 \langle f,g\rangle_{L^2_{\rm rad}}
 =2\pi\int_0^\infty f(r)g(r)\,r\,dr.
\end{equation}
The common form domain $\mathcal H$ is the completion of
$C_c^\infty(0,\infty;\mathbb R)$ for the norm
\[
 \|f\|_{\mathcal H}^2
 =\int_0^\infty
 \left(|f'|^2+\left(\frac{m^2}{r^2}+1\right)|f|^2\right)2\pi r\,dr.
\]
All spectral statements can equivalently be read on the complexification.
First define the azimuthal form on test functions by
\[
 \mathfrak q_\phi[f,g]
 =\int_0^\infty(f'g'+\widetilde\lambda fg)2\pi r\,dr,
 \qquad \mathfrak q_\phi[f]=\mathfrak q_\phi[f,f].
\]
Here the previously defined multiplier is
\[
 \widetilde\lambda={}-(\widetilde\theta')^2+
 \widetilde k\cos\widetilde\theta
 +\frac{(\widetilde h-m)^2}{r^2}\cos^2\widetilde\theta
 +\{V'(1-\cos\widetilde\theta)+\omega\}\cos\widetilde\theta.
\]
The difference $\widetilde\lambda-m^2/r^2-1$ is bounded, by the
origin regularity and the decay of the fixed profile.  Hence this form
extends to a closed lower-bounded form on $\mathcal H$.
The associated self-adjoint operator is defined by
\begin{equation}\label{eq:Lphi-operator-domain}
 \begin{aligned}
 L_\phi f&=g, \\
 D(L_\phi)&=\bigl\{f\in\mathcal H:\ \exists\,g\in L^2_{\rm rad}
 \text{ such that }\mathfrak q_\phi[f,\varphi]
 =\langle g,\varphi\rangle_{L^2_{\rm rad}}
 \ \forall\varphi\in\mathcal H\bigr\} 
 \end{aligned}
\end{equation}
On test functions its differential expression is
\begin{equation}\label{eq:Lphi-differential-expression}
 L_\phi f=-f''-\frac1r f'+\widetilde\lambda f.
\end{equation}
For a general $f\in\mathcal H$ we use $\mathfrak q_\phi[f]$, not the
possibly undefined operator pairing $\langle L_\phi f,f\rangle$.

For the polar component set
\[
 \begin{split}
 \mathcal U_\theta={}&\widetilde k\cos\widetilde\theta
 +\frac{(\widetilde h-m)^2}{r^2}\cos(2\widetilde\theta)
 +V''(1-\cos\widetilde\theta)\sin^2\widetilde\theta\\
 &+\{V'(1-\cos\widetilde\theta)+\omega\}\cos\widetilde\theta,
 \end{split}
\]
and use the following initially formal operator on test functions:
\[
 \begin{split}
 L_\theta f={}&-\Delta_rf+\mathcal U_\theta f
 +\sin\widetilde\theta(r)\int_r^\infty
   \frac{\sin^2\widetilde\theta(s)}s
       (\mathcal I_{\widetilde\theta}f)(s)\,ds\\
 &+\frac{4\sin\widetilde\theta\cos\widetilde\theta
              (\widetilde h-m)}{r^2}\mathcal I_{\widetilde\theta}f.
 \end{split}
\]
Define its symmetric form by
\begin{equation}\label{eq:polar-Hessian-form}
 \mathfrak q_\theta[f,g]
 =\tfrac12\{\langle L_\theta f,g\rangle_{L^2}
                +\langle f,L_\theta g\rangle_{L^2}\},
 \qquad \mathfrak q_\theta[f]=\mathfrak q_\theta[f,f].
\end{equation}
Origin regularity and exponential decay give
\begin{equation}\label{eq:polar-form-bounded-perturbation}
 \mathfrak q_\theta[f,g]
 =\int_0^\infty
 \left(f'g'+\left(\frac{m^2}{r^2}+1\right)fg\right)2\pi r\,dr
 +\langle B_\theta f,g\rangle_{L^2},
\end{equation}
The perturbation can be written explicitly.  Put
$S(r)=\sin\widetilde\theta(r)$ and $a(r)=\widetilde h(r)-m$ and
introduce the two nonlocal maps
\begin{equation}\label{eq:polar-nonlocal-operators}
 \begin{aligned}
 (\mathcal N_1f)(r)
 &=S(r)\int_r^\infty\frac{S(s)^2}{s}
                     (\mathcal I_{\widetilde\theta}f)(s)\,ds,\\
 (\mathcal N_2f)(r)
 &=\frac{4S(r)\cos\widetilde\theta(r)a(r)}{r^2}
                     (\mathcal I_{\widetilde\theta}f)(r).
 \end{aligned}
\end{equation}
Their kernels are displayed in \eqref{ess-spectrum-Ltheta} below;
$\mathcal N_1$ is symmetric, while $\mathcal N_2^*$ is the adjoint
in $L^2_{\rm rad}$.  In the decomposition
\eqref{eq:polar-form-bounded-perturbation},
\begin{equation}\label{eq:explicit-Btheta}
 B_\theta=
 \left(\mathcal U_\theta-\frac{m^2}{r^2}-1\right)I
 +\mathcal N_1+\frac12(\mathcal N_2+\mathcal N_2^*).
\end{equation}
The local multiplier is bounded; the two nonlocal operators are
Hilbert--Schmidt by the profile estimates used in the proof of
\eqref{ess-spectrum-Ltheta}, independent of any spectral conclusion.
Hence $B_\theta$ is bounded and self-adjoint and the polar form extends
to a closed lower-bounded form on $\mathcal H$.
Define its self-adjoint realization by
\begin{equation}\label{def-Ltheta-s}
 \begin{aligned}
 D(L_\theta^s)&=\bigl\{f\in\mathcal H:\ \exists\,g\in L^2_{\rm rad}
 \text{ such that }\mathfrak q_\theta[f,\varphi]
 =\langle g,\varphi\rangle_{L^2_{\rm rad}}
 \ \forall\varphi\in\mathcal H\bigr\},\\
 L_\theta^sf&=g.
 \end{aligned}
\end{equation}
In particular, on its operator domain it acts as
\begin{equation}\label{eq:Ltheta-realization-explicit}
 L_\theta^s=-\partial_r^2-\frac1r\partial_r
                 +\frac{m^2}{r^2}+1+B_\theta.
\end{equation}
For $f$ merely in $\mathcal H$, its quadratic expression is
$\mathfrak q_\theta[f]$.  The formal $L_\theta$ above is generally
not symmetric; its realization in the Hessian is
$L_\theta^s$, namely the symmetrized form operator.

Substitution of the component decomposition into the preceding second
variation gives
\begin{equation}\label{decomp-H}
 D^2H(\widetilde u)
 [\xi_\theta e_\theta+\xi_\phi e_\phi,
  \xi_\theta e_\theta+\xi_\phi e_\phi]
 =\frac{d^2}{ds^2}\bigg|_{s=0}H(u_s)
 =\mathfrak q_\phi[\xi_\phi]+\mathfrak q_\theta[\xi_\theta].
\end{equation}
This identity is first computed on test sections and extends to the
form domain.  The mass direction in the polar tangent component can be read off
directly from \eqref{eq:mass-radial}: for the tangent direction
$\xi_\theta e_\theta+\xi_\phi e_\phi$ one has
\begin{equation}\label{eq:mass-first-variation-polar}
 DM(\widetilde u)[\xi_\theta e_\theta+\xi_\phi e_\phi]
 =-2\pi\int_0^\infty\xi_\theta(r)
                  \sin\widetilde\theta(r)\,r\,dr
 =-\langle\xi_\theta,\sin\widetilde\theta\rangle_{L^2_{\rm rad}}.
\end{equation}
Thus the linearized mass constraint is precisely orthogonality to
$\sin\widetilde\theta$ in the polar component.  We impose
\begin{equation}\label{theta-ND}
 \ker L_\theta^s\cap
 \{f\in\mathcal H:\langle f,\sin\widetilde\theta\rangle_{L^2}=0\}
 =\{0\}.
\end{equation}

\begin{lemma}
There exists a constant $c > 0$ such that 
\[
\mathfrak q_\phi[\xi_\phi] \geq c\|\xi_\phi\|_{L^2}^2
\]
for any $\xi_\phi \in \mathcal{H}$ which is $L^2$ orthogonal to $\sin\widetilde\theta$.
\end{lemma}
\begin{proof}
The equation \eqref{EL-v} is rewritten as
\begin{equation}\label{eqn-theta}
\Delta_r \widetilde\theta -\left(\widetilde k+\frac{(\widetilde h-m)^2}{r^2}\cos\widetilde\theta+V'(1-\cos\widetilde\theta)+\omega\right)\sin\widetilde\theta = 0.
\end{equation}
Multiplying \eqref{eqn-theta} by $\cos\widetilde\theta$ gives the
corresponding distributional equation for $\sin\widetilde\theta$.
This function belongs to $\mathcal H$ by its origin behavior and decay.
Testing against a dense set in $\mathcal H$ gives
$\mathfrak q_\phi[\sin\widetilde\theta,g]=0$ for all $g\in\mathcal H$;
the representation theorem therefore places it in $D(L_\phi)$ and
$L_\phi\sin\widetilde\theta=0$.  Since $0<\widetilde\theta<\pi$ on
$(0,\infty)$, the function $\sin\widetilde\theta$ is positive.  The standard
ground-state representation therefore yields
\[
\mathfrak q_\phi[f]
=2\pi\int_0^\infty \sin^2\widetilde\theta
\left|\left(\frac{f}{\sin\widetilde\theta}\right)'\right|^2r\,dr\geq0
\]
first for compactly supported $f$ and then, by density, on the form domain.
It also shows that the $L^2$ kernel is
$\operatorname{span}\{\sin\widetilde\theta\}$.

Finally, the exponential decay of $\widetilde\theta$
\cite[Theorem~1.2]{JHS} and the convergence of $\widetilde h$ imply
\[
\sigma_{\mathrm{ess}}(L_\phi)=[\omega,\infty).
\]
Thus zero is an isolated simple eigenvalue and the spectral theorem gives
the claimed positive lower bound on its orthogonal complement.
\end{proof}

\begin{lemma}
There exists $c>0$ such that
\[
 \mathfrak q_\theta[f]\ge c\|f\|_{L^2}^2
\]
for every $f\in\mathcal H$ satisfying
$\langle f,\sin\widetilde\theta\rangle_{L^2}=0$.
\end{lemma}
\begin{proof}
Since $\widetilde u$ is a minimizer of $H$, its polar second variation is
nonnegative, and hence
$\mathfrak q_\theta[f]\ge0$ on $\mathcal H$.

We next note that
\begin{equation}\label{ess-spectrum-Ltheta}
\sigma_{\mathrm{ess}}(L_\theta^s)=[\omega,\infty).
\end{equation}
Indeed, at infinity the local part of $L_\theta^s$ is
$-\Delta_r+\omega+O(r^{-2})+o(1)$.  More explicitly, put
$S=\sin\widetilde\theta$ and $a=\widetilde h-m$.  With respect to the measure
$r\,dr$, the two nonlocal terms in $L_\theta$ have kernels
\[
K_1(r,\sigma)
= S(r)S(\sigma)\int_{\max\{r,\sigma\}}^\infty
\frac{S(\rho)^2}{\rho}\,d\rho
\]
and
\[
K_2(r,\sigma)
=\frac{4S(r)\cos\widetilde\theta(r)a(r)}{r^2}
S(\sigma)\boldsymbol 1_{\{\sigma<r\}}.
\]
The regularity $S(r)=O(r^m)$ at the origin and its exponential decay at
infinity show that $K_1,K_2\in
L^2((r\,dr)\otimes(\sigma\,d\sigma))$.  Hence these operators, as well as
their adjoints, are Hilbert--Schmidt.  The remaining difference between the
local part of $L_\theta^s$ and
$-\Delta_r+m^2/r^2+\omega$ is a bounded multiplication operator tending to
zero at infinity.  It is relatively compact by the local Rellich theorem.
Consequently, \eqref{ess-spectrum-Ltheta} follows from Weyl's theorem.

Suppose that the asserted estimate were false.  Then there would be a
sequence $f_n\in\mathcal H$ such that
\[
\|f_n\|_{L^2}=1,\qquad
\langle f_n,\sin\widetilde\theta\rangle_{L^2}=0,
\qquad
\mathfrak q_\theta[f_n]\longrightarrow0.
\]
By \eqref{ess-spectrum-Ltheta}, zero is an isolated eigenvalue of finite
multiplicity whenever it belongs to the spectrum.  If $P$ denotes the
$L^2$-orthogonal projection onto $\ker L_\theta^s$, the spectral theorem
therefore gives
\[
\|(I-P)f_n\|_{L^2}\longrightarrow0.
\]
After passing to a subsequence, $Pf_n$ converges in the finite-dimensional
space $\ker L_\theta^s$ to a function $f$ satisfying
\[
\|f\|_{L^2}=1,\qquad
\langle f,\sin\widetilde\theta\rangle_{L^2}=0.
\]
This contradicts the constrained nondegeneracy condition
\eqref{theta-ND}.  Hence the claimed positive constant $c$ exists.
\end{proof}

\begin{prop}\label{H1-lowerbound}
For a radial tangent vector field $\xi(r)\cdot\widetilde v(r)=0$ with
$R_{-m\chi}\xi\in H^1(\mathbb R^2;\mathbb R^3)$, let
$u_s=R_{-m\chi}v_s$ be an $m$--equivariant differentiable curve through
$u_0=\widetilde u$ with
$\xi=\frac d{ds}|_{s=0}v_s$.  There is a constant $c>0$ such that
\[
\frac{d^2}{ds^2}\bigg|_{s=0}H(u_s)
\ge c\bigl(\|\xi_\theta\|_{\mathcal H}^2
             +\|\xi_\phi\|_{\mathcal H}^2\bigr)
\ge c'\|R_{-m\chi}\xi\|_{H^1}^2
\]
whenever
$\xi=\xi_\theta(-\cos\widetilde\theta,0,\sin\widetilde\theta)
+\xi_\phi(0,1,0)$ and
\[
\langle\xi_\phi,\sin\widetilde\theta\rangle_{L^2}
=\langle\xi_\theta,\sin\widetilde\theta\rangle_{L^2}=0.
\]
Here the constants depend only on the fixed profile and $m$.
\end{prop}
\begin{proof}
Take an arbitrary small $\varepsilon > 0$. Then
\[
\begin{aligned}
&\mathfrak q_\phi[\xi_\phi] \\
=& \varepsilon\mathfrak q_\phi[\xi_\phi] + (1-\varepsilon)\mathfrak q_\phi[\xi_\phi] \\
\geq&\; \varepsilon\int_0^\infty\Bigl[(\xi_\phi')^2
 +\Bigl(-(\widetilde\theta')^2+\widetilde k\cos\widetilde\theta
 +\frac{(\widetilde h-m)^2}{r^2}\cos^2\widetilde\theta
+V'(1-\cos\widetilde\theta)\cos\widetilde\theta
 +\omega\cos\widetilde\theta\Bigr)\xi_\phi^2\Bigr]2\pi r\,dr\\
&\qquad+c(1-\varepsilon)\|\xi_\phi\|_{L^2}^2.
\end{aligned}
\]
Defining
\[
\mathcal V_\phi \coloneqq{}-(\widetilde\theta')^2+\widetilde k\cos\widetilde\theta
 +\frac{(\widetilde h-m)^2}{r^2}\cos^2\widetilde\theta-\frac{m^2}{r^2}
+V'(1-\cos\widetilde\theta)\cos\widetilde\theta+\omega\cos\widetilde\theta.
\]
the regularity at the origin and the decay at infinity of $\widetilde\theta$
show that $\|\mathcal V_\phi\|_\infty < \infty$.
Then one has
\[
\mathfrak q_\phi[\xi_\phi]
\geq \varepsilon\int_0^\infty \left((\xi_\phi')^2 +\frac{m^2}{r^2}\xi_\phi^2\right)\,2\pi rdr
-\varepsilon \|\mathcal V_\phi\|_\infty \|\xi_\phi\|_{L^2}^2
+c(1-\varepsilon)\|\xi_\phi\|_{L^2}^2.
\]
Therefore by taking sufficiently small $\varepsilon > 0$, we see that there exists a small constant $c_1 > 0$ such that
\[
\mathfrak q_\phi[\xi_\phi] \geq c_1\|\xi_\phi\|_{\mathcal{H}}^2.
\]
For the polar component, write
\[
L_\theta^s=-\Delta_r+\frac{m^2}{r^2}+1+B_\theta.
\]
The local part of $B_\theta$ is bounded by the regularity of
$\widetilde\theta$ at zero and its exponential decay at infinity.  The kernel
calculation in the proof of \eqref{ess-spectrum-Ltheta} shows that its
nonlocal part is Hilbert--Schmidt.  Thus $B_\theta$ is bounded on $L^2$.
Splitting off a sufficiently small multiple of
$-\Delta_r+m^2/r^2$ and applying the preceding $L^2$ coercivity estimate
to the remaining part gives
\[
\mathfrak q_\theta[\xi_\theta]
\geq c_2\|\xi_\theta\|_{\mathcal H}^2.
\]

Finally, for the fixed profile $\widetilde u$ with the regularity stated above one has the norm equivalence
\begin{equation}\label{eq:tangent-H1-equivalence}
\|R_{-m\chi}\xi\|_{H^1}^2
\simeq_{\widetilde u,m}
\|\xi_\theta\|_{\mathcal H}^2+\|\xi_\phi\|_{\mathcal H}^2.
\end{equation}
Indeed, radial differentiation produces only the additional bounded
coefficient $\widetilde\theta'\xi_\theta$.  On the two tails,
$\widetilde\theta$ is close to a pole, so the angular derivative is a
uniformly nondegenerate multiple of $m/r$ on the tangent plane; on a fixed
annulus the $r^{-2}$ weight is controlled by the ordinary $L^2$ term.
This proves both directions of \eqref{eq:tangent-H1-equivalence} and
completes the proof.
\end{proof}

\begin{lemma}[Equivariant tubular control]\label{lem:equivariant-tube-Linfty}
For every $\varepsilon>0$ there is $\rho>0$ such that the following
holds.  If $u$ is a sphere-valued $m$--equivariant map in the same
topological component as $\widetilde u$ and
$d_{H^1}(u,\mathcal O(\widetilde u))\le\rho$, then one can choose a comparison phase $\gamma$ so that, with
$\eta=R_{\gamma}u-\widetilde u$,
\[
\|\eta\|_{H^1}\le2\rho,
\qquad
\|\eta\|_{L^\infty}\le\varepsilon.
\]
\end{lemma}
\begin{proof}
It is enough to prove the assertion after a fixed target rotation.  Write
\[
u=R_{-m\chi}v(r),\qquad \theta=\arccos v_3.
\]
Since $|\theta_r|\le|v_r|$ and $\sqrt{v_1^2+v_2^2}=\sin\theta$, the endpoint
conditions give, for every $r>0$,
\begin{align*}
|v(r)+\vec n|^2
&=2(1+v_3(r))
 \le \frac1m\int_0^r
 \left(|v_s|^2+\frac{m^2(v_1^2+v_2^2)}{s^2}\right)s\,ds,\\
|v(r)-\vec n|^2
&=2(1-v_3(r))
 \le \frac1m\int_r^\infty
 \left(|v_s|^2+\frac{m^2(v_1^2+v_2^2)}{s^2}\right)s\,ds.
\end{align*}
Choose $0<r_0<R_0$ so that the corresponding Dirichlet tails of
$\widetilde u$ are as small as desired.  If
$\|u-\widetilde u\|_{H^1}$ is small, the same two tail energies for $u$
are small by the triangle inequality for the gradients, and the displayed
estimates make both maps uniformly close to the same pole on
$(0,r_0)$ and $(R_0,\infty)$.  On the fixed annulus
$[r_0,R_0]$, the radial profiles are one-dimensional and the weighted
$H^1$ norm is equivalent to the usual $H^1$ norm; hence the one-dimensional
Sobolev embedding gives uniform convergence there.  Combining the three
regions proves the claim.
\end{proof}

\begin{lemma}[Weighted control in the sphere chart]
\label{lem:weighted-sphere-chart}
In a sufficiently small equivariant $H^1$ tube, choose a rotated
representative and write $\eta=R_{\gamma}u-\widetilde u$.  Then
\begin{align}
 \|\eta/r\|_2&\le C_{\widetilde u,m}\|\eta\|_{H^1},
                                                   \label{eq:weighted-chord}\\
 \|(\eta-P_{\widetilde u}\eta)/r\|_2
 &\le C_{\widetilde u,m}\|\eta\|_\infty\|\eta\|_{H^1}.
                                                   \label{eq:weighted-normal-error}
\end{align}
The second right-hand side is $o(\|\eta\|_{H^1})$ as the tube radius
tends to zero.
\end{lemma}
\begin{proof}
By Lemma~\ref{lem:equivariant-tube-Linfty}, there are fixed
$0<a<A$ such that both profiles lie in the same south-pole chart on
$(0,a)$ and the same north-pole chart on $(A,\infty)$.
On either tail, with $v$ the profile of $R_{\gamma}u$, one has
$|v_3+\widetilde v_3|\ge c_*>0$ and the exact identity
\[
 v_3-\widetilde v_3
 =-\frac{\sum_{j=1}^2(v_j+\widetilde v_j)(v_j-\widetilde v_j)}
          {v_3+\widetilde v_3}.
\]
Thus $|\eta_3|\le C\sqrt{\eta_1^2+\eta_2^2}$ on the two tails.
Equivariance gives
$m\left(\sum_{j=1}^2\|\eta_j/r\|_2^2\right)^{1/2}\le\|\nabla\eta\|_2$, and on $[a,A]$ the
weight $r^{-1}$ is bounded.  This proves \eqref{eq:weighted-chord}.
The sphere constraint gives
\[
 \eta-P_{\widetilde u}\eta=-\tfrac12|\eta|^2\widetilde u,
\]
so \eqref{eq:weighted-normal-error} follows immediately.  The tubular
$L^\infty$ control makes its coefficient tend to zero.
\end{proof}

\begin{lemma}[Second-order expansion in the equivariant $H^1$ tube]
\label{lem:cubic-remainder}
There is a modulus $\omega_{\widetilde u}:[0,1]\to[0,\infty)$ with
$\omega_{\widetilde u}(\rho)\to0$ as $\rho\downarrow0$ such that,
for every $u$ with $d_{H^1}(u,\mathcal O(\widetilde u))\le\rho$,
if $\gamma$ is chosen as in Lemma~\ref{lem:equivariant-tube-Linfty}
and $\eta=R_\gamma u-\widetilde u$, then
\begin{equation}\label{eq:cubic-remainder}
\left|H(u)-H(\widetilde u)-\frac12D^2H(\widetilde u)
[P_{\widetilde u}\eta,P_{\widetilde u}\eta]\right|
\le \omega_{\widetilde u}(\rho)\|\eta\|_{H^1}^2.
\end{equation}
\end{lemma}
\begin{proof}
By Lemma~\ref{lem:equivariant-tube-Linfty}, after shrinking the tube we
may assume that $\eta=R_{\gamma}u-\widetilde u$ is uniformly small.
Its tangent projection is explicitly
\[
P_{\widetilde u}\eta
=\eta-(\eta\cdot\widetilde u)\widetilde u.
\]
The sphere constraint gives the exact identity
$\eta\cdot\widetilde u=-|\eta|^2/2$.  Consequently
\[
\|\eta-P_{\widetilde u}\eta\|_{H^1}
\le C\bigl(\|\eta\|_\infty+\|\eta\|_{H^1}\bigr)\|\eta\|_{H^1}
=o_{\rho\downarrow0}(\|\eta\|_{H^1}).
\]
Thus the usual local tangent chart of the sphere is valid uniformly in the
present equivariant $H^1$ tube.  The Dirichlet energy has an exact
extrinsic quadratic expansion in $\eta$.  The first variation of $H$ at
$\widetilde u$ is normal to the sphere, and the identity
$\eta\cdot\widetilde u=-|\eta|^2/2$ converts its normal linear contribution
into the second-fundamental-form term contained in the constrained Hessian.
This replacement must also be checked for the singular angular energy.
We retain that energy in the transverse form
$\frac{m^2}{2}\int (v_1^2+v_2^2)r^{-2}\,dx$ rather than treating a
radial third component as a nonzero angular mode.  By
Lemma~\ref{lem:weighted-sphere-chart},
\[
\begin{aligned}
 \left|\sum_{j=1}^2\int\frac{\eta_j^2-(P_{\widetilde u}\eta)_j^2}{r^2}\,dx\right|
 &\le\sum_{j=1}^2\|(\eta_j-(P_{\widetilde u}\eta)_j)/r\|_2
       (\|\eta_j/r\|_2+\|(P_{\widetilde u}\eta)_j/r\|_2)\\
 &\le C\|\eta\|_\infty\|\eta\|_{H^1}^2.
\end{aligned}
\]
Also $\big||\eta|^2-|P_{\widetilde u}\eta|^2\big|=|\eta|^4/4$.  The normal
first-variation multiplier satisfies
$|\widetilde\lambda(r)|\le C(1+r^{-2})$; hence replacing its quadratic
$|\eta|^2$ by $|P_{\widetilde u}\eta|^2$ costs at most
\[
 C\|\eta\|_\infty^2
       (\|\eta\|_2^2+\|\eta/r\|_2^2)
       =o(\|\eta\|_{H^1}^2).
\]
Together with the ordinary $H^1$ replacement above, this justifies the
second-fundamental-form contribution and the angular quadratic form.
No unweighted $H^1$ Hardy estimate for $\eta_3$ is used.

For the gauge term let $v$ be the radial profile of $R_{\gamma}u$, put $\delta v=v-\widetilde v$, and
\[
\delta h(r)=-\int_0^r s\,\delta v_3(s)\,ds.
\]
The exact averaging identity
\[
\frac{\delta h(r)}{r^2}=-\int_0^1t\,\delta v_3(tr)\,dt
\]
shows, for every $1<p<\infty$,
\begin{equation}\label{eq:delta-h-Lp}
\|\delta h/r^2\|_{L^p(\mathbb R^2)}
\le C_p\|\delta v_3\|_{L^p(\mathbb R^2)}.
\end{equation}
The associated vector field
\[
 (\delta h/r)e_\chi
 =\nabla^\perp\Delta^{-1}(-\delta v_3)
\]
is controlled by Hardy--Littlewood--Sobolev.  We record the integrability
needed for that estimate.  Since both maps are sphere-valued and target
rotation leaves the third component unchanged,
\[
 \delta v_3
 =-\frac12\left(
 |R_{\gamma}u-\vec n|^2-|\widetilde u-\vec n|^2\right).
\]
Consequently
\begin{equation}\label{eq:delta-v3-L1}
 \|\delta v_3\|_1
 \le C_{\widetilde u,\rho}\|\eta\|_2,
 \qquad
 \|\delta v_3\|_2\le\|\eta\|_2.
\end{equation}
Interpolation gives
$\|\delta v_3\|_{4/3}\le C_{\widetilde u,\rho}\|\eta\|_{H^1}$, and hence
\begin{equation}\label{eq:delta-h-over-r-L4}
 \|\delta h/r\|_4
 \lesssim\|\delta v_3\|_{4/3}
 \le C_{\widetilde u,\rho}\|\eta\|_{H^1}.
\end{equation}

Expand
\[
 \frac{(h-m)^2}{2r^2}(1-v_3^2)
 =\frac12\left(\frac{m^2}{r^2}-2m\frac h{r^2}
              +\left|\frac h r\right|^2\right)(1-v_3^2)
\]
about $(\widetilde h,\widetilde v_3)$.  This is a finite polynomial in
$v_3$, $h/r^2$, and $h/r$, with $h$ linear in $v_3$.  Every remainder
term beyond second order contains at least three perturbative factors.
The fixed singular angular term has already been handled in its
transverse form by Lemma~\ref{lem:weighted-sphere-chart}; it creates no
uncontrolled third-component quotient.  The remaining representative products are
bounded, for example, by
\[
 \|\delta h/r\|_4^2\|\delta v_3\|_2,
 \qquad
 \|\delta h/r^2\|_2\|\delta v\|_4^2,
 \qquad
 \|\delta h/r^2\|_4\|\delta v\|_4\|\delta v\|_2.
\]
Equations \eqref{eq:delta-h-Lp},
\eqref{eq:delta-v3-L1}, and
\eqref{eq:delta-h-over-r-L4}, together with
$H^1(\mathbb R^2)\hookrightarrow L^p$ for $2\le p<\infty$, therefore
give the gauge remainder bound
$C_{\widetilde u}\|\eta\|_{H^1}^3$.

Finally, Taylor's theorem and $V\in C^3([0,2])$ give a potential remainder
bounded by $C\|\delta v_3\|_{L^3}^3\le C\|\eta\|_{H^1}^3$.
Combining these estimates with the tubular control proves
\eqref{eq:cubic-remainder} after enlarging the modulus
$\omega_{\widetilde u}$.
\end{proof}

\begin{lemma}[Local phase chart for the ground-state orbit]
\label{lem:local-phase-chart}
There exist $\rho_*>0$ and constants $c_*,C_*>0$ such that every
sphere-valued equivariant map with
$d_{H^1}(u,\mathcal O(\widetilde u))<\rho_*$ admits a unique sufficiently close phase
$\gamma(u)\in\mathbb T:=\mathbb R/(2\pi\mathbb Z)$ for which, with
\[
 \eta=R_{\gamma(u)}u-\widetilde u,
 \qquad
 P_{\widetilde u}\eta
      =\xi_\theta e_\theta+\xi_\phi e_\phi,
\]
one has
\begin{equation}\label{eq:phase-chart-orthogonality}
 \langle\xi_\phi,\sin\widetilde\theta\rangle_{L^2}=0
\end{equation}
and
\begin{equation}\label{eq:phase-chart-equivalence}
 c_*d_{H^1}(u,\mathcal O(\widetilde u))
 \le\|\eta\|_{H^1}
 \le C_*d_{H^1}(u,\mathcal O(\widetilde u)).
\end{equation}
The circle-valued map $u\mapsto\gamma(u)$ is continuous in the $H^1$
topology.  Along a continuous path on a time interval it admits a
continuous real lift.
\end{lemma}

\begin{proof}
Consider
\[
 \Psi(\gamma,u)
 =\left\langle
 P_{\widetilde u}(R_{\gamma}u-\widetilde u),(\vec n\times\widetilde u)
 \right\rangle_{L^2}.
\]
At $(0,\widetilde u)$,
\[
 \partial_\gamma\Psi(0,\widetilde u)
 =\|(\vec n\times\widetilde u)\|_2^2
 =\|\sin\widetilde\theta\|_2^2>0.
\]
The implicit function theorem gives a real phase on each local chart
and \eqref{eq:phase-chart-orthogonality}.  On overlapping charts the
phases differ by $2\pi\mathbb Z$.  The branch close to the corresponding
orbit point is selected by positivity of
\[
 \operatorname{Re}\int_{\mathbb R^2}
 e^{i\gamma}(u_1+iu_2)
       \overline{(\widetilde u_1+i\widetilde u_2)}\,dx.
\]
The complex integral is nonzero in a sufficiently small orbit tube, so
these branches define one circle-valued phase.  A continuous path on an
interval has a real lift after fixing its value at one time.  The lower bound in
\eqref{eq:phase-chart-equivalence} is immediate from the definition of the
orbit distance.  For the upper bound, choose a phase $\sigma$ with
$\|R_{\sigma}u-\widetilde u\|_{H^1}\le2d_{H^1}(u,\mathcal O)$.
The implicit-function estimate gives
$|\gamma(u)-\sigma|\le C\|R_{\sigma}u-\widetilde u\|_{H^1}$, and the smooth
rotation action on $H^1$ then yields the upper bound.  Continuity follows
from the same implicit-function argument.
\end{proof}

\begin{proof}[End of proof of Theorem \ref{thm:conditional-stability}]
Let
\[
\mathcal O\coloneqq\{R_\gamma\widetilde u:\gamma\in\mathbb R\},\qquad
d(u,\mathcal O)\coloneqq d_{H^1}(u,\mathcal O).
\]
For $u$ in the neighborhood of Lemma~\ref{lem:local-phase-chart}, choose
$\gamma(u)$ there and write
\[
\eta=R_{\gamma(u)}u-\widetilde u,\qquad
P_{\widetilde u}\eta=\eta-(\eta\cdot\widetilde u)\widetilde u
=\xi_\theta e_\theta+\xi_\phi e_\phi.
\]
Then
\begin{equation}\label{phase-orthogonality}
\langle\xi_\phi,\sin\widetilde\theta\rangle_{L^2}=0.
\end{equation}
Moreover, \eqref{eq:phase-chart-equivalence} identifies
$\|\eta\|_{H^1}$ with the orbit distance up to fixed constants.

Set $\Delta M=M(u)-M(\widetilde u)$.  Since $|u|=|\widetilde u|=1$, one has
$\eta\cdot\widetilde u=-|\eta|^2/2$.  Rotations about the third axis preserve
the mass, and hence
\[
|\langle\xi_\theta,\sin\widetilde\theta\rangle_{L^2}|
\leq |\Delta M|+C\|\eta\|_{L^2}^2.
\]
Consequently, if
\[
b=\frac{\langle\xi_\theta,\sin\widetilde\theta\rangle_{L^2}}
        {\|\sin\widetilde\theta\|_{L^2}^2},
\qquad
\xi_\theta^\perp=\xi_\theta-b\sin\widetilde\theta,
\]
then
\begin{equation}\label{mass-direction-bound}
|b|\leq C\bigl(|\Delta M|+\|\eta\|_{H^1}^2\bigr),
\qquad
\langle\xi_\theta^\perp,\sin\widetilde\theta\rangle_{L^2}=0.
\end{equation}

Applying Proposition \ref{H1-lowerbound} to
$\xi_\theta^\perp e_\theta+\xi_\phi e_\phi$, using
\eqref{mass-direction-bound} for the remaining bounded Hessian terms, and
then applying Lemma~\ref{lem:cubic-remainder}, we obtain
\[
H(u)-H(\widetilde u)
\geq c\,d(u,\mathcal O)^2
-C\,|M(u)-M(\widetilde u)|^2
-\omega_{\widetilde u}(d(u,\mathcal O))d(u,\mathcal O)^2.
\]
After shrinking the tubular neighborhood so that the last coefficient is
at most $c/2$, this gives the local Lyapunov estimate
\begin{equation}\label{local-Lyapunov-bound}
H(u)-H(\widetilde u)
\geq c_0\,d(u,\mathcal O)^2
-C\,|M(u)-M(\widetilde u)|^2.
\end{equation}

Now let $u(t)$ be the solution in the theorem.  Energy and mass conservation
give
\[
H(u(t))-H(\widetilde u)=H(u_0)-H(\widetilde u),\qquad
M(u(t))-M(\widetilde u)=M(u_0)-M(\widetilde u).
\]
Since $\widetilde u$ is a critical point of $H$ and
$\sin\widetilde\theta\in L^2$, the initial closeness assumption implies
\[
|H(u_0)-H(\widetilde u)|\leq C\delta^2,\qquad
|M(u_0)-M(\widetilde u)|\leq C\delta.
\]
Combining these estimates with \eqref{local-Lyapunov-bound} and using the
continuity of $t\mapsto u(t)$ gives, by the standard bootstrap argument,
\[
\sup_{t\in[0,T]}d(u(t),\mathcal O)<\varepsilon
\]
provided $\delta>0$ is sufficiently small.  This proves the theorem.
\end{proof}

\section{Continuation in an \texorpdfstring{$H^1$}{H1} tube}\label{sec:continuation}

Recall that $\n=(0,0,1)$.  An $m$--equivariant map is written
as
\[
 u(r,\chi,t)=R_{-m\chi}v(r,t),\qquad m\ge2.
\]
Set
\begin{equation}\label{eq:hk}
 h(r)=\int_0^r s(1-v_3(s))\,ds,
 \qquad
 k(r)=\int_r^\infty \frac{(v_1(s)^2+v_2(s)^2)}{s}\,(h(s)-m)\,ds .
\end{equation}
In the Coulomb gauge, the radial CSLL equation is
\begin{equation}\label{eq:radial-CSLL}
 v_t=v\times\Delta_rv-\left\{
 \frac{(h-m)^2}{r^2}v_3+V'(1-v_3)+k
 \right\}(\vec n\times v) .
\end{equation}

Let $\widetilde u$ be an $m$--equivariant ground-state profile for
$H_\omega=E+\omega M$, where $\omega>0$, and let
\[
 \mathcal O=\{R_\gamma\widetilde u:\gamma\in\R\}.
\]

\begin{prop}[Continuation in an $H^1$ tube]\label{prop:continuation}
Assume $m\ge2$, $V\in C^{3,1}([0,2])$, and $V'(0)=0$.  There is a number
$\rho_*>0$ with the following property.  Let $u$ be an $H^2$ solution of
\eqref{eq:radial-CSLL} on its maximal forward interval $[0,T^{+})$ and assume
\begin{equation}\label{eq:tube}
 \sup_{0\le t<T^{+}}d_{H^1}(u(t),\mathcal O)\le\rho_*.
\end{equation}
If $T^{+}<\infty$, then
\begin{equation}\label{eq:finite-bound}
 \sup_{0\le t<T^{+}}\|u(t)-\n\|_{H^2}
 \le C_{T^{+},\widetilde u}
 \bigl(1+\|u(0)-\n\|_{H^2}\bigr)<\infty.
\end{equation}
In particular, $T^{+}=\infty$.  The analogous conclusion holds
backward in time.
\end{prop}

\subsection{The Chern--Simons part is a tame zeroth-order perturbation}

The full $m$--equivariant Schr\"odinger-map equation has radial form
\[
 v_t=v\times\Delta_rv-\frac{m^2}{r^2}v_3(\vec n\times v).
\]
Therefore \eqref{eq:radial-CSLL} is exactly
\begin{equation}\label{eq:forced-SM}
 u_t=u\times\Delta u+\Gamma[u](\vec n\times u),
\end{equation}
where, in radial variables,
\begin{equation}\label{eq:Gamma}
 \Gamma[v]=\frac{h(2m-h)}{r^2}v_3-V'(1-v_3)-k.
\end{equation}
Thus no spatial derivative of $v$ occurs in the additional vector field
$\Gamma[u](\vec n\times u)$.

\begin{lemma}[Tame bounds for the CSLL perturbation]\label{lem:tame}
For every $B<\infty$ there is $C_B<\infty$, depending only on
$m$, $B$, and $\|V\|_{C^3([0,2])}$, such that every
$m$--equivariant $u:\R^2\to S^2$ with $u-\n\in H^2$ and
$\|u-\n\|_{H^1}\le B$ satisfies
\begin{align}
 \|\Gamma[u](\vec n\times u)\|_{H^1}&\le C_B,                                      \label{eq:P-H1}\\
 \|\Gamma[u](\vec n\times u)\|_{H^2}&\le C_B\bigl(1+\|u-\n\|_{H^2}\bigr).          \label{eq:P-H2}
\end{align}
\end{lemma}

\begin{proof}
Put $f=1-v_3$ and regard it as a radial function on $\R^2$.  Since
$0\le f\le2$ and
\[
 \int_{\R^2}f=\frac12\|u-\n\|_2^2,
\]
the $L^p$ norm of $f$ is bounded by $C_B$ for every finite $p$.

First consider the coefficient containing $h$.  If
\[
 a(r)=\frac{h(r)}{r^2},
\]
then the exact averaging identity
\begin{equation}\label{eq:average}
 a(x)=\int_0^1 t f(tx)\,dt
\end{equation}
gives, for $j=0,1,2$,
\begin{equation}\label{eq:average-Sobolev}
 \|D^ja\|_2
 \le \int_0^1t^{j+1}\|(D^jf)(t\,\cdot)\|_2\,dt
 =\frac1{j+1}\|D^jf\|_2.
\end{equation}
Let
\[
 A(x)=x^\perp a(|x|)=\frac{h(r)}r e_\chi.
\]
This is precisely $\nabla^\perp\Delta^{-1}f$.  Since
$h(r)\le\min\{r^2,C_B\}$,
\begin{equation}\label{eq:A-infty}
 \|A\|_\infty\le C_B.
\end{equation}
The Hardy--Littlewood--Sobolev and Calder\'on--Zygmund estimates give
\begin{equation}\label{eq:A-estimates}
 \|A\|_4\lesssim\|f\|_{4/3},\qquad
 \|DA\|_p\lesssim_p\|f\|_p,\qquad
 \|D^2A\|_2\lesssim\|Df\|_2.
\end{equation}
Because
\begin{equation}\label{eq:c-decomp}
 \frac{h(2m-h)}{r^2}=2ma-|A|^2,
\end{equation}
equations \eqref{eq:average-Sobolev}--\eqref{eq:A-estimates} and the
product rule imply
\begin{align}
 \left\|\frac{h(2m-h)}{r^2}\right\|_{L^\infty\cap H^1}&\le C_B,
                                                                    \label{eq:c-H1}\\
 \left\|D^2\left(\frac{h(2m-h)}{r^2}\right)\right\|_2
 &\le C_B\bigl(1+\|D^2u\|_2\bigr).                                \label{eq:c-H2}
\end{align}
For example,
$\|D^2(|A|^2)\|_2\lesssim
\|A\|_\infty\|D^2A\|_2+\|DA\|_4^2$, so this term is already
controlled by the $H^1$ bound.

We next estimate $k$.  By \eqref{eq:k-def-main},
$k=T((h-m)(v_1^2+v_2^2))$.  The backward radial Hardy operator
$Tf(r)=\int_r^\infty f(s)\,ds/s$ is bounded on radial
$L^2(\R^2)$.  Indeed, after $r=e^t$, the function
$e^t(Tf)(e^t)$ is the convolution of $e^tf(e^t)$ with
$e^{-s}{\bf1}_{s>0}$, whose $L^1$ norm is one.  Hence
\begin{equation}\label{eq:k-L2}
 \|k\|_2\le\|(h-m)(v_1^2+v_2^2)\|_2\le C_B.
\end{equation}
Moreover,
\begin{equation}\label{eq:k-first}
 k'=-\frac{(h-m)(v_1^2+v_2^2)}{r},
\end{equation}
and therefore
\[
 \|k\|_\infty
 \le C_B\int_0^\infty\frac{(v_1^2+v_2^2)}{r}\,dr
 \le C_B,
 \qquad
 \|Dk\|_2\le C_B\left\|\frac{\sqrt{v_1^2+v_2^2}}{r}\right\|_2\le C_B.
\]
The last inequality uses
$m\sqrt{v_1^2+v_2^2}/r\le|\nabla u|$.
The apparent $r^{-2}$ term cancels at second order.  In fact,
\begin{equation}\label{eq:k-Laplacian}
 \Delta k=k''+\frac1rk'
 =-(1-v_3)(v_1^2+v_2^2)-(h-m)\frac{\partial_r(v_1^2+v_2^2)}r.
\end{equation}
Since
\[
 \frac{|\partial_r(v_1^2+v_2^2)|}{r}
 \le2\frac{\sqrt{v_1^2+v_2^2}}{r}|v_r|
 \le C_m|\nabla u|^2,
\]
the two-dimensional Gagliardo--Nirenberg inequality gives
\begin{equation}\label{eq:k-H2}
 \|D^2k\|_2=\|\Delta k\|_2
 \le C_B+C_B\|\nabla u\|_4^2
 \le C_B\bigl(1+\|D^2u\|_2\bigr).
\end{equation}

Finally, $V'(0)=0$ and the chain rule give
\begin{align*}
 \|V'(f)\|_{L^\infty\cap H^1}&\le C_B,\\
 \|D^2(V'(f))\|_2
 &\le \|V''\|_\infty\|D^2f\|_2
   +\|V'''\|_\infty\|Df\|_4^2
 \le C_B(1+\|D^2u\|_2).
\end{align*}
Equations \eqref{eq:c-H1}, \eqref{eq:c-H2}, and
\eqref{eq:k-L2}--\eqref{eq:k-H2} show
\begin{equation}\label{eq:Gamma-tame}
 \|\Gamma\|_{L^\infty\cap H^1}\le C_B,
 \qquad
 \|D^2\Gamma\|_2\le C_B(1+\|D^2u\|_2).
\end{equation}
Using $|(\vec n\times u)|\le1$ and the product rule,
\[
 \|D^2(\Gamma (\vec n\times u))\|_2
 \lesssim \|D^2\Gamma\|_2
 +\|D\Gamma\|_4\|Du\|_4
 +\|\Gamma\|_\infty\|D^2u\|_2.
\]
Interpolation between the two bounds in \eqref{eq:Gamma-tame}, followed
by Gagliardo--Nirenberg for $Du$, bounds the middle term by
$C_B(1+\|D^2u\|_2)$.  This proves \eqref{eq:P-H1} and
\eqref{eq:P-H2}.
\end{proof}

\subsection{The moving-frame source and the correct \texorpdfstring{$H^2$}{H2} recovery}

Let $J=-u\times$, choose a radial parallel frame $(e,Je)$ on $u^*TS^2$,
and set
\[
 W=u_r+\frac1rJ u_\chi,
 \qquad W=(\Re q)e+(\Im q)Je.
\]
We keep exactly the convention of \eqref{eq:W-local}: $J=-u\times$
and the plus sign in $W=u_r+r^{-1}Ju_\chi$.  With this fixed choice the
angular lift of $q$ has index $m-1$, so its planar field is
$Q(r,\chi)=e^{i(m-1)\chi}q(r)$.

\begin{lemma}[Contribution of a tangent forcing]\label{lem:frame-source}
Suppose $u_t=u\times\Delta u+Z$, where $Z\cdot u=0$ is equivariant.
The contribution of $Z$ to the equation for $q$ is
\begin{equation}\label{eq:frame-source}
 \mathfrak G_u[Z]
 =[\nabla_r^u Z+r^{-1}J\nabla_\chi^u Z]_e-i \alpha_Zq,
 \qquad
 \alpha_Z(r)=-\int_r^\infty\ip{JZ(s)}{u_s(s)}\,ds,                     
\end{equation}
where $[\cdot]_e$ denotes complex frame coordinates.  If
$Z=\Gamma[u](\vec n\times u)$ and
$\|u-\n\|_{H^1}\le B$, then
\begin{align}
 \|e^{i(m-1)\chi}\mathfrak G_u[\Gamma[u](\vec n\times u)]\|_2&\le C_B,             \label{eq:G-L2}\\
 \|\nabla(e^{i(m-1)\chi}\mathfrak G_u[\Gamma[u](\vec n\times u)])\|_2
 &\le C_B\bigl(1+\|u-\n\|_{H^2}\bigr).                          \label{eq:G-H1}
\end{align}
\end{lemma}

\begin{proof}
Formula \eqref{eq:frame-source} follows by applying $\nabla_t^u$ to $W$ and
using compatibility of the pullback connection.  The first term is the
first covariant Cauchy--Riemann derivative of $Z$.  The second term is the
change of the temporal frame connection (with the sign coming from moving it
to the right-hand side of the coordinate compatibility equation): the curvature identity gives
$(\alpha_Z)_r=\ip{JZ}{u_r}$, and $\alpha_Z(\infty)=0$ gives the displayed integral.

For the estimates we now use the actual source $Z=\Gamma (\vec n\times u)$ rather
than an arbitrary tangent field.  Since $(\vec n\times u)=J(P_u\vec n)$,
$[\vec n\times u]_e=i\zeta$ and $[\nabla_r^u(\vec n\times u)+r^{-1}J\nabla_\chi^u(\vec n\times u)]_e=-iu_3q$.  Therefore
\begin{equation}\label{eq:CS-source-specialized}
 \mathfrak G_u[\Gamma (\vec n\times u)]
   =i\Gamma_r\zeta-i\Omega q,\qquad
 \Omega=\Gamma u_3+\alpha_Z,\qquad
 \alpha_Z(r)=\int_r^\infty\Gamma(s)\partial_s u_3(s)\,ds.
\end{equation}
In particular the radial derivative has the exact cancellation
\begin{equation}\label{eq:Omega-cancellation}
 \Omega_r=\Gamma_r u_3.
\end{equation}
The estimate for the connection term is
\begin{equation}\label{eq:a-infty}
 \|\alpha_Z\|_\infty
 \le \|\Gamma\|_\infty
       \left(\sum_{j=1}^2\left\|u_j/r\right\|_{L^2(rdr)}^2\right)^{1/2}\|u_r\|_{L^2(rdr)}
 \le C_B.
\end{equation}
Thus $\|\Omega\|_\infty\le C_B$ by Lemma~\ref{lem:tame}.
Since $|\zeta|\le1$, the undifferentiated estimate follows from
\[
 \|\Gamma_r\zeta\|_2+\|\Omega q\|_2
 \le \|\nabla\Gamma\|_2+\|\Omega\|_\infty\|q\|_2\le C_B.
\]

We make the origin estimate explicit.  For a radial scalar $\Gamma$,
\begin{equation}\label{eq:radial-scalar-Hessian}
 \|D^2\Gamma\|_2^2
 =2\pi\int_0^\infty
      \left(|\Gamma_{rr}|^2+|\Gamma_r/r|^2\right)r\,dr.
\end{equation}
For any radial complex amplitude $y$,
\[
 \|\nabla(e^{i(m-1)\chi}y)\|_2^2
 =2\pi\int_0^\infty
        \left(|y_r|^2+(m-1)^2|y/r|^2\right)r\,dr.
\]
As $\zeta_r=-u_3[u_r]_e$, the product $y=\Gamma_r\zeta$ satisfies
\begin{equation}\label{eq:CS-source-origin-estimate}
 \|\nabla(e^{i(m-1)\chi}\Gamma_r\zeta)\|_2
 \le C_m\|D^2\Gamma\|_2
           +\|\nabla\Gamma\|_4\|u_r\|_4.
\end{equation}
In particular its angular quotient is bounded by
$|\Gamma_r\zeta/r|\le|\Gamma_r/r|$, which is precisely the second
term in \eqref{eq:radial-scalar-Hessian}.  For the other product,
\eqref{eq:Omega-cancellation} gives
\begin{equation}\label{eq:CS-source-Omega-estimate}
 \|\nabla(e^{i(m-1)\chi}\Omega q)\|_2
 \le \|\Omega\|_\infty\|\nabla Q\|_2
              +\|\nabla\Gamma\|_4\|q\|_4.
\end{equation}

The direct derivative-field bound, obtained by differentiating $W$ and
using the angular Hessian identity of
Lemma~\ref{lem:forward-modified-field}, is
\[
 \|\nabla Q\|_2
 \le C_m\bigl(\|D^2u\|_2+\|du\|_4^2\bigr).
\]
Indeed the transverse modes have absolute index $m\ge2$, which controls
$u_{1,r}/r, u_{2,r}/r$ and $u_1/r^2, u_2/r^2$ by $D^2u$; the radial third
component contributes only $\partial_r u_3/r$.  The remaining terms are quadratic
in $du$.  Since $|q|\le C_m|du|$, Gagliardo--Nirenberg and
\eqref{eq:Gamma-tame} yield
\[
 \|u_r\|_4+\|q\|_4\le C_B\|D^2u\|_2^{1/2},\qquad
 \|\nabla\Gamma\|_4\le C_B(1+\|D^2u\|_2)^{1/2}.
\]
Equations \eqref{eq:CS-source-origin-estimate} and
\eqref{eq:CS-source-Omega-estimate} now prove \eqref{eq:G-H1}.
All estimates are first on smooth maps and extend by density.  We have
estimated the lifted amplitude and its angular quotient directly;
no globally smooth extension of the equivariant radial frame at the
origin is assumed.
\end{proof}

We next record the local elliptic recovery estimate in a form that uses the Cartesian Hessian and preserves the cancellations at the origin.

\begin{lemma}[Local elliptic recovery from $q$]\label{lem:recovery}
There exists $\rho_0>0$, depending on $m$ and $\widetilde u$, such that
the following holds.  If $d_{H^1}(u,\mathcal O)\le\rho_0$, choose a constant
$\gamma\in\mathbb R$ for which
\[
 \|R_{\gamma}u-\widetilde u\|_{H^1}\le2\rho_0,
 \qquad \eta_q=e^{-i\gamma}q-\widetilde q.
\]
Then
\begin{equation}\label{eq:recovery}
 \|u-\n\|_{H^2}
 \le C_{\widetilde u,m}
 \left(1+\|u-\n\|_2+
 \|e^{i(m-1)\chi}\eta_q\|_{H^1}\right).
\end{equation}
\end{lemma}

\begin{proof}
By the covariance formula above, $e^{-i\gamma}q$ is the modified field
of the rotated representative $R_\gamma u$.  We replace $u$ by this
representative and write $w=u-\widetilde u$; then
$\|w\|_{H^1}\le2\rho_0$ and its derivative-field difference is $\eta_q$.
The Bogomolny derivative has the exact identity
\begin{equation}\label{eq:tension-W}
 P_u\Delta u=\nabla_r^u W+\frac1r(1-J\nabla_\chi^u)W.
\end{equation}
Equivariance reduces the last term to a bounded scalar multiple of $W/r$.
Since $|\nabla_r^u W|=|q_r|$ and $m-1\ne0$,
\begin{equation}\label{eq:tension-q}
 \|P_u\Delta u\|_2
 \le C_m\|e^{i(m-1)\chi}q\|_{\dot H^1}
 \le C_{\widetilde u,m}
 \left(1+\|e^{i(m-1)\chi}\eta_q\|_{H^1}\right).
\end{equation}
No estimate of $q/r$ is available from this argument when $m=1$; this
is the precise reason for the restriction $m\ge2$.

Since $u$ is sphere-valued,
\[
 \Delta u=P_u\Delta u-|\nabla u|^2u.
\]
For componentwise functions on $\R^2$,
$\|D^2u\|_2=\|\Delta u\|_2$.  Furthermore,
\begin{align*}
 \|\nabla u\|_4^2
 &\le C\|\nabla\widetilde u\|_4^2+C\|\nabla w\|_4^2\\
 &\le C_{\widetilde u}+C\|\nabla w\|_2\|D^2w\|_2\\
 &\le C_{\widetilde u}+C\rho_0
 \bigl(\|D^2u\|_2+\|D^2\widetilde u\|_2\bigr).
\end{align*}
Combining this with \eqref{eq:tension-q} and choosing $\rho_0$ so that
the coefficient of $\|D^2u\|_2$ is at most $1/2$ proves
\eqref{eq:recovery}.  Notice that the argument uses the Cartesian
identity for the full Hessian and therefore retains the cancellations at
$r=0$.
\end{proof}

\subsection{Exact difference equation and short-time Strichartz iteration}

We now write explicitly the equation to which the Gustafson--Shatah
iteration is applied.  Throughout this subsection all derivative fields
are first written in the canonical frame
\eqref{eq:canonical-frame-infinity}.  Thus the exact equation and the
fixed-time frame comparison use the same gauge.  A target rotation about
$\vec n$ acts on the canonical modified field by a constant unit complex
factor.  When the solution is compared with the ground-state orbit we
introduce the corresponding time-dependent factor explicitly, and its
time derivative remains in the difference equation.

Throughout the continuation argument $u$ denotes the solution of
the original CSLL equation \eqref{eq:forced-SM}, with no additional
time-dependent change of target coordinates.  Under a constant target
rotation $R_\sigma$, the scalar quantities $h,k,\Gamma$ are invariant,
whereas the canonical complex fields $q,\zeta,\frac{\zeta}{r}$ are multiplied by
$e^{-i\sigma}$, by the covariance formula in
Section~\ref{sec:lwp-proof}.  We retain the canonical-frame coordinate $\zeta=[P_u\vec n]_e$ and
its nonlocal temporal-connection coefficient
$\alpha_{\mathrm{nl}}$ from \eqref{eq:alpha-exact}.
The reconstruction identities from Section~\ref{sec:lwp-proof} give
\begin{equation}\label{eq:reconstruction-cont}
 \frac{\partial_r u_3}{r}=\Re(\frac{\zeta}{r}\overline q)-m\frac{|\zeta|^2}{r^2},
 \qquad
 \partial_r\bigl(\frac{\zeta}{r}\bigr)+\frac{m+1}{r^2}\zeta
 =-u_3\frac qr+m(1+u_3)\frac{\zeta}{r^2}.
\end{equation}

\begin{lemma}[Exact canonical and modulated equations]
\label{lem:exact-GS-CS}
In the canonical frame,
\begin{equation}\label{eq:rotation-source}
 \mathfrak G_u[\vec n\times u]=-iq.
\end{equation}
Consequently the modified field of a smooth solution of
\eqref{eq:forced-SM} satisfies
\begin{equation}\label{eq:q-original-exact}
 (i\partial_t-\partial_r^2-r^{-1}\partial_r+(m-1)^2r^{-2})q
 =\mathcal V[u,q]q
  +i\mathfrak G_u[\Gamma[u](\vec n\times u)],
\end{equation}
where
\begin{equation}\label{eq:NSM-exact}
 \mathcal V[u,q]q
 =\left\{\frac12|q|^2
 -2m\frac{1+u_3}{r^2}+\alpha_{\mathrm{nl}}\right\}q.
\end{equation}

Let $\gamma\in W^{1,1}_{\rm loc}$ be real-valued and define the
modulated coordinate fields by
\begin{equation}\label{eq:modulated-coordinate-dictionary}
 (\widehat q,\widehat\zeta,\widehat Q)
   =e^{-i\gamma}(q,\zeta,Q),
\end{equation}
where $Q(r,\chi)=e^{i(m-1)\chi}q(r)$.
For the fixed ground-state profile $\widetilde u$, the fields
$\widetilde q=q[\widetilde u]$ and
$\widetilde Q=e^{i(m-1)\chi}\widetilde q$ satisfy the stationary
profile equation
\begin{equation}\label{eq:stationary-modified-ground-field}
 (-\Delta+\omega)\widetilde Q
 =e^{i(m-1)\chi}\bigl\{\mathcal V[\widetilde u,\widetilde q]\widetilde q
       +i\mathfrak G_{\widetilde u}[\Gamma[\widetilde u](\vec n\times\widetilde u)]\bigr\}.
\end{equation}
Set
\[
 \beta=\widehat Q-\widetilde Q.
\]
In the expressions $I,II,III$ below, $q$ and $\zeta$ stand for
$\widehat q$ and $\widehat\zeta$; the unmodulated fields are used
in $\mathfrak G_{\rm CS}$.  Then
\begin{equation}\label{eq:beta-exact-CSLL}
 (i\partial_t-\Delta+\omega)\beta
 =e^{i(m-1)\chi}(I+II+III)+\mathfrak G_{\rm CS}
  +\mathfrak G_{\rm ph},
\end{equation}
where
\begin{align}
 I&=\frac12\bigl(|q|^2q-|\widetilde q|^2\widetilde q\bigr),
                                                               \label{eq:I-GS}\\
 II&=-2m\left\{
 \frac{1+u_3}{r^2}q-
 \frac{1+\widetilde u_3}{r^2}\widetilde q\right\},             \label{eq:II-GS}\\
 III&=q\alpha_{\mathrm{nl}}-\widetilde q\widetilde\alpha_{\mathrm{nl}},
                                                               \label{eq:III-GS}\\
 \mathfrak G_{\rm CS}
 &=ie^{i(m-1)\chi}\left(
 e^{-i\gamma}\mathfrak G_u[\Gamma[u](\vec n\times u)]-
 \mathfrak G_{\widetilde u}[\Gamma[\widetilde u](\vec n\times\widetilde u)]\right),       \label{eq:GCS-exact}\\
 \mathfrak G_{\rm ph}&=(\dot\gamma+\omega)(\beta+\widetilde Q).
                                                               \label{eq:Gphase-exact}
\end{align}
\end{lemma}

\begin{proof}
Since $\vec n\times u=JP_u\vec n$, one has $[\vec n\times u]_e=i\zeta$.  The second identity in
\eqref{eq:local-reconstruction-identities} and the equivariant identity
$J\nabla_\chi^u Y=-mu_3Y$ give
\[
 [\nabla_r^u(\vec n\times u)+r^{-1}J\nabla_\chi^u(\vec n\times u)]_e=-iu_3q.
\]
Moreover $J(\vec n\times u)=-P_u\vec n$, and therefore, using $u_3(\infty)=1$,
\[
 \alpha_{\vec n\times u}(r)
 =-\int_r^\infty\ip{J(\vec n\times u)}{u_s}\,ds
 =\int_r^\infty \partial_s u_3(s)\,ds=1-u_3(r).
\]
Lemma~\ref{lem:frame-source} now yields
$\mathfrak G_u[\vec n\times u]=-iu_3q-i(1-u_3)q=-iq$, proving
\eqref{eq:rotation-source}.  Since we use the original equation
\eqref{eq:forced-SM}, there is no additional rotation-frequency term
in the exact $q$ equation \eqref{eq:q-original-exact}.

By \eqref{eq:alpha-exact} and \eqref{eq:reconstruction-cont},
\[
 \alpha+\frac{m\partial_r u_3}{r}
 =\frac12|q|^2-m^2\frac{|\zeta|^2}{r^2}+\alpha_{\mathrm{nl}}.
\]
Since $\frac{|\zeta|^2}{r^2}=(1-u_3^2)/r^2$,
\[
 -m^2\frac{|\zeta|^2}{r^2}+\frac{(m-1)^2-(1+m u_3)^2}{r^2}
 =-2m\frac{1+u_3}{r^2}.
\]
This proves \eqref{eq:NSM-exact} and \eqref{eq:q-original-exact}.

Under a target rotation $R_\sigma$, the canonical complex fields
$q,\zeta,Q$ acquire the factor $e^{-i\sigma}$, while
$u_3,h,k,\Gamma$ and $\alpha_{\mathrm{nl}}$ remain invariant.
Indeed, every product in the definition of $\alpha_{\mathrm{nl}}$
is phase neutral.  In particular, $R_{\omega t}\widetilde u$
solves the original CSLL equation and
$Q[R_{\omega t}\widetilde u]=e^{-i\omega t}\widetilde Q$.
Inserting this in \eqref{eq:q-original-exact} gives
\eqref{eq:stationary-modified-ground-field}.  The aligned map
corresponding to the modulation $e^{-i\gamma}$ is $R_\gamma u$.

Finally,
$i\partial_t(e^{-i\gamma}Q)=e^{-i\gamma}iQ_t+
\dot\gamma e^{-i\gamma}Q$.  Write the right-hand side of the
original equation as
$\mathcal F=e^{i(m-1)\chi}
 (\mathcal V[u,q]q
  +i\mathfrak G_u[\Gamma[u](\vec n\times u)])$.
When $(i\partial_t-\Delta)Q=\mathcal F$ is written with the stationary
operator $-\Delta+\omega$, this gives
\[
 (i\partial_t-\Delta+\omega)(e^{-i\gamma}Q)
 =e^{-i\gamma}\mathcal F+(\dot\gamma+\omega)e^{-i\gamma}Q.
\]
Phase covariance of all coordinate expressions and subtraction of the
stationary profile equation give
\eqref{eq:beta-exact-CSLL}--\eqref{eq:Gphase-exact}.
\end{proof}

We next choose the phase by a horizontal orthogonality condition.  For
complex-valued functions write
\[
 (f,g)_{\mathbb R}:=\Re\int_{\mathbb R^2}f\overline g\,dx.
\]

\begin{lemma}[Horizontal phase modulation]
\label{lem:horizontal-phase}
Assume first that $\widetilde Q\ne0$.  There exists
$\varepsilon_{\widetilde Q}>0$ such that, whenever
\[
 \inf_{\sigma\in\mathbb R}
 \|e^{-i\sigma}Q-\widetilde Q\|_2
 <\varepsilon_{\widetilde Q},
\]
one may choose a locally unique continuous phase $\gamma$ so that
\begin{equation}\label{eq:horizontal-phase-condition}
 \beta=e^{-i\gamma}Q-\widetilde Q,
 \qquad
 (\beta,i\widetilde Q)_{\mathbb R}=0,
 \qquad
 (\beta+\widetilde Q,\widetilde Q)_{\mathbb R}
 \ge\frac12\|\widetilde Q\|_2^2.
\end{equation}
If
\begin{equation}\label{eq:abstract-phase-equations}
 (i\partial_t-\Delta)Q=\mathcal F,
 \qquad
 (-\Delta+\omega)\widetilde Q=\widetilde{\mathcal F},
\end{equation}
then $\gamma\in W^{1,1}_{\rm loc}$ and, with
$\widehat{\mathcal F}=e^{-i\gamma}\mathcal F$,
\begin{equation}\label{eq:gamma-dot-formula}
 \dot\gamma+\omega
 =\frac{(\beta,\widetilde{\mathcal F})_{\mathbb R}
 -(\widehat{\mathcal F}-\widetilde{\mathcal F},
    \widetilde Q)_{\mathbb R}}
 {(\beta+\widetilde Q,\widetilde Q)_{\mathbb R}}.
\end{equation}
If $\widetilde Q=0$, the phase is arbitrary for alignment; we choose
$\gamma(t)=-\omega t$ so that $\dot\gamma+\omega=0$ and the phase
source vanishes.
\end{lemma}

\begin{proof}
Apply the implicit function theorem to
\[
 \Phi(\gamma,Q)
 =(e^{-i\gamma}Q-\widetilde Q,i\widetilde Q)_{\mathbb R}.
\]
At a point on the phase orbit of $\widetilde Q$,
$\partial_\gamma\Phi=-\|\widetilde Q\|_2^2$.  This gives
\eqref{eq:horizontal-phase-condition}, and continuation along a connected
time interval gives a continuous branch.

Differentiate the orthogonality condition.  The equation for
$\widehat Q=e^{-i\gamma}Q$ is
\[
 (i\partial_t-\Delta+\omega)\widehat Q
 =\widehat{\mathcal F}+(\dot\gamma+\omega)\widehat Q.
\]
Thus, with $\beta=\widehat Q-\widetilde Q$,
\[
 (i\partial_t-\Delta+\omega)\beta
 =\widehat{\mathcal F}-\widetilde{\mathcal F}
  +(\dot\gamma+\omega)(\beta+\widetilde Q).
\]
Taking the real inner product of $\partial_t\beta$ with
$i\widetilde Q$, using self-adjointness of $-\Delta+\omega$, and then
using the stationary equation in
\eqref{eq:abstract-phase-equations}, gives exactly
\eqref{eq:gamma-dot-formula}.  The denominator is bounded below by
\eqref{eq:horizontal-phase-condition}; the asserted time regularity
follows by approximation from the displayed identity.
\end{proof}

We next estimate the pure terms in the exact difference equation.  The
coefficient bounds retain their dependence on the modified field, and the
small quantities are spacetime norms on a short interval.  In particular,
no higher spatial integrability of the coefficients is inferred from the
$H^1$ tube alone.

\begin{lemma}[Gustafson--Shatah estimates in spacetime]
\label{lem:GS-explicit-estimates}
Assume the tube condition \eqref{eq:tube} with sufficiently small
radius $\rho_*$.  Use the modulated coordinates of
Lemma~\ref{lem:exact-GS-CS}, so that in $I,II,III$ one has
$Q=e^{i(m-1)\chi}q=\beta+\widetilde Q$.  For
$I_0=[t_0,t_0+\tau]\subset[0,T^{+})$, $0<\tau\le1$, put
\begin{align*}
 S&=\|\beta\|_{L_t^\infty L_x^2(I_0)}
       +\|\beta\|_{L^4_{t,x}(I_0)},\\
 X&=\|\nabla\beta\|_{L_t^\infty L_x^2(I_0)}
       +\|\nabla\beta\|_{L^4_{t,x}(I_0)},\\
 \mathcal A&=\tau^{1/4}+S,\qquad
 \mathcal B=\tau^{1/4}+X.
\end{align*}
There is a constant $C=C(m,\widetilde u,\rho_*)$, independent of
any $H^2$ bound, such that
\begin{align}
 \|e^{i(m-1)\chi}(I+II+III)\|_{L^{4/3}_{t,x}(I_0)}
 &\le C\mathcal A^3,                                      \label{eq:pure-GS0}\\
 \|\nabla(e^{i(m-1)\chi}(I+II+III))\|_{L^{4/3}_{t,x}(I_0)}
 &\le C\mathcal A^2\mathcal B.                            \label{eq:pure-GS1}
\end{align}
The inequalities hold for equivariant paths in the tube whenever the
displayed norms are finite.
\end{lemma}

\begin{proof}
All unqualified $L^p$ norms in a fixed-time calculation are spatial norms
on $\mathbb R^2$; equivalently they use the measure $2\pi r\,dr$ for a
radial amplitude.  All spacetime norms in this proof are over
$I_0\times\mathbb R^2$.  Recall the reconstruction identities \eqref{eq:local-reconstruction-identities}
and the nonlocal temporal-connection formula \eqref{eq:alpha-exact}.
The spatial identities used below are unaffected by a phase depending
only on time.

\smallskip
\noindent\emph{Step 1: an inner radius fixed by the $H^1$ tube.}
Fix $\varepsilon_0\in(0,1/4]$.  There exist $r_0>0$ and a sufficiently
small tube radius, depending only on $m$, $\widetilde u$, and
$\varepsilon_0$, such that
\[
 0\le1+u_3(r,t)\le\varepsilon_0\qquad(0<r<r_0,\ t\in I_0).
\]
Here the radius is independent of $\|u(t)-\n\|_{H^2}$.
Indeed, $u_3(0)=-1$ and the pointwise identity $|u|=1$ gives
\[
 |\partial_r u_3|\le \sqrt{u_1^2+u_2^2}\,|u_r|.
\]
Consequently,
\[
 2(1+u_3(r))\le\frac1m\int_0^r
 \left(|u_s|^2+\frac{m^2(u_1^2+u_2^2)}{s^2}\right)s\,ds.
\]
The corresponding inner energy of the fixed profile tends to zero with
$r_0$.  A representative $R_\sigma u$ that is close to $\widetilde u$
in $H^1$ has equally small inner energy, by the triangle inequality for
the gradients.  Since $u_3$ and this energy are rotation invariant, the
assertion holds for the coordinates used here.  This is also the inner
part of Lemma~\ref{lem:equivariant-tube-Linfty}.

\smallskip
\noindent\emph{Step 2: fixed-time Hardy bounds with the $q$ dependence retained.}
The regular solution of \eqref{eq:reconstruction-cont} satisfies on
$(0,r_0)$
\[
 (\zeta(r)/r)=r^{-m-1}\int_0^r s^m
 \{-u_3(s)q(s)+m(1+u_3(s))(\zeta(s)/s)\}\,ds.
\]
For $p>2$, Minkowski's inequality and the two-dimensional dilation
law give directly
\[
 \left\|\int_0^1t^m f(t\,\cdot)\,dt\right\|_{L^p(0,r_0;rdr)}
 \le\left(\int_0^1t^{m-2/p}\,dt\right)
       \|f\|_{L^p(0,r_0;rdr)}
 =\frac{\|f\|_{L^p(0,r_0;rdr)}}{m+1-2/p}.
\]
Since $m\varepsilon_0/(m+1-2/p)<1/2$, the term containing $\frac{\zeta}{r}$ is absorbed.
Thus
\[
 \|\frac{\zeta}{r}\|_{L^p(0,r_0)}\le C_{m,p}\|q\|_{L^p(0,r_0)}.
\]
For $r\ge r_0$, one has $\frac{|\zeta|}{r}=\sqrt{u_1^2+u_2^2}/r\le r^{-1}$.
This gives an exterior $L^p$ bound for $p>2$.
Also,
\[
 \frac{1+u_3}{r^2}=\frac{|\zeta|^2}{r^2(1-u_3)}\quad(0<r<r_0),\qquad
 0\le \frac{1+u_3}{r^2}\le\frac2{r^2}\quad(r\ge r_0).
\]
Consequently the fixed-time coefficient estimates are
\begin{equation}\label{eq:GS-absolute-coefficients}
 \begin{aligned}
  \|\zeta(t)/r\|_{L^p}&\le C_p\bigl(1+\|q(t)\|_{L^p}\bigr),\\
  \|(1+u_3(t))/r^2\|_{L^{p/2}}&\le C_p\bigl(1+\|q(t)\|_{L^p}^2\bigr),
  \qquad p>2.
 \end{aligned}
\end{equation}
The $L^2$ norm of $\zeta/r$ is controlled separately by the angular energy.
No $L^1$ bound for $(1+u_3)/r^2$ is asserted:
$u_3(r)\to1$ implies $(1+u_3(r))/r^2\sim2/r^2$ at infinity.

We will also need a derivative bound, only for $p=4$.  The same integral
identity, divided by $r$, reads
\[
 \frac{\zeta(r)}{r^2}=\int_0^1 t^{m+1}
 \left\{-u_3(tr)\frac{q(tr)}{tr}
       +m(1+u_3(tr))\frac{\zeta(tr)}{(tr)^2}\right\}dt.
\]
The corresponding dilation constant is $(m+3/2)^{-1}$.  Absorbing again
gives
\[
 \|\zeta/r^2\|_{L^4(0,r_0)}\le C_m\|q/r\|_{L^4(0,r_0)}.
\]
On the exterior, $|\zeta/r^2|\le r^{-2}$.  The first-order equation and
$|u_3|\le1$ then yield, globally,
\begin{equation}\label{eq:GS-b-radial-derivative}
 \|\partial_r\bigl(\frac{\zeta}{r}\bigr)\|_4+\|\zeta/r^2\|_4
 \le C\bigl(1+\|q/r\|_4\bigr)
 \le C\bigl(1+\|\nabla Q\|_4\bigr).
\end{equation}
The last step follows directly from the nonzero angular index $m-1$:
$|q|/r\le|\nabla Q|/(m-1)$.
The integral identities are first applied to smooth maps.  Their regular
boundary term vanishes because $r^{m+1}\frac{\zeta}{r}=r^m\zeta\to0$ as $r\to0$.
The estimates extend by approximation (or truncated Hardy estimates),
so no boundary value of $\zeta/r^2$ is presupposed.

\smallskip
\noindent\emph{Step 3: spacetime coefficient estimates.}
The fixed profile satisfies $\widetilde Q\in W^{1,4}(\mathbb R^2)$.
This follows from its regular behavior at zero, $m\ge2$, and its
exponential decay.  Thus
\[
 \|Q\|_{L^4_{t,x}}\le C\mathcal A,\qquad
 \|\nabla Q\|_{L^4_{t,x}}\le C\mathcal B.
\]
Integrating \eqref{eq:GS-absolute-coefficients} with $p=4$ and
\eqref{eq:GS-b-radial-derivative} in time gives
\begin{align}
 \|\frac{\zeta}{r}\|_{L^4_{t,x}}&\le C\mathcal A,                         \label{eq:GS-b-spacetime}\\
 \|\frac{1+u_3}{r^2}\|_{L^2_{t,x}}+\|\alpha_{\mathrm{nl}}\|_{L^2_{t,x}}
 &\le C\mathcal A^2,                                      \label{eq:GS-dH-spacetime}\\
 \|q/r\|_{L^4_{t,x}}+\|\zeta/r^2\|_{L^4_{t,x}}
       +\|\partial_r\bigl(\frac{\zeta}{r}\bigr)\|_{L^4_{t,x}}
 &\le C\mathcal B.                                        \label{eq:GS-radial-spacetime}
\end{align}
For $\alpha_{\mathrm{nl}}$, the pointwise-in-time boundedness of $T$ on radial $L^2$ gives
explicitly
\[
 \|\alpha_{\mathrm{nl}}\|_{L^2_{t,x}}\le\||q|^2-\frac{2m}{r}\Re(q\overline\zeta)\|_{L^2_{t,x}}
 \le\|q\|_{L^4_{t,x}}^2
       +2m\|q\|_{L^4_{t,x}}\|\frac{\zeta}{r}\|_{L^4_{t,x}}
 \le C\mathcal A^2.
\]
In particular these are spacetime bounds, not uniform-in-time higher
$L^p$ bounds inferred from the energy.

For the differentiated coefficients use the exact identities
\begin{equation}\label{eq:GS-derivative-identities}
 \partial_r\bigl(\frac{1+u_3}{r^2}\bigr)
 =\frac{\Re(\zeta\overline q)}{r^2}-\frac{m|\zeta|^2}{r^3}
 -\frac{2(1+u_3)}{r^3},
 \qquad
 (\alpha_{\mathrm{nl}})_r=\frac{|q|^2-\frac{2m}{r}\Re(q\overline\zeta)}{r}.
\end{equation}
On $(0,r_0)$ one has $(1+u_3)/r^3\le \frac{|\zeta|}{r}|\zeta/r^2|$, whereas on
$(r_0,\infty)$ one has $(1+u_3)/r^3\le2r^{-3}$.  Therefore
\[
 |\partial_r\bigl(\frac{1+u_3}{r^2}\bigr)|\le \frac{|\zeta|}{r}|q/r|+(m+2)\frac{|\zeta|}{r}|\zeta/r^2|
                +4\boldsymbol1_{\{r\ge r_0\}}r^{-3},
 \qquad
 |(\alpha_{\mathrm{nl}})_r|\le |q/r|(|q|+2m\frac{|\zeta|}{r}).
\]
Since $\boldsymbol1_{\{r\ge r_0\}}r^{-3}\in L^2(\mathbb R^2)$,
H\"older's inequality and \eqref{eq:GS-b-spacetime}--\eqref{eq:GS-radial-spacetime}
yield
\begin{equation}\label{eq:GS-dH-derivative-spacetime}
 \|\nabla\bigl(\frac{1+u_3}{r^2}\bigr)\|_{L^2_{t,x}}+\|\nabla\alpha_{\mathrm{nl}}\|_{L^2_{t,x}}
 \le C\{\tau^{1/2}+\mathcal A\mathcal B\}
 \le C\mathcal A\mathcal B.
\end{equation}
Both $(1+u_3)/r^2$ and $\alpha_{\mathrm{nl}}$ are radial scalars,
so their Cartesian gradient magnitudes equal the absolute values of
their radial derivatives.  The last inequality uses
$\mathcal A\mathcal B\ge\tau^{1/2}$.  This argument does not estimate
$\nabla u_3$ in $L^\infty$ and has no hidden $H^2$ coefficient.

\smallskip
\noindent\emph{Step 4: estimate the full products and subtract the fixed ones.}
In Cartesian variables the three terms are exactly
\[
\begin{aligned}
 e^{i(m-1)\chi}I
 &=\tfrac12(|Q|^2Q-|\widetilde Q|^2\widetilde Q),\\
 e^{i(m-1)\chi}II
 &=-2m(\frac{1+u_3}{r^2}Q-\frac{1+\widetilde u_3}{r^2}\widetilde Q),\\
 e^{i(m-1)\chi}III
 &=\alpha_{\mathrm{nl}}Q-\widetilde\alpha_{\mathrm{nl}}\widetilde Q.
\end{aligned}
\]
Thus no Lipschitz estimate for the difference of the scalar coefficients
$(1+u_3)/r^2$ or $\alpha_{\mathrm{nl}}$ is needed here.
For the nonstationary products, H\"older in spacetime gives
\begin{align*}
 \||Q|^2Q\|_{4/3}&\le\|Q\|_4^3,\\
 \|\frac{1+u_3}{r^2}Q\|_{4/3}&\le\|\frac{1+u_3}{r^2}\|_2\|Q\|_4,\\
 \|\alpha_{\mathrm{nl}}Q\|_{4/3}&\le\|\alpha_{\mathrm{nl}}\|_2\|Q\|_4.
\end{align*}
Each right-hand side is at most $C\mathcal A^3$.  The Cartesian product
rule, with the same spacetime norms, gives separately
\begin{align*}
 \|\nabla(|Q|^2Q)\|_{4/3}
 &\le3\|Q\|_4^2\|\nabla Q\|_4,\\
 \|\nabla(\frac{1+u_3}{r^2}Q)\|_{4/3}
 &\le\|\frac{1+u_3}{r^2}\|_2\|\nabla Q\|_4+\|\nabla\bigl(\frac{1+u_3}{r^2}\bigr)\|_2\|Q\|_4,\\
 \|\nabla(\alpha_{\mathrm{nl}}Q)\|_{4/3}
 &\le\|\alpha_{\mathrm{nl}}\|_2\|\nabla Q\|_4+\|\nabla\alpha_{\mathrm{nl}}\|_2\|Q\|_4.
\end{align*}
Each of these is at most $C\mathcal A^2\mathcal B$ by
\eqref{eq:GS-dH-spacetime} and
\eqref{eq:GS-dH-derivative-spacetime}.  Apply the same spatial estimates
to the fixed profile.  Its three products belong to
$W^{1,4/3}(\mathbb R^2)$ and therefore contribute at most
$C\tau^{3/4}$ after integration over $I_0$.  Since
$\mathcal A^3,\mathcal A^2\mathcal B\ge\tau^{3/4}$, the triangle
inequality proves \eqref{eq:pure-GS0}--\eqref{eq:pure-GS1}.
The first-order identities and the product bounds also justify these
estimates for $H^2$ maps with $Q\in L_t^4W_x^{1,4}$ by approximation.

For clarity, expansion of the two bounds shows the powers relevant to
absorption:
\begin{align*}
 \mathcal A^3
 &\le C\{\tau^{3/4}+\tau^{1/2}S+\tau^{1/4}S^2+S^3\},\\
 \mathcal A^2\mathcal B
 &\le C\{\tau^{3/4}+\tau^{1/2}S+\tau^{1/4}S^2
           +(\tau^{1/2}+\tau^{1/4}S+S^2)X\}.
\end{align*}
In particular the cubic contribution $S^3$, and its differentiated
counterpart $S^2X$, do not carry a small time factor; they are absorbed
using smallness of $S$.
\end{proof}

Before applying Strichartz, note that
Corollary~\ref{cor:local-spacetime} gives
$Q\in C_tH^1\cap L^4_{\rm loc}W_x^{1,4}$ independently of the present
continuation argument.  Multiplication by $e^{-i\gamma(t)}$ preserves
all spatial Sobolev norms.  Thus $S(I),X(I)<\infty$ for every compact
$I$ in the maximal interval.  For a fixed left endpoint, the $L^4$
norms are continuous as the right endpoint varies, by absolute
continuity of the integral, and the supremum norms are continuous by
$C_tH^1$.  This supplies the a priori finiteness and continuity needed
for the bootstrap below; absorption is not used to establish them.

Applying the two-dimensional Strichartz estimate to
\eqref{eq:beta-exact-CSLL} (the constant potential $\omega$ only adds a
unit scalar phase) and using Lemma~\ref{lem:GS-explicit-estimates} gives
\begin{align}
 S(I)&\le C\|\beta(t_0)\|_2+C(\tau^{1/4}+S(I))^3\notag\\
 &\quad+\|\mathfrak G_{\rm CS}\|_{L_t^1L_x^2(I)}
       +\|\mathfrak G_{\rm ph}\|_{L_t^1L_x^2(I)},             \label{eq:GS0}\\
 X(I)&\le C\|\nabla\beta(t_0)\|_2
       +C(\tau^{1/4}+S(I))^2(\tau^{1/4}+X(I))\notag\\
 &\quad+\|\nabla\mathfrak G_{\rm CS}\|_{L_t^1L_x^2(I)}
       +\|\nabla\mathfrak G_{\rm ph}\|_{L_t^1L_x^2(I)}.     \label{eq:GS1}
\end{align}

The Chern--Simons contribution is now completely explicit through
\eqref{eq:GCS-exact}.  Lemma~\ref{lem:frame-source} and the fixed
regularity of $\widetilde u$ give
\begin{align}
 \|\mathfrak G_{\rm CS}\|_{L^1_tL^2_x(I)}&\le C_B\tau,          \label{eq:CS0}\\
 \|\nabla\mathfrak G_{\rm CS}\|_{L^1_tL^2_x(I)}
 &\le C_B\tau\left(1+
 \sup_{t\in I}\|u(t)-\n\|_{H^2}\right).                       \label{eq:CS1}
\end{align}
By Lemma~\ref{lem:recovery}, \eqref{eq:CS1} is bounded by
\begin{equation}\label{eq:CS1-X}
 C_B\tau(1+X(I)).
\end{equation}

We now estimate the explicit phase source.  If $\widetilde Q=0$, choose
$\gamma(t)=-\omega t$ as above; then the phase source vanishes.  Otherwise
$\widetilde{\mathcal F}=(-\Delta+\omega)\widetilde Q$ is fixed and
\[
 \widehat{\mathcal F}-\widetilde{\mathcal F}
 =e^{i(m-1)\chi}(I+II+III)+\mathfrak G_{\rm CS}.
\]
The denominator in \eqref{eq:gamma-dot-formula} is at least
$\|\widetilde Q\|_2^2/2$.  Pairing the pure source with the fixed
$\widetilde Q\in L_x^4$, and the CS source with
$\widetilde Q\in L_x^2$, gives
\begin{align}
 \|\dot\gamma+\omega\|_{L^1(I)}
 &\le C_{\widetilde u}\Bigl[
  \tau S(I)+\tau^{1/4}
       \|e^{i(m-1)\chi}(I+II+III)\|_{L^{4/3}_{t,x}(I)}
       +\|\mathfrak G_{\rm CS}\|_{L_t^1L_x^2(I)}\Bigr]\notag\\
 &\le C_\gamma\left\{
       \tau(1+S(I))+\tau^{1/4}(\tau^{1/4}+S(I))^3\right\}
   =:\mathfrak g_\tau(S(I)).                              \label{eq:gamma-L1}
\end{align}
The constant $C_\gamma$ depends only on the fixed profile, the $H^1$
tube, and the potential.  In particular,
$\sup_{0\le s\le1}\mathfrak g_\tau(s)\to0$ as $\tau\downarrow0$.
Consequently
\begin{align}
 \|\mathfrak G_{\rm ph}\|_{L_t^1L_x^2(I)}
 &\le\mathfrak g_\tau(S(I))
       \{S(I)+\|\widetilde Q\|_2\},                       \label{eq:Gphase0}\\
 \|\nabla\mathfrak G_{\rm ph}\|_{L_t^1L_x^2(I)}
 &\le\mathfrak g_\tau(S(I))
       \{X(I)+\|\nabla\widetilde Q\|_2\}.                \label{eq:Gphase1}
\end{align}
For $\widetilde Q=0$, we use the same notation with
$\mathfrak g_\tau\equiv0$.  Thus the phase source has a coefficient
uniformly small in time on the $S$--bootstrap region, independently of
$X$.

\begin{proof}[Proof of Proposition~\ref{prop:continuation}]
For each fixed $t$, choose $\sigma$ so that
$\|R_\sigma u(t)-\widetilde u\|_{H^1}\le2\rho$.
The covariance formula and
Lemma~\ref{lem:canonical-frame-H1-tube}, applied with
$u_*=\widetilde u$, give
\begin{equation}\label{eq:q-small}
 d_{H^1}(u(t),\mathcal O)\le\rho
 \quad\Longrightarrow\quad
 \inf_{\sigma\in\mathbb R}
 \|e^{-i\sigma}Q(t)-\widetilde Q\|_2
 \le\varepsilon(\rho):=2C_*\rho.
\end{equation}
Here $C_*$ and the allowed tube radius depend only on the fixed profile
and $m$, not on the current $H^2$ norm.  No time-dependent change of
frame is made in deriving this fixed-time geometric bound.
If $\widetilde Q\ne0$, choose the continuous horizontal phase of
Lemma~\ref{lem:horizontal-phase} on the maximal interval.  The implicit
function theorem also gives
$|\gamma-\sigma|\le C_{\widetilde Q}
\|e^{-i\sigma}Q-\widetilde Q\|_2$ for the comparison phase above.
Hence $R_\gamma u$ remains in the $H^1$ recovery tube of
$\widetilde u$.  If $\widetilde Q=0$, set $\gamma(t)=-\omega t$ in the equation;
for the phase-invariant recovery estimate one may still use any constant
comparison phase $\sigma(t)$ at each fixed time.  In either case
$\sup_t\|\beta(t)\|_2$ is as small as required once the tube radius is
small.

Choose $s_*\in(0,1/2)$ so small that, for the constants in
\eqref{eq:GS0}--\eqref{eq:GS1},
\[
 C(3s_*)^3\le s_*/4,\qquad C(3s_*)^2\le1/4.
\]
Next choose $\rho_*\le\rho_0$ so that \eqref{eq:q-small} and the phase
choice above give $C\sup_t\|\beta(t)\|_2\le s_*/4$.
Finally choose $0<\tau\le1$ such that $\tau^{1/4}\le s_*$ and
\[
 \begin{aligned}
 C_B\tau+\mathfrak g_\tau(2s_*)
       (2s_*+\|\widetilde Q\|_2)&\le s_*/4,\\
 C_B\tau+\mathfrak g_\tau(2s_*)&\le1/4.
 \end{aligned}
\]
All these choices depend only on the fixed profile, the potential, and
the $H^1$ tube.  On the continuity-bootstrap region $S(I)\le2s_*$,
\eqref{eq:GS0}, \eqref{eq:CS0}, and \eqref{eq:Gphase0} give
$S(I)\le3s_*/4$.  Hence
\begin{equation}\label{eq:S-small}
 S(I)\le s_*
\end{equation}
on every interval of length at most $\tau$, uniformly in its left
endpoint.  The bootstrap starts from the small $L^2$ norm at that endpoint
and the absolute continuity of the spacetime $L^4$ norm.

In the differentiated estimate, \eqref{eq:GS1}, \eqref{eq:CS1-X}, and
\eqref{eq:Gphase1} bound the coefficient of $X(I)$ by
\[
 C(\tau^{1/4}+S(I))^2+C_B\tau+\mathfrak g_\tau(S(I))\le1/2.
\]
The remaining terms are bounded by a constant depending only on the same
fixed quantities.  Absorption therefore yields
\begin{equation}\label{eq:X-recurrence}
 X([t_j,t_{j+1}])
 \le C_0\|\nabla\beta(t_j)\|_2+C_1,
 \qquad t_{j+1}-t_j\le\tau,
\end{equation}
where $C_0,C_1$ are independent of $j$ and of the current $H^2$ norm.
In particular,
\[
 \|\nabla\beta(t_{j+1})\|_2
 \le C_0\|\nabla\beta(t_j)\|_2+C_1.
\]
Covering $[0,T^{+})$ by at most
$N=\lceil T^{+}/\tau\rceil+1$ intervals and iterating yields
\begin{equation}\label{eq:X-finite}
 \sup_{0\le t<T^{+}}\|\nabla\beta(t)\|_2
 \le C_{T^{+},\widetilde u}
 \bigl(1+\|\nabla\beta(0)\|_2\bigr)<\infty.
\end{equation}
The definition of $q$ gives
$\|\nabla\beta(0)\|_2\le
C_B(1+\|u(0)-\n\|_{H^2})$.  Lemma~\ref{lem:recovery} and conservation
of mass, namely
\[
 \|u(t)-\n\|_2^2=2M(u(t))=2M(u(0)),
\]
now turn \eqref{eq:X-finite} into \eqref{eq:finite-bound}.  The flow
only depends on $V'$, so replacing $V$ by $V-V(0)$, if necessary,
puts the equation under the normalization $V(0)=0$ of
Theorem~\ref{thm:lwp}.  Its $H^2$ blow-up alternative contradicts
maximality when $T^{+}<\infty$.
\end{proof}

\subsection{Ground-state regularity and application}

For a finite-energy ground state with $\omega>0$, the regularity and decay
needed above are automatic.  Write
\[
\widetilde v=(\sin\widetilde\theta,0,\cos\widetilde\theta).
\]
Its Euler--Lagrange equation is
\[
 \widetilde\theta''+\frac1r\widetilde\theta'
 -\left\{
 \widetilde k+\frac{(\widetilde h-m)^2}{r^2}
 \cos\widetilde\theta+V'(1-\cos\widetilde\theta)+\omega
 \right\}\sin\widetilde\theta=0.
\]
At the origin the regular-singular alternatives are $r^m$ and $r^{-m}$;
finite energy selects
$\pi-\widetilde\theta=O(r^m)$.  At infinity the equation is
\[
 \widetilde\theta''+r^{-1}\widetilde\theta'
 -\left(\omega+O(r^{-2})+o(1)\right)\widetilde\theta=0.
\]
Comparison with $e^{-cr}$, for any $0<c<\sqrt\omega$, followed by the
ODE and the formulas for $h,k$, gives exponential decay of all derivatives
needed in \eqref{eq:GS0}--\eqref{eq:GS1}.  Standard ODE bootstrapping
gives the finite regularity needed here under $V\in C^{3,1}$:
in particular $\widetilde u-\vec n\in H^2$,
$\widetilde Q\in H^2\cap W^{1,4}$, and
$(-\Delta+\omega)\widetilde Q\in L^2$.
No $C^\infty$ regularity of a ground state for a merely $C^{3,1}$
potential is assumed.  The fixed-profile comparison lemma requires only
$u_*\in\mathcal X_m^2$, and hence applies with this regularity.

If the constrained nondegeneracy condition (ND) yields orbital stability,
then the latter supplies \eqref{eq:tube} on the maximal $H^2$ existence
interval.  Proposition~\ref{prop:continuation} then proves global $H^2$
existence.  Thus (ND) is used only to keep the solution in the $H^1$ tube;
the continuation argument itself is independent of the Hessian.

\section{BPS and near-BPS nondegeneracy}\label{sec:bps}

\subsection{Nondegeneracy of the critical radial BPS ground state}

We now specialize to
\[
V(t)=t^2,\qquad \omega=0,
\]
which corresponds to the self-dual case $\kappa=1$.  Let
\[
\widetilde u=R_{-m\chi}(\sin\widetilde\theta,0,\cos\widetilde\theta)
\]
be the critical radial BPS ground state obtained from the work of Han and
Huh \cite{HanHuh}, as recalled and reconstructed for CSLL in
\cite[Theorem~1.1 and Lemma~3.3]{JHS}.  Here ``critical'' refers to the endpoint of the self-dual
asymptotic-slope classification, made explicit in
\eqref{eq:BPS-slope-definition} below.  This slope parameter is
unrelated to the difference field $\beta$ in Section~\ref{sec:continuation}.
The acronym BPS is used in the
Bogomolny--Prasad--Sommerfield sense; this solution is not a
Belavin--Polyakov instanton.

To specify the external result being used, the common spatial self-dual
system of the Maxwell gauged $O(3)$ model and the present CSLL model is
\[
 D_1\phi+\phi\times D_2\phi=0,\qquad F_{12}=1-\phi_3.
\]
For the south-pole stereographic convention the corresponding real
logarithmic variable is $w=2\log\tan(\theta/2)$.
The asymptotic slope used in the self-dual classification is
\begin{equation}\label{eq:BPS-slope-definition}
 \beta_{\rm BPS}:=-\lim_{r\to\infty}r w_r(r).
\end{equation}
The critical solution has $\beta_{\rm BPS}=2$; for an interior BPS
member $\Theta_\nu$ the slope is $2(1+\nu)$.
With the single south-pole
point of multiplicity $m$ at zero and no north-pole zero, this gives
\begin{equation}\label{eq:BPS-distributional-scalar}
 \Delta w=\frac{4e^w}{1+e^w}-4\pi m\delta_0,
       \qquad w(x)\longrightarrow-\infty\quad(|x|\to\infty).
\end{equation}
This is precisely \cite[(1.13), Theorem~1.1, and Lemma~3.3]{JHS},
which recalls the self-dual existence and critical asymptotic from
\cite{HanHuh}.  Only the spatial self-dual equations and their scalar
classification are used here; the Maxwell and CSLL temporal potentials
and dynamics are not identified.

Before using spectral notation at $\omega=0$, we extend the form
realization to the algebraic BPS tails.  This step does not use
nondegeneracy or a positive spectral gap.

\begin{lemma}[Form realization for BPS profiles]
\label{lem:BPS-form-realization}
Let $\theta$ be a finite-mass radial BPS profile for $V(t)=t^2$,
$\omega=0$, with $0<\theta<\pi$ and
\[
 \theta'+\frac{m-h_\theta}{r}\sin\theta=0,\qquad
 S(r):=\sin\theta(r)=O(r^m)\quad(r\downarrow0).
\]
Suppose its tail is either $S(r)=O((r\log r)^{-1})$ or
$S(r)=O(r^{-1-\nu})$ for some $\nu>0$.  The azimuthal and polar forms
of Section~\ref{sec:stability} are then closed lower-bounded forms on
the same $\mathcal H$ and define self-adjoint operators $L_\phi$ and
$L_\theta^s$ on radial $L^2$.
\end{lemma}
\begin{proof}
Put $a=h_\theta-m$.  Finite mass bounds $a$, and the assumed endpoint
bounds give $S\in L^2(r\,dr)$ and $S/r\in L^2(r\,dr)$.
The formula for $k$ therefore bounds it uniformly.  Near zero
$h_\theta=O(r^2)$, $\theta'=O(r^{m-1})$, and
$\cos(2\theta)-1=O(r^{2m})$.  Thus each local Hessian multiplier minus
$m^2/r^2+1$ is bounded.  At infinity boundedness follows directly from
the indicated tails and the BPS equation.

For the nonlocal kernels in the proof of
\eqref{ess-spectrum-Ltheta}, put
$\int_0^\infty S(\rho)^2\,d\rho/\rho<\infty$.  Then
\[
 |K_1(r,s)|\le
 \left(\int_0^\infty\frac{S(\rho)^2}{\rho}\,d\rho\right)|S(r)S(s)|,
\]
so $K_1$ is Hilbert--Schmidt.  Similarly its companion satisfies
\[
 \|K_2\|_{\mathrm{HS}}^2
 \le C\int_0^\infty\frac{S(r)^2a(r)^2}{r^3}
              \left(\int_0^r S(s)^2s\,ds\right)dr<\infty.
\]
The integrand is $O(r^{4m-1})$ at zero.  At infinity the inner integral
is bounded and either assumed tail makes the remaining integral finite.
Hence the same bounded-perturbation and representation-theorem argument
as in Section~\ref{sec:stability} applies, without exponential decay.

We also record the continuity needed when extending the BPS square
identity from test functions.  With $\mathcal I_\theta$ defined in
Section~\ref{sec:stability},
\[
 \|\mathcal I_\theta f\|_\infty
 \le\|S\|_{L^2(rdr)}\|f\|_{L^2(rdr)}.
\]
Consequently the first-order map
\[
 \mathcal B_\theta f
 =f'+\frac{m-h_\theta}{r}\cos\theta\,f
        -\frac Sr\mathcal I_\theta f
\]
is continuous from $\mathcal H$ to $L^2(r\,dr)$.  Thus a square identity
for test functions involving $\mathcal B_\theta$ extends by density to
all of $\mathcal H$; its polarization is equally valid there.
\end{proof}

The critical tail and the origin behavior used in this lemma are checked
below directly from the BPS equation and the quoted critical asymptotic.
For interior members the corresponding endpoint bounds are established
in Lemma~\ref{lem-BPS-family}.  Neither verification uses a Hessian
spectral assertion.

In this subsection, $\widetilde h=h_{\widetilde\theta}$ and $L_\theta^s$ denotes
the self-adjoint polar Hessian \eqref{def-Ltheta-s}, evaluated at this
ground state.  The following result is stronger than the constrained
nondegeneracy condition \eqref{theta-ND}: no mass orthogonality is needed.

\begin{prop}[Polar nondegeneracy of the critical radial BPS ground state]
\label{prop-BPS-nondegeneracy}
Let $m\geq2$.  For the critical radial BPS ground state corresponding to
$V(t)=t^2$ and $\omega=0$, one has
\[
\ker L_\theta^s=\{0\}
\]
as a self-adjoint operator on $L^2((0,\infty),2\pi r\,dr)$ with form
domain $\mathcal H$.  In particular, \eqref{theta-ND} holds without
imposing $\langle f,\sin\widetilde\theta\rangle_{L^2}=0$.
\end{prop}

\begin{proof}
We give the factorization and the kernel argument in detail.  For a scalar
angle $\theta$, set
\[
h_\theta(r)=\int_0^r s(1-\cos\theta(s))\,ds.
\]
In the present case the energy of the corresponding equivariant map is
\begin{equation}\label{BPS-radial-energy}
E(\theta)=2\pi\int_0^\infty
\left\{
\frac12(\theta')^2
+\frac{(m-h_\theta)^2}{2r^2}\sin^2\theta
+(1-\cos\theta)^2
\right\}r\,dr.
\end{equation}
We first rewrite this energy as a square.  Since
\[
h_\theta'(r)=r(1-\cos\theta(r)),
\]
direct differentiation gives
\[
\frac{d}{dr}\bigl[(m-h_\theta)(1-\cos\theta)\bigr]
=-r(1-\cos\theta)^2
+(m-h_\theta)\sin\theta\,\theta'.
\]
The boundary conditions $\theta(0)=\pi$, $\theta(\infty)=0$ and
$h_\theta(0)=0$ imply
\[
\left[(m-h_\theta)(1-\cos\theta)\right]_{r=0}^{r=\infty}
=-2m.
\]
Consequently,
\begin{equation}\label{BPS-potential-identity}
2\pi\int_0^\infty r(1-\cos\theta)^2\,dr
=4\pi m
+2\pi\int_0^\infty
(m-h_\theta)\sin\theta\,\theta'\,dr.
\end{equation}
Substituting \eqref{BPS-potential-identity} into
\eqref{BPS-radial-energy} and completing the square, we obtain
\begin{equation}\label{BPS-factorization}
E(\theta)-4\pi m
=\pi\int_0^\infty r
\left(
\theta'+\frac{m-h_\theta}{r}\sin\theta
\right)^2\,dr.
\end{equation}
This is the BPS factorization in the radial variables.  In particular,
the critical radial BPS profile satisfies the first-order equation
\begin{equation}\label{BPS-first-order-theta}
\widetilde\theta'
+\frac{m-\widetilde h}{r}\sin\widetilde\theta=0.
\end{equation}

We next differentiate the square in \eqref{BPS-factorization}.  Initially
take $f\in C_c^\infty((0,\infty))$ and put
$\theta_s=\widetilde\theta+sf$.  The variation of $h_{\theta_s}$ at $s=0$ is
\begin{equation}\label{BPS-h-variation}
q_f(r):=D_\theta h_{\widetilde\theta}[f](r)
=(\mathcal I_{\widetilde\theta}f)(r)
=\int_0^r \sigma\sin\widetilde\theta(\sigma)f(\sigma)\,d\sigma.
\end{equation}
It follows that
\begin{align}
&\left.\frac{d}{ds}\right|_{s=0}
\left(
\theta_s'+\frac{m-h_{\theta_s}}{r}\sin\theta_s
\right) \notag\\
&\qquad =f'
+\frac{m-\widetilde h}{r}\cos\widetilde\theta\,f
-\frac{\sin\widetilde\theta}{r}q_f.             \label{BPS-linearized-first-order}
\end{align}
Because the expression being squared in \eqref{BPS-factorization}
vanishes at $\widetilde\theta$, differentiating twice gives
\begin{equation}\label{BPS-Hessian-square}
\mathfrak q_\theta[f]
=2\pi\int_0^\infty r
\left|
f'+\frac{m-\widetilde h}{r}\cos\widetilde\theta\,f
-\frac{\sin\widetilde\theta}{r}q_f
\right|^2\,dr.
\end{equation}
Recall that the vector $e_\theta$ used above equals the negative of the
derivative of $(\sin\theta,0,\cos\theta)$ with respect to $\theta$.
Thus the scalar variation $f$ here corresponds to the polar coefficient
$-f$ in the earlier convention; this change of sign does not alter the
quadratic identity.  Lemma~\ref{lem:BPS-form-realization} proves that
both the Hessian form and the first-order map in the square are continuous
on $\mathcal H$.  Density therefore extends
\eqref{BPS-Hessian-square}, and its polarized identity, to that domain.

Suppose now that $f\in\ker L_\theta^s$.  By the representation theorem,
$\mathfrak q_\theta[f]=\langle L_\theta^sf,f\rangle_{L^2}=0$.
The right-hand side of \eqref{BPS-Hessian-square} must therefore vanish,
and hence
\begin{equation}\label{BPS-kernel-first-order}
f'+\frac{m-\widetilde h}{r}\cos\widetilde\theta\,f
-\frac{\sin\widetilde\theta}{r}q_f=0.
\end{equation}
For brevity, write
\[
S=\sin\widetilde\theta,\qquad z=\frac{f}{S}.
\]
The profile satisfies $0<\widetilde\theta<\pi$ on $(0,\infty)$, so $S>0$
there.  Differentiating $S$ and using \eqref{BPS-first-order-theta}, we
find
\begin{equation}\label{BPS-S-derivative}
S'=\cos\widetilde\theta\,\widetilde\theta'
=-\frac{m-\widetilde h}{r}\cos\widetilde\theta\,S.
\end{equation}
Substituting $f=Sz$ and \eqref{BPS-S-derivative} into
\eqref{BPS-kernel-first-order}, the terms containing
$(m-\widetilde h)\cos\widetilde\theta/r$ cancel and give
\[
z'=\frac{q_f}{r}.
\]
On the other hand, differentiating \eqref{BPS-h-variation} gives
\[
q_f'=rSf=rS^2z.
\]
Thus every kernel element gives a solution of
\begin{equation}\label{BPS-zq-system}
z'=\frac{q_f}{r},\qquad q_f'=rS^2z,
\end{equation}
or, equivalently,
\begin{equation}\label{BPS-z-equation}
(rz')'=rS^2z.
\end{equation}

We record the endpoint behavior needed for this equation.  Near the
origin, \eqref{BPS-first-order-theta} and $\widetilde h(r)=O(r^2)$ give
\[
S(r)=c_0r^m(1+o(1))
\qquad (r\longrightarrow0)
\]
for some $c_0>0$.  Indeed,
\[
\frac{d}{dr}\log\tan\frac{\widetilde\theta}{2}
=\frac{\widetilde\theta'}{\sin\widetilde\theta}
=-\frac{m-\widetilde h}{r}.
\]
Since $\widetilde h(r)=O(r^2)$, integration shows that
$\tan(\widetilde\theta(r)/2)=c r^{-m}(1+o(1))$.  Using
$\widetilde\theta(r)\to\pi$ then gives the stated behavior of $S$.

The regular solution of
\eqref{BPS-z-equation} therefore satisfies
\begin{equation}\label{BPS-z-origin}
z(r)\longrightarrow z_0\in\mathbb R,
\qquad q_f(r)\longrightarrow0
\qquad (r\longrightarrow0).
\end{equation}
Here is a direct explanation of this assertion.  Since
$S^2=O(r^{2m})$, the two possible leading behaviors for the
regular-singular equation \eqref{BPS-z-equation} are a constant and
$\log r$.  The coefficient of the logarithmic term is
$\lim_{r\to0}rz'(r)$.  But $rz'=q_f$, whereas the integral definition
\eqref{BPS-h-variation} and Cauchy--Schwarz give $q_f(r)\to0$.  The
logarithmic term is therefore absent, which proves
\eqref{BPS-z-origin}.

If $z_0=0$, integrating \eqref{BPS-zq-system} twice gives
\[
z(r)=\int_0^r \tau S(\tau)^2z(\tau)
\log\frac r\tau\,d\tau.
\]
For sufficiently small $R>0$, the integral of the nonnegative kernel on
the right is smaller than one uniformly for $0<r\leq R$.  Taking the
supremum of $|z|$ on $[0,R]$ therefore gives $z=0$ there.  Standard
uniqueness away from the origin then gives $z=q_f=0$ on $(0,\infty)$,
and hence $f=0$.  Suppose instead that $z_0>0$.  Then $z>0$ near zero,
and \eqref{BPS-zq-system} gives successively
\[
q_f'(r)=rS(r)^2z(r)>0,\qquad q_f(r)>0,
\qquad z'(r)=\frac{q_f(r)}r>0.
\]
These inequalities prevent $z$ from reaching zero and hence hold for all
$r>0$.  Fixing any $R>0$, we obtain for $r\geq R$
\begin{equation}\label{BPS-z-log-lower}
z(r)=z(R)+\int_R^r\frac{q_f(s)}s\,ds
\geq q_f(R)\log\frac rR.
\end{equation}
If $z_0<0$, the same argument applied to $-z$ gives the corresponding
estimate for $|z|$.

It remains to use the fact that the BPS state is the critical
$\beta_{\rm BPS}=2$ solution.  If $w$ denotes the scalar stereographic variable in
\cite{HanHuh}, then
\[
w(r)=-2\log r-2\log\log r+O(1)
\qquad (r\longrightarrow\infty)
\]
and
\[
S(r)=\frac{2e^{w(r)/2}}{1+e^{w(r)}}.
\]
Consequently, there exist $c,C>0$ and $R_0>1$ such that
\begin{equation}\label{BPS-critical-decay}
\frac{c}{r\log r}\leq S(r)\leq\frac{C}{r\log r}
\qquad (r\geq R_0).
\end{equation}
Combining the lower bounds in \eqref{BPS-z-log-lower} and
\eqref{BPS-critical-decay}, every nonzero regular solution of
\eqref{BPS-zq-system} satisfies
\[
|f(r)|=S(r)|z(r)|\geq\frac{c_1}{r}
\]
for all sufficiently large $r$.  This contradicts $f\in L^2((0,\infty),
r\,dr)$, since
\[
\int_{R_0}^\infty |f(r)|^2r\,dr
\geq c_1^2\int_{R_0}^\infty\frac{dr}{r}=\infty.
\]
Therefore $z_0$ cannot be nonzero, and the only element of
$\ker L_\theta^s$ is $f=0$.
\end{proof}

\begin{remark}
Proposition \ref{prop-BPS-nondegeneracy} concerns the radial polar
operator.  The full Hessian still has the phase zero mode
$L_\phi\sin\widetilde\theta=0$, which comes from rotations about the third
axis.  Moreover, the proposition asserts the absence of an $L^2$ zero
mode, not a positive spectral gap.  Since $\omega=0$, test functions
supported in annuli tending to infinity rule out a uniform estimate of
the form
\[
\langle L_\theta^sf,f\rangle_{L^2}\geq c\|f\|_{L^2}^2,
\qquad c>0,
\]
even after imposing one mass-orthogonality condition.
\end{remark}

\subsection{Constrained nondegeneracy in an interior near-BPS sector}

We next consider the standard potential
\[
V(t)=\kappa t^2,
\qquad 0<\kappa<1,
\]
when $(\kappa,\omega)$ approaches the self-dual endpoint $(1,0)$.
There is an important restriction on how this endpoint can be approached.
Put
\[
\delta=1-\kappa,
\qquad \lambda=\frac{\omega}{\delta}.
\]
For clarity, in this subsection $H_{\kappa,\omega}(\theta)$ denotes
the scalar functional \eqref{eq:JHS-scalar-functional} with
$V(t)=\kappa t^2$ and frequency $\omega$.  For an angle $\theta$, define
\begin{align*}
\mathsf M(\theta)&=2\pi\int_0^\infty
 r(1-\cos\theta)\,dr,\\
\mathsf P(\theta)&=2\pi\int_0^\infty
 r(1-\cos\theta)^2\,dr.
\end{align*}
We also use the Belavin--Polyakov radial energy
\begin{equation}\label{eq:BP-energy-E0}
E_0(w)=\pi\int_0^\infty
\left(|w'|^2+\frac{m^2}{r^2}\sin^2w\right)r\,dr.
\end{equation}
The factorization \eqref{BPS-factorization} gives the exact identity
\begin{equation}\label{near-BPS-decomposition}
H_{1-\delta,\delta\lambda}(\theta)-4\pi m
=\pi\int_0^\infty r\left|
\theta'+\frac{m-h_\theta}{r}\sin\theta\right|^2\,dr
+\delta\{\lambda\mathsf M(\theta)-\mathsf P(\theta)\}.
\end{equation}
Thus it is the ratio $\lambda=\omega/(1-\kappa)$, rather than the two
parameters separately, that selects the limiting BPS state.

We first record the radial BPS family and the two properties of this
family that will be used below.  The parameter $\nu$ is related to the
asymptotic slope \eqref{eq:BPS-slope-definition} by
$\beta_{\rm BPS}=2(1+\nu)$.

\begin{lemma}[The radial BPS family]\label{lem-BPS-family}
Let $m\geq2$.  There is a one-parameter family
\[
\{\Theta_\nu:0\leq\nu<m-1\}
\]
of radial BPS profiles satisfying
\begin{equation}\label{BPS-family-first-order}
\Theta_\nu'
+\frac{m-h_{\Theta_\nu}}r\sin\Theta_\nu=0,
\qquad
\Theta_\nu(0)=\pi,\qquad \Theta_\nu(\infty)=0.
\end{equation}
It exhausts the radial solutions of \eqref{BPS-family-first-order} in
the admissible class.  Moreover,
\begin{align}
h_{\Theta_\nu}(\infty)&=m-1-\nu,
                                                        \label{BPS-flux-nu}\\
\mathsf M(\Theta_\nu)&=2\pi(m-1-\nu),                  \label{BPS-mass-nu}\\
\sin\Theta_\nu(r)&=c_\nu r^{-1-\nu}(1+o(1))
       \quad (\nu>0),                                  \label{BPS-power-tail}
\end{align}
where $c_\nu>0$.  The member $\Theta_0$ is the critical BPS state in
Proposition \ref{prop-BPS-nondegeneracy}; in particular,
\[
\sin\Theta_0(r)\asymp\frac1{r\log r}.
\]
If
\[
P(\nu)=\mathsf P(\Theta_\nu),
\qquad P_*=P(0),\qquad M_*=2\pi(m-1),
\]
then $P$ is continuous on $[0,m-1]$, after setting $P(m-1)=0$, and
\begin{equation}\label{P-endpoint-small-o}
P(\nu)=P_*+o(\nu)\qquad(\nu\downarrow0).
\end{equation}
\end{lemma}

\begin{proof}
We include the endpoint calculation because \eqref{P-endpoint-small-o} is exactly what prevents the critical member
$\Theta_0$ from being selected when $\lambda>0$.

Set
\[
u=2\log\tan\frac{\Theta}{2}.
\]
Then \eqref{BPS-family-first-order}, together with
$h_\Theta'=r(1-\cos\Theta)$, is equivalent to
\begin{equation}\label{BPS-scalar-u}
u''+\frac1r u'=\frac{4e^u}{1+e^u},
\qquad
u(r)=-2m\log r+\alpha+o(1)\quad(r\downarrow0).
\end{equation}
For each initial constant $\alpha$, the local solution is equivalently
defined by the Volterra equation
\begin{equation}\label{BPS-u-Volterra}
u_\alpha(r)=-2m\log r+\alpha
+\int_0^r s\log\frac rs\,
       \frac{4e^{u_\alpha(s)}}{1+e^{u_\alpha(s)}}\,ds.
\end{equation}
When indexing by the shooting constant, write
$\Theta_\alpha(r)=2\arctan(e^{u_\alpha(r)/2})$; after labeling its
asymptotic slope by $\nu$, this is the same profile as $\Theta_\nu$.
It follows directly from \eqref{BPS-u-Volterra} that the solution and
its first derivative depend continuously differentiably on $\alpha$
on every compact subinterval of $(0,\infty)$.  Its variation
$z_\alpha=\partial_\alpha u_\alpha$ satisfies
\begin{equation}\label{BPS-u-variation}
z_\alpha''+\frac1r z_\alpha'
=\frac{4e^{u_\alpha}}{(1+e^{u_\alpha})^2}z_\alpha,
\qquad z_\alpha(0)=1,\qquad rz_\alpha'(r)\longrightarrow0.
\end{equation}
The Volterra form of \eqref{BPS-u-variation} shows that
$z_\alpha>0$ and $z_\alpha'>0$, so the shooting family is strictly
ordered.  All shooting solutions exist for every finite $r>0$:
the nonlinearity is bounded between $0$ and $4$, and the Volterra
formula gives
\begin{equation}\label{eq:shooting-elementary-bounds}
 -2m\log r+\alpha\le u_\alpha(r)
       \le-2m\log r+\alpha+r^2.
\end{equation}
Let $u_*$ be the critical radial solution supplied by
\cite[Theorem~1.1 and Lemma~3.3]{JHS}.  The regular origin expansion in
that lemma gives its finite shooting constant $\alpha_*$, so uniqueness
of the Volterra problem identifies $u_*=u_{\alpha_*}$.
Its known critical tail has slope $ru_*'\to-2$ and total slope integral
$2m-2$.

For $\alpha<\alpha_*$, ordering gives $u_\alpha<u_*$ and an integrable
majorant in the slope identity
\[
 ru_\alpha'(r)=-2m+\int_0^r\frac{4s e^{u_\alpha(s)}}{1+e^{u_\alpha(s)}}\,ds.
\]
Strict inequality on a compact annulus makes the limiting integral
strictly smaller than $2m-2$, and positivity makes it greater than zero.
Consequently the limiting slope is $-2(1+\nu)$ with
$0<\nu<m-1$.  Conversely, if $\alpha>\alpha_*$ then the lower limit
of that slope is greater than $-2$.  Such an orbit cannot tend to
$-\infty$: if its increasing slope ever becomes nonnegative the
conclusion is immediate, while a finite negative limiting slope
$L>-2$ would imply $u_\alpha(r)\ge C+(-2+\varepsilon)\log r$ for some
$\varepsilon>0$.  The slope integral would then diverge, a contradiction.
Thus the admissible shooting interval is exactly $(-\infty,\alpha_*]$.
Any admissible radial BPS solution has a finite origin constant because
$h=O(r^2)$ and $u'=-2(m-h)/r$; it is therefore one of these orbits.

The slope identity and strict ordering make $\nu(\alpha)$ strictly
decreasing.  On an interior parameter neighborhood, any fixed larger
admissible orbit supplies an integrable majorant, so dominated
convergence gives continuity.  As $\alpha\uparrow\alpha_*$, monotone
convergence gives $\nu(\alpha)\downarrow0$.
As $\alpha\to-\infty$, \eqref{eq:shooting-elementary-bounds} implies
$u_\alpha(r)\to-\infty$ for every fixed $r>0$.  Domination by one fixed
admissible orbit gives a slope integral tending to zero and hence
$\nu(\alpha)\uparrow m-1$.  These facts prove the asserted bijection,
uniqueness, and exhaustion without attributing radial uniqueness to a
more general nonradial existence theorem.

For completeness, the noncritical power tail follows from the limiting
slope.  First $e^{u_\alpha(r)}\le Cr^{-2-\nu}$ for large $r$.  Therefore
\[
 ru_\alpha'(r)+2(1+\nu)
  =-\int_r^\infty\frac{4s e^{u_\alpha(s)}}{1+e^{u_\alpha(s)}}\,ds
  =O(r^{-\nu}).
\]
Dividing by $r$ and integrating shows that
$u_\alpha(r)+2(1+\nu)\log r$ tends to a finite constant.  This proves
the power law for $\sin\Theta_\nu$.  The critical logarithmic
correction is the cited asymptotic for the critical orbit, not an
inference from monotonicity alone.

Since
\[
h_\Theta=m+\frac12ru',
\]
relations \eqref{BPS-flux-nu} and \eqref{BPS-mass-nu} follow at once.
The same identities and \eqref{BPS-scalar-u} give
\eqref{BPS-power-tail}.  The critical asymptotic is the one used in
\eqref{BPS-critical-decay}.

We now prove the endpoint estimate \eqref{P-endpoint-small-o}.
The argument only uses a one-sided quantitative comparison with the
critical orbit and avoids differentiating the scattering parameter at
infinity.

Write $\alpha(\nu)<\alpha_*$ for the shooting parameter that produces
$\Theta_\nu$.  The map $\alpha\mapsto\nu(\alpha)$ is strictly
decreasing on the admissible shooting interval.  Indeed, integrating
\eqref{BPS-scalar-u} gives
\begin{equation}\label{shooting-slope-integral}
-2\{1+\nu(\alpha)\}
=-2m+\int_0^\infty
 \frac{4r e^{u_\alpha(r)}}{1+e^{u_\alpha(r)}}\,dr,
\end{equation}
and $u_{\alpha_1}<u_{\alpha_2}$ pointwise whenever
$\alpha_1<\alpha_2$.  Hence the right-hand side of
\eqref{shooting-slope-integral} is strictly increasing in $\alpha$.
We justify differentiation on an interior shooting interval without
assuming differentiability of the asymptotic parameter.  Fix
$[\alpha^{-},\alpha^{+}]\Subset(-\infty,\alpha_*)$, and let
$\nu_0=\nu(\alpha^{+})>0$.  Ordering and the already established tail of
the single orbit $u_{\alpha^{+}}$ give, uniformly in this interval,
\[
 e^{u_\alpha(r)}\le Cr^{-2-2\nu_0}\quad(r\ge R_0).
\]
Put $Z_\alpha(t)=z_\alpha(e^t)$.  Its variation equation is
\[
 Z_\alpha''=\mathfrak a_\alpha(t)Z_\alpha,\qquad
 \mathfrak a_\alpha(t)=\frac{4e^{2t+u_\alpha(e^t)}}{(1+e^{u_\alpha(e^t)})^2}
 \le Ce^{-2\nu_0t}\quad(t\ge T_0).
\]
Values and first derivatives at the fixed $T_0$ are uniformly bounded
by local shooting dependence.  Integrating twice,
\[
 Z_\alpha(t)=Z_\alpha(T_0)+(t-T_0)Z_\alpha'(T_0)
       +\int_{T_0}^t(t-s)\mathfrak a_\alpha(s)Z_\alpha(s)\,ds.
\]
Divide by $1+t-T_0$ and use
$(t-s)/(1+t-T_0)\le1$.  Gronwall with the integrable coefficient
$(1+s-T_0)\mathfrak a_\alpha(s)$ gives
\[
 |Z_\alpha(t)|\le C(1+t-T_0),\qquad
 |Z_\alpha'(t)|\le C\quad(t\ge T_0),
\]
with constants uniform on the shooting interval.  Near zero the Volterra
equation gives a uniform bound for $z_\alpha$, and its variation
coefficient is $O(r^{2m})$.  The derivative of the integrand in
\eqref{shooting-slope-integral} is therefore dominated by
$Cr^{2m+1}$ at zero and by
$Cr^{-1-2\nu_0}(1+\log r)$ at infinity.  Both are integrable.
Dominated differentiation, and local continuous dependence of
$z_\alpha$, give a continuous derivative and
\eqref{BPS-u-variation} yields
\begin{equation}\label{eq:nu-alpha-derivative}
-2\nu'(\alpha)
=\int_0^\infty
\frac{4r e^{u_\alpha(r)}}{(1+e^{u_\alpha(r)})^2}
 z_\alpha(r)\,dr>0.
\end{equation}
Thus $\nu'(\alpha)<0$ and the inverse shooting parameter
$\alpha=\alpha(\nu)$ is $C^1$ on every compact subinterval of
$(0,m-1)$.  Since every fixed $\bar\alpha<\alpha_*$ corresponds to a
subcritical orbit with $\nu(\bar\alpha)>0$, it follows immediately that
\begin{equation}\label{alpha-to-critical}
\alpha(\nu)\uparrow\alpha_*\qquad(\nu\downarrow0).
\end{equation}
Put
\[
\varepsilon=\alpha_* - \alpha(\nu)>0.
\]
We first compare $\varepsilon$ with the asymptotic slope $\nu$.

Use again $t=\log r$ and set
\[
p_\alpha(t)=\frac12\{u_\alpha(e^t)+2t\},
\qquad a_\alpha(t)=-p_\alpha'(t).
\]
Equation \eqref{BPS-scalar-u} is then equivalent to the exact system
\begin{equation}\label{BPS-pa-system}
p_\alpha'=-a_\alpha,
\qquad
a_\alpha'=-\frac{2e^{2p_\alpha}}
 {1+e^{2p_\alpha-2t}}.
\end{equation}
In fact the identity $h_\Theta=m+ru'/2$ gives
\begin{equation}\label{a-as-flux-defect}
a_\alpha(t)=m-1-h_{\Theta_\alpha}(e^t).
\end{equation}
Thus $a_\alpha(t)>0$ on every admissible branch and
$a_\alpha(t)\downarrow\nu$ on the branch labelled by $\nu$.
Define
\[
J_\alpha(t)=a_\alpha(t)^2-2e^{2p_\alpha(t)}.
\]
A direct differentiation using \eqref{BPS-pa-system} gives
\begin{equation}\label{BPS-J-derivative}
J_\alpha'(t)
=\frac{4a_\alpha(t)e^{4p_\alpha(t)-2t}}
 {1+e^{2p_\alpha(t)-2t}}\geq0.
\end{equation}
Moreover,
\[
J_{\alpha(\nu)}(\infty)=\nu^2,
\qquad J_{\alpha_*}(\infty)=0.
\]
For the critical orbit, the asymptotic quoted before
\eqref{BPS-critical-decay} gives
$u_{\alpha_*}(e^t)=-2t-2\log t+O(1)$.  Hence
$p_{\alpha_*}(t)=-\log t+O(1)$.  In addition,
\eqref{a-as-flux-defect} and \eqref{BPS-critical-decay} give
\begin{equation}\label{critical-pa-rough}
p_{\alpha_*}(t)=-\log t+O(1),
\qquad
a_{\alpha_*}(t)\asymp\frac1t
\qquad(t\to\infty).
\end{equation}
Consequently, after increasing a fixed $T_0$ if necessary,
\begin{equation}\label{critical-a-lower}
a_{\alpha_*}(t)\geq \frac{c_0}{1+t},
\qquad t\geq T_0,
\end{equation}
and, since $J_{\alpha_*}(\infty)=0$,
\begin{equation}\label{critical-J-tail}
0\geq J_{\alpha_*}(T)
=-\int_T^\infty J_{\alpha_*}'(t)\,dt,
\qquad
|J_{\alpha_*}(T)|\leq C e^{-2T}(1+T)^{-5}.
\end{equation}
Here the last estimate follows directly from
\eqref{BPS-J-derivative} and \eqref{critical-pa-rough}.

Let
\[
Z_\beta(t):=\partial_\beta u_\beta(e^t).
\]
The variation equation \eqref{BPS-u-variation} becomes
\begin{equation}\label{BPS-Z-t-equation}
Z_\beta''(t)=c_\beta(t)Z_\beta(t),
\qquad
c_\beta(t)=
\frac{4e^{2p_\beta(t)}}
     {(1+e^{2p_\beta(t)-2t})^2}\geq0.
\end{equation}
Thus $Z_\beta'>0$ and is increasing.  Since
$Z_{\alpha_*}'(T_0)>0$, continuous dependence at the fixed time $T_0$
gives numbers $c_1>0$ and $\varepsilon_0>0$ such that
\begin{equation}\label{Zt-fixed-lower}
Z_\beta'(T_0)\geq c_1
\qquad
\text{for }\alpha_*-\varepsilon_0\leq\beta\leq\alpha_*.
\end{equation}
Hence $Z_\beta'(T)\geq c_1$ for every $T\geq T_0$ in this parameter
range.  Since $\partial_\beta a_\beta=-Z_\beta'/2$, for
$0<\varepsilon<\varepsilon_0$ we obtain
\begin{equation}\label{a-separation}
a_{\alpha_* -\varepsilon}(T)-a_{\alpha_*}(T)
=\frac12\int_{\alpha_* -\varepsilon}^{\alpha_*}
 Z_\beta'(T)\,d\beta
\geq \frac{c_1}{2}\varepsilon,
\qquad T\geq T_0.
\end{equation}
The shooting ordering also gives
$p_{\alpha_* -\varepsilon}\leq p_{\alpha_*}$ and
$a_{\alpha_* -\varepsilon}\geq a_{\alpha_*}$.  Therefore
\begin{align}
J_{\alpha_* -\varepsilon}(T)-J_{\alpha_*}(T)
&=(a_{\alpha_* -\varepsilon}-a_{\alpha_*})
  (a_{\alpha_* -\varepsilon}+a_{\alpha_*}) \notag\\
&\quad
 +2\bigl(e^{2p_{\alpha_*}}-e^{2p_{\alpha_* -\varepsilon}}\bigr)
\notag\\
&\geq \frac{c_2\varepsilon}{1+T},
\qquad T\geq T_0,
\label{J-separation}
\end{align}
where \eqref{critical-a-lower} was used in the last line.

Choose
\[
T_\varepsilon=\log\frac1\varepsilon.
\]
For sufficiently small $\varepsilon$, combining monotonicity of
$J_{\alpha(\nu)}$, \eqref{critical-J-tail}, and
\eqref{J-separation} yields
\begin{equation}\label{nu-epsilon-lower}
\nu^2=J_{\alpha(\nu)}(\infty)
\geq J_{\alpha(\nu)}(T_\varepsilon)
\geq \frac{c_3\varepsilon}{\log(1/\varepsilon)}.
\end{equation}
In particular,
\begin{equation}\label{epsilon-o-nu}
\frac{\varepsilon}{\nu}
\leq C\sqrt{\varepsilon\log(1/\varepsilon)}
\longrightarrow0.
\end{equation}
This weaker relation is exactly what is needed below; no differentiability
of $\nu$ with respect to the shooting parameter at the critical orbit is
required.

It remains to compare $P(\nu)$ with $P_*$.  In the $t$ variable,
\begin{equation}\label{P-t-variable}
P(\alpha)=8\pi\int_{-\infty}^{\infty}
 e^{2t}\Phi(u_\alpha(e^t))\,dt,
\qquad
\Phi(s)=\frac{e^{2s}}{(1+e^s)^2}.
\end{equation}
Because the shooting family is ordered,
$u_{\alpha(\nu)}\leq u_{\alpha_*}$.  Since $\Phi$ is increasing and
\begin{equation}\label{Phi-derivative-bound}
0\leq \Phi'(s)
=\frac{2e^{2s}}{(1+e^s)^3}
\leq2\Phi(s),
\end{equation}
we have $0\leq P_*-P(\nu)$.

We next record a uniform polynomial bound for the shooting variation.
For $\beta\leq\alpha_*$, the ordering gives $p_\beta\leq p_{\alpha_*}$,
and hence, by \eqref{critical-pa-rough},
\[
0\leq c_\beta(t)\leq \frac{C}{(1+t)^2},
\qquad t\geq T_0.
\]
An elementary comparison in \eqref{BPS-Z-t-equation} with
$(1+t)^K$, where $K>1$ is chosen so that $K(K-1)>C$, therefore gives
\begin{equation}\label{Z-polynomial-bound}
0<Z_\beta(t)\leq C(1+t)^K,
\qquad t\geq T_0,
\end{equation}
uniformly for
$\alpha_*-\varepsilon_0\leq\beta\leq\alpha_*$.  On $(-\infty,T_0]$ the Volterra equation
\eqref{BPS-u-variation}, together with
$Z_\beta(t)\to1$ and $Z_\beta'(t)\to0$ as $t\to-\infty$, gives the
uniform bound $0<Z_\beta(t)\le C$ on the indicated compact parameter
interval.  Combining the two regions, one has
$Z_\beta(t)\le C\{1+\max(t,0)\}^K$.  Thus, for
$t\leq T_\varepsilon$,
\begin{equation}\label{u-difference-log-window}
0\leq u_{\alpha_*}(e^t)-u_{\alpha(\nu)}(e^t)
=\int_{\alpha(\nu)}^{\alpha_*}Z_\beta(t)\,d\beta
\leq C\varepsilon(1+T_\varepsilon)^K.
\end{equation}
Using \eqref{Phi-derivative-bound} and
\eqref{u-difference-log-window} on $t\leq T_\varepsilon$ gives
\begin{align}
&\int_{-\infty}^{T_\varepsilon} e^{2t}
\left\{
\Phi(u_{\alpha_*}(e^t))-
\Phi(u_{\alpha(\nu)}(e^t))
\right\}\,dt
\notag\\
&\hspace{2cm}\leq
C\varepsilon(1+T_\varepsilon)^K
\int_{-\infty}^{\infty}e^{2t}
\Phi(u_{\alpha_*}(e^t))\,dt
\leq C\varepsilon\{\log(1/\varepsilon)\}^K.
\label{P-inner-difference}
\end{align}
At the critical orbit, \eqref{critical-pa-rough} implies
\[
e^{2t}\Phi(u_{\alpha_*}(e^t))
\leq C e^{-2t}(1+t)^{-4}
\qquad(t\geq T_0).
\]
Since the integrand for the subcritical orbit is smaller, the remaining
tail satisfies
\begin{equation}\label{P-tail-difference}
0\leq
\int_{T_\varepsilon}^{\infty}e^{2t}
\left\{
\Phi(u_{\alpha_*})-
\Phi(u_{\alpha(\nu)})
\right\}\,dt
\leq C e^{-2T_\varepsilon}
\leq C\varepsilon^2.
\end{equation}
Combining \eqref{P-inner-difference}, \eqref{P-tail-difference}, and
\eqref{nu-epsilon-lower}, we conclude that
\begin{align*}
0\leq\frac{P_*-P(\nu)}{\nu}
&\leq
C\frac{\varepsilon\{\log(1/\varepsilon)\}^K+\varepsilon^2}{\nu}\\
&\leq
C\sqrt\varepsilon\{\log(1/\varepsilon)\}^{K+1/2}
 +C\varepsilon^{3/2}\sqrt{\log(1/\varepsilon)}
\longrightarrow0.
\end{align*}
This proves \eqref{P-endpoint-small-o}.  Continuity for
$0<\nu<m-1$ follows from the shooting dependence and the uniform
power tail on compact $\nu$-intervals.  Continuity at $\nu=m-1$
follows from \eqref{BPS-mass-nu} and
$0\leq\mathsf P\leq2\mathsf M$.
\end{proof}

For $\lambda>0$, introduce the reduced BPS energy
\begin{equation}\label{reduced-BPS-energy}
g_\lambda(\nu)
=\lambda\mathsf M(\Theta_\nu)-\mathsf P(\Theta_\nu)
=2\pi\lambda(m-1-\nu)-P(\nu).
\end{equation}
The next lemma makes precise the selection of the limiting BPS state.

\begin{lemma}[Selection of an interior BPS limit]
\label{lem-interior-BPS-selection}
Set
\begin{equation}\label{lambda-star-definition}
\lambda_*:=\frac{P_*}{M_*}>0.
\end{equation}
Let $\mathcal K=[\lambda^{-},\lambda^{+}]$ be a compact subinterval of
$(0,\lambda_*)$.  If
\[
\delta_n\downarrow0,\qquad \lambda_n\in \mathcal K,
\]
and $\theta_n$ is a ground state for
$(\kappa_n,\omega_n)=(1-\delta_n,\delta_n\lambda_n)$, then, after
passing to a subsequence,
\begin{equation}\label{ground-state-strong-BPS-limit}
\lambda_n\longrightarrow\lambda\in \mathcal K,
\qquad
\theta_n\longrightarrow\Theta_\nu
\quad\hbox{strongly in }H^1_{\rm rad}(\mathbb R^2),
\end{equation}
where
\[
\nu\in\mathop{\rm argmin}_{0\leq\mu\leq m-1}g_\lambda(\mu).
\]
In fact, if $w_n=\theta_n-\Theta_\nu$, then
\begin{equation}\label{ground-state-angular-strong-limit}
E_0(w_n)\longrightarrow0,
\qquad
\|w_n\|_{L^\infty(0,\infty)}\longrightarrow0.
\end{equation}
There is a number $\nu_{\mathcal K}>0$, independent of the sequence, such that
every possible limiting parameter satisfies
\begin{equation}\label{nu-away-endpoints}
\nu_{\mathcal K}\leq\nu\leq m-1-\nu_{\mathcal K}.
\end{equation}
\end{lemma}

\begin{proof}
Define the compact-sector margin
\begin{equation}\label{eq:compact-sector-margin}
 c_{\mathcal K}:=P_*-\lambda^{+}M_*>0.
\end{equation}
Then
\[
 g_\lambda(0)=\lambda M_*-P_*\le-c_{\mathcal K}<0
       \qquad(\lambda\in \mathcal K).
\]
All uniform constants in this proof may depend on $\mathcal K$, in particular
on $c_{\mathcal K}$ and $\lambda^{-}>0$.  No uniform assertion is made as $\mathcal K$
approaches either endpoint of $(0,\lambda_*)$.
Thus \eqref{near-BPS-decomposition}, tested at $\Theta_0$, gives
\[
\inf H_{1-\delta,\delta\lambda}<4\pi m.
\]
The existence criterion in \cite[Theorem~1.3]{JHS} therefore
gives a ground state.  By the truncation argument
\cite[Lemma~3.1]{JHS}, it may be chosen with
$0\leq\theta_n\leq\pi$.

Write $E_{\rm BPS}=H_{1,0}$.  Since
$\mathsf P(\theta)\leq E_{\rm BPS}(\theta)$ and all the terms in
$H_{1-\delta,\delta\lambda}$ in its original form are nonnegative,
minimality and the preceding competitor give
\begin{align*}
(1-\delta_n)E_{\rm BPS}(\theta_n)
&\leq H_{1-\delta_n,\delta_n\lambda_n}(\theta_n)\leq C,\\
\lambda_n\mathsf M(\theta_n)
&\leq C+\mathsf P(\theta_n)\leq C.
\end{align*}
Hence $E_{\rm BPS}(\theta_n)$ and $\mathsf M(\theta_n)$ are uniformly
bounded.  Since $0\leq\theta_n\leq\pi$ and
$\theta_n^2\leq C(1-\cos\theta_n)$, the sequence is bounded in
$H^1_{\rm rad}(\mathbb R^2)$.  Moreover,
\begin{equation}\label{BPS-energy-convergence}
4\pi m\leq E_{\rm BPS}(\theta_n)
\leq H_{1-\delta_n,\delta_n\lambda_n}(\theta_n)
   +\delta_n\mathsf P(\theta_n)
=4\pi m+O(\delta_n).
\end{equation}

We also verify the angular integrability required by the splitting
result.  Choose a fixed $r_0>0$ with $r_0^2\le m/2$.  Since
$h_{\theta_n}(r)\le r^2$, the gauge part of the energy gives
\[
 \int_0^{r_0}\frac{\sin^2\theta_n}{r}\,dr
 \le C_m E_{\rm BPS}(\theta_n).
\]
On the complement,
\[
 \int_{r_0}^\infty\frac{\sin^2\theta_n}{r}\,dr
 \le r_0^{-2}\int_{r_0}^\infty\sin^2\theta_n\,r\,dr
 \le C_{r_0}\mathsf M(\theta_n).
\]
Thus these angular integrals are uniformly bounded.  After extraction,
let $\theta_n\rightharpoonup\theta$ in $H^1$, with almost-everywhere and
local strong convergence.  The radial bound
$|\theta_n(r)|\le C r^{-1/2}$ also gives, for each finite $p>2$,
\[
 \int_R^\infty|\theta_n|^p r\,dr
 \le C_p R^{-(p-2)/2}\int_R^\infty|\theta_n|^2r\,dr
 \le C_p R^{-(p-2)/2}.
\]
Together with local compactness this yields strong $L^p$ convergence for
$2<p<\infty$.  All assumptions of \cite[Proposition~3.5]{JHS} are
therefore satisfied.  If $\theta=0$, then
\cite[Lemma~3.4]{JHS}, applied to the potential $t^2$, gives
$\mathsf P(\theta_n)\to0$.  Consequently,
\[
\liminf_{n\to\infty}
\frac{H_{1-\delta_n,\delta_n\lambda_n}(\theta_n)-4\pi m}{\delta_n}
\geq0,
\]
because the square term and $\lambda_n\mathsf M(\theta_n)$ in
\eqref{near-BPS-decomposition} are nonnegative.  This contradicts the
uniformly negative value supplied by the competitor $\Theta_0$.
Therefore $\theta\ne0$.

Put $w_n=\theta_n-\theta$.  Proposition~3.5 of
\cite{JHS} gives
\begin{equation}\label{BPS-energy-splitting-nearby}
E_{\rm BPS}(\theta_n)
=E_{\rm BPS}(\theta)+E_0(w_n)+o(1),
\end{equation}
where $E_0$ is the Belavin--Polyakov energy defined in \eqref{eq:BP-energy-E0}, equivalently the functional in \cite[Lemma~3.2]{JHS}.  The same topological argument as in the
proof of \cite[Theorem~1.3]{JHS} shows that $\theta(0)$ is an
odd multiple of $\pi$: otherwise $w_n$ carries the degree and
$E_0(w_n)\geq4\pi m$, which is incompatible with
\eqref{BPS-energy-convergence} and $\theta\ne0$.  Hence
$E_{\rm BPS}(\theta)\geq4\pi m$.  Equations
\eqref{BPS-energy-convergence}--\eqref{BPS-energy-splitting-nearby}
then imply
\[
E_{\rm BPS}(\theta)=4\pi m,
\qquad E_0(w_n)\longrightarrow0.
\]
The square identity \eqref{BPS-factorization} shows that
$\theta=\Theta_\nu$ for some $0\leq\nu<m-1$.

Lemma~3.4 of \cite{JHS}, now applied respectively to $t^2$ and
$t$, gives
\begin{align}
\mathsf P(\theta_n)&=\mathsf P(\Theta_\nu)+o(1),
                                                        \label{P-strong-limit}\\
\mathsf M(\theta_n)&=\mathsf M(\Theta_\nu)
  +\mathsf M(w_n)+o(1).                                \label{M-splitting-limit}
\end{align}
For every BPS competitor $\Theta_\mu$ with $0\le\mu<m-1$, minimality gives
\[
\frac{H_{1-\delta_n,\delta_n\lambda_n}(\theta_n)-4\pi m}{\delta_n}
\leq g_{\lambda_n}(\mu).
\]
Using \eqref{near-BPS-decomposition},
\eqref{P-strong-limit}, and \eqref{M-splitting-limit}, we conclude
\[
g_\lambda(\nu)
\leq\liminf_n\frac{H_{1-\delta_n,\delta_n\lambda_n}(\theta_n)-4\pi m}
 {\delta_n}
\leq\min_\mu g_\lambda(\mu).
\]
Passing to $\mu\uparrow m-1$ by the continuous extension of $g_\lambda$ shows that the same comparison holds on the closed interval $[0,m-1]$.  Thus $\nu$ is a minimizer.  Equality also forces
$\mathsf M(w_n)\to0$.  Since $w_n^2\leq C(1-\cos w_n)$ and
$E_0(w_n)\to0$, this proves the strong convergence in
\eqref{ground-state-strong-BPS-limit}.  We also retain the angular part
of the latter convergence.  Since $w_n(0)=0$ and $|w_n|\leq\pi$,
\[
1-\cos w_n(r)=\int_0^r \sin w_n(s)w_n'(s)\,ds,
\]
and therefore
\[
\sup_{r>0}\{1-\cos w_n(r)\}
\leq
\left(\int_0^\infty\frac{\sin^2 w_n}{r}\,dr\right)^{1/2}
\left(\int_0^\infty |w_n'|^2r\,dr\right)^{1/2}
\longrightarrow0.
\]
On $[-\pi,\pi]$, $1-\cos t\geq2t^2/\pi^2$, and hence
\eqref{ground-state-angular-strong-limit} follows.

It remains to locate the minimizers.  At the concentrating endpoint
$\nu=m-1$, one has $g_\lambda(m-1)=0$, whereas
$g_\lambda(0)\le-c_{\mathcal K}$; hence a minimizer stays uniformly away from
$m-1$.  At the critical endpoint, \eqref{P-endpoint-small-o} gives
\[
g_\lambda(\nu)-g_\lambda(0)
=-2\pi\lambda\nu+P_*-P(\nu)<0
\]
for all sufficiently small $\nu>0$, uniformly for $\lambda\in \mathcal K$.
Compactness of $\mathcal K$ and continuity of $g_\lambda$ now give
\eqref{nu-away-endpoints}.
\end{proof}

We next identify the kernel at an interior BPS member.  Unlike the
critical state, an interior member has a one-dimensional $L^2$ polar
kernel, but this kernel is transverse to the mass constraint.

\begin{lemma}[Kernel at an interior BPS state]
\label{lem-interior-BPS-kernel}
For $0<\nu<m-1$,
\begin{equation}\label{interior-BPS-kernel}
\ker L_{\Theta_\nu}^s
=\operatorname{span}\{\partial_\nu\Theta_\nu\}.
\end{equation}
Moreover,
\begin{equation}\label{BPS-kernel-mass-pairing}
\left\langle\partial_\nu\Theta_\nu,
\sin\Theta_\nu\right\rangle_{L^2}=-2\pi.
\end{equation}
In particular,
\[
\ker L_{\Theta_\nu}^s\cap
\{f:\langle f,\sin\Theta_\nu\rangle_{L^2}=0\}=\{0\}.
\]
\end{lemma}

\begin{proof}
The factorization \eqref{BPS-Hessian-square} and the reduction
\eqref{BPS-zq-system} hold for every BPS member.  The regularity at the
origin again leaves only one solution, up to multiplication by a
constant.  Normalize the nonzero regular solution by $z(0)>0$.  Then
\eqref{BPS-zq-system} gives
\[
 q_f'(r)=r\sin^2\Theta_\nu(r)z(r)>0,
 \qquad z'(r)=q_f(r)/r>0,
\]
so $q_f(r)>0$ for $r>0$.  Put $t=\log r$ and
$a(t)=e^{2t}\sin^2\Theta_\nu(e^t)$.  By
\eqref{BPS-power-tail},
$a(t)\le Ce^{-2\nu t}$ for large $t$.  The system becomes
\[
 z_t=q_f,
 \qquad \partial_t q_f=a(t)z,
\]
and hence, for $t\ge T$,
\[
 z(t)=z(T)+(t-T)q_f(T)+\int_T^t(t-s)a(s)z(s)\,ds.
\]
Since $\int_T^\infty(1+s)a(s)\,ds<\infty$, Volterra--Gronwall gives
$|z(t)|\le C(1+t)$.  Therefore $\partial_t q_f=a(t)z(t)$ is integrable,
$q_f(t)$ converges to a finite limit $q_\infty$, and monotonicity gives
$q_\infty>0$.  A second integration, using
$\int_T^\infty s|\partial_t q_f(s)|\,ds<\infty$, yields
\[
 q_f(r)\longrightarrow q_\infty\ne0,
 \qquad z(r)=q_\infty\log r+O(1).
\]
Consequently
\[
f(r)=\sin\Theta_\nu(r)z(r)
=O(r^{-1-\nu}\log r),
\]
which belongs to $L^2((0,\infty),r\,dr)$.  The first-order system also
gives $f'=O(r^{-2-\nu}\log r)$ at infinity.  At zero its regular solution
has $f=O(r^m)$ and $f'=O(r^{m-1})$.  Consequently $f\in\mathcal H$.
Moreover $q_f=\mathcal I_{\Theta_\nu}f$ and
$\mathcal B_{\Theta_\nu}f=0$.  The polarized square identity, extended
by Lemma~\ref{lem:BPS-form-realization}, gives
\[
 \mathfrak q_{\Theta_\nu}[f,g]
 =2\pi\int_0^\infty
      (\mathcal B_{\Theta_\nu}f)(\mathcal B_{\Theta_\nu}g)r\,dr=0
 \quad(g\in\mathcal H).
\]
The representation theorem therefore implies
$f\in D(L_{\Theta_\nu}^s)$ and $L_{\Theta_\nu}^sf=0$.
Conversely every operator-kernel element satisfies this first-order
system by the square identity.  Its regular solution space at zero is
one-dimensional, so the operator kernel is exactly one-dimensional.

By \eqref{eq:nu-alpha-derivative}, the interior shooting branch is $C^1$
in $\nu$.  The relation
$\partial_\alpha\Theta_\alpha=(\sin\Theta_\alpha)z_\alpha/2$,
together with the uniform variation bounds used to prove that equation,
shows that $\partial_\nu\Theta_\nu\in\mathcal H$ locally uniformly in
the interior parameter.  Differentiating the BPS equation gives a
nonzero regular solution of the same first-order system.  The preceding
form-domain argument places it in the operator kernel, proving
\eqref{interior-BPS-kernel}.

Finally, on a compact interior $\nu$-interval the same bounds dominate
$r\sin\Theta_\nu\,\partial_\nu\Theta_\nu$ at infinity by
$Cr^{-1-2\nu_0}(1+\log r)$ for some $\nu_0>0$, and at zero by an
integrable power.  Dominated differentiation of the mass formula
\eqref{BPS-mass-nu} is therefore justified:
\[
-2\pi=\frac d{d\nu}\mathsf M(\Theta_\nu)
=2\pi\int_0^\infty
r\sin\Theta_\nu\,\partial_\nu\Theta_\nu\,dr.
\]
This is \eqref{BPS-kernel-mass-pairing}.
\end{proof}

The following compactness statement is the point at which the strict
inequality $\nu>0$ is essential.

\begin{lemma}[No escape of normalized zero modes]
\label{lem-zero-mode-no-escape}
Let $\mathcal K\Subset(0,\lambda_*)$, let $\delta_n\downarrow0$ and
$\lambda_n\in \mathcal K$, and let $\theta_n$ be ground states for
\[
\kappa_n=1-\delta_n,
\qquad
\omega_n=\delta_n\lambda_n.
\]
Assume that, after extraction, the convergence in
Lemma~\ref{lem-interior-BPS-selection} is
\[
\theta_n\longrightarrow\Theta_\nu
\quad\hbox{strongly in }H^1_{\rm rad}(\mathbb R^2)
\]
for some $\nu>0$.  Put $L_n=L_{\theta_n}^s$.  If
\[
L_nf_n=0,
\qquad
\|f_n\|_{L^2}=1,
\]
then $\{f_n\}$ is relatively compact in radial $L^2$.
\end{lemma}

\begin{proof}
We use throughout the measure $r\,dr$; the harmless factor $2\pi$ is
suppressed.  Set
\[
S_n=\sin\theta_n,
\qquad C_n=\cos\theta_n,
\qquad a_n=h_{\theta_n}-m,
\]
and
\[
H_n=h_{\theta_n}(\infty),
\qquad
p_n=m-H_n,
\qquad
d_n(r)=H_n-h_{\theta_n}(r)\geq0.
\]
Thus
\begin{equation}\label{pn-dn-relation}
m-h_{\theta_n}(r)=p_n+d_n(r).
\end{equation}
The proof of Lemma~\ref{lem-interior-BPS-selection} gives in addition
\eqref{ground-state-angular-strong-limit}.  In particular,
$\theta_n\to\Theta_\nu$ uniformly on $(0,\infty)$.

\smallskip
\noindent\emph{Step 1: a uniform form bound for the zero modes.}
Strong $H^1$ convergence implies convergence of the masses.  Indeed,
for $0\leq s,t\leq\pi$,
\[
|\cos s-\cos t|\leq\frac12(|s|+|t|)|s-t|,
\]
so Cauchy--Schwarz gives
\begin{equation}\label{pn-convergence-rigorous}
p_n\longrightarrow m-h_{\Theta_\nu}(\infty)=1+\nu.
\end{equation}
Moreover, from $\theta_n=\Theta_\nu+w_n$ and
$|\sin(\alpha+\beta)|\leq|\sin\alpha|+|\sin\beta|$,
\eqref{ground-state-angular-strong-limit} gives
\begin{equation}\label{uniform-angular-profile-bound}
\sup_n\int_0^\infty\frac{S_n(r)^2}{r}\,dr<\infty.
\end{equation}
Together with $h_{\theta_n}(r)\leq r^2$ and the uniform mass bound,
this implies
\begin{equation}\label{uniform-k-bound-noescape}
\sup_n\|k_{\theta_n}\|_{L^\infty}<\infty.
\end{equation}

For a test function $f\in\mathcal H$, put
\[
\mathcal H_{n,f}(r)=\int_0^r sS_n(s)f(s)\,ds.
\]
Here $\mathcal H_{n,f}=\mathcal I_{\theta_n}f$.
For the polar section with radial coefficient $f$ the variation is
$\delta h[f e_{\theta_n}^0]=-\mathcal I_{\theta_n}f$
(the scalar angle variation has the opposite sign).  Thus the two
nonlocal pieces of
$L_{\theta_n}$ are
\begin{align}
(\mathcal N_{1,n}f)(r)
&=S_n(r)\int_r^\infty\frac{S_n(\rho)^2}{\rho}
       \mathcal H_{n,f}(\rho)\,d\rho,                 \label{N1-noescape}\\
(\mathcal N_{2,n}f)(r)
&=\frac{4S_n(r)C_n(r)a_n(r)}{r^2}\mathcal H_{n,f}(r).
                                                               \label{N2-noescape}
\end{align}
The nonlocal part of the symmetrized operator is
$\mathcal N_{1,n}+\frac12(\mathcal N_{2,n}+\mathcal N_{2,n}^*)$.
Notice that
\begin{equation}\label{Hnf-pointwise}
|\mathcal H_{n,f}(r)|
\leq
\left(\int_0^r\frac{|f(s)|^2}{s}\,ds\right)^{1/2}
\left(\int_0^r s^3S_n(s)^2\,ds\right)^{1/2}
\leq Cr^2\|f\|_{\mathcal H}.
\end{equation}
The $r^2$ bound is valid for all $r$ and is the one needed near zero.
On $r\ge1$ it can also be sharpened to
$|\mathcal H_{n,f}(r)|\le Cr\|f\|_{\mathcal H}$, since
$\int_0^r s^3S_n^2\,ds\le r^2\|S_n\|_{L^2(rdr)}^2$.
Since $\sup_n\|S_n\|_{L^2(rdr)}<\infty$, formulas
\eqref{N1-noescape}--\eqref{Hnf-pointwise} imply
\[
\|\mathcal N_{1,n}f\|_{L^2}
 +\|\mathcal N_{2,n}f\|_{L^2}
\leq C\|f\|_{\mathcal H}.
\]
Consequently, for every $\varepsilon>0$,
\begin{equation}\label{nonlocal-form-small}
\left|
\left\langle
\left(\mathcal N_{1,n}+\frac12(\mathcal N_{2,n}+\mathcal N_{2,n}^*)
\right)f,f\right\rangle
\right|
\leq
\varepsilon\left(\|f'\|_2^2+m^2\|f/r\|_2^2\right)
+C_\varepsilon\|f\|_2^2,
\end{equation}
uniformly in $n$.

We next treat the local part.  The only term which is not immediately a
uniformly bounded multiplier after subtracting $m^2/r^2$ is
\[
\frac{a_n^2\cos(2\theta_n)-m^2}{r^2}
=\frac{h_{\theta_n}^2-2mh_{\theta_n}}{r^2}
 -\frac{2a_n^2S_n^2}{r^2}.
\]
The first term on the right is uniformly bounded because
$h_{\theta_n}(r)\leq r^2$ for $r\leq1$ and $H_n$ is uniformly bounded.
By the uniform convergence in \eqref{ground-state-angular-strong-limit}
and $\Theta_\nu(0)=\pi$, we may choose $r_0>0$ so small that, after
discarding finitely many indices,
\[
\sup_{0<r<r_0}S_n(r)^2\leq\varepsilon.
\]
Hence the second term is bounded below in quadratic-form sense by
\[
-C\varepsilon\int_0^{r_0}\frac{f^2}{r^2}\,r\,dr
-C_{\varepsilon,r_0}\|f\|_2^2.
\]
All remaining local coefficients are uniformly bounded by
\eqref{uniform-k-bound-noescape}, $0\leq\kappa_n\leq1$, and
$0<\omega_n\leq1$.  Choosing first $r_0$ and then the $\varepsilon$ in
\eqref{nonlocal-form-small} sufficiently small, the identity
$\langle L_nf_n,f_n\rangle=0$ and $\|f_n\|_2=1$ give
\begin{equation}\label{uniform-kernel-H-bound}
\sup_n\|f_n\|_{\mathcal H}<\infty.
\end{equation}

\smallskip
\noindent\emph{Step 2: a uniform polynomial tail for the profiles.}
The radial lemma and the uniform $H^1$ bound for $\theta_n$ give
\[
\sup_n\sup_{r\geq R}|\theta_n(r)|\longrightarrow0
\qquad(R\to\infty).
\]
Moreover,
\begin{equation}\label{uniform-dn-tail}
\sup_n d_n(R)
\leq C\sup_n\int_R^\infty \theta_n(r)^2r\,dr
\longrightarrow0.
\end{equation}
Since $S_n^2\leq2(1-C_n)$ and $|a_n|\leq C$,
\begin{equation}\label{uniform-k-tail-preliminary}
r^2|k_{\theta_n}(r)|
\leq C d_n(r),
\qquad r>0.
\end{equation}

Fix $0<\eta<\nu/2$.  The profile equation can be written on the
exterior region as
\begin{equation}\label{profile-linear-comparison-equation}
-\theta_n''-\frac1r\theta_n'+c_n(r)\theta_n=0,
\end{equation}
where
\[
c_n(r)=
\left[
 k_{\theta_n}(r)+\omega_n
 +\frac{(m-h_{\theta_n}(r))^2}{r^2}C_n(r)
 +2\kappa_n(1-C_n(r))
\right]\frac{S_n(r)}{\theta_n(r)}.
\]
Here $0<\theta_n<\pi$ by the ground-state maximum principle.  By
\eqref{pn-convergence-rigorous}, \eqref{uniform-dn-tail}, and
\eqref{uniform-k-tail-preliminary}, there exists $R_1$ such that
\begin{equation}\label{profile-exterior-indicial-gap}
c_n(r)\geq\frac{(1+\eta)^2}{r^2},
\qquad r\geq R_1,
\end{equation}
for all sufficiently large $n$.  The function $r^{-1-\eta}$ is therefore
a supersolution of \eqref{profile-linear-comparison-equation}.  Comparing
on $(R_1,T)$ and then sending $T\to\infty$ (each $\theta_n$ has the
exponential tail supplied by $\omega_n>0$) yields
\begin{equation}\label{uniform-ground-state-power-tail}
0<\theta_n(r)\leq Cr^{-1-\eta},
\qquad r\geq R_1,
\end{equation}
uniformly in $n$.  Consequently,
\begin{equation}\label{uniform-hk-tail}
d_n(r)\leq Cr^{-2\eta},
\qquad
|k_{\theta_n}(r)|\leq Cr^{-2-2\eta}.
\end{equation}

\smallskip
\noindent\emph{Step 3: the exact exterior Hessian and its weighted
remainder.}
On $(R_1,\infty)$ the symmetrized polar Hessian has the exact form
\begin{equation}\label{exact-exterior-hessian-decomposition}
L_n=A_{p_n,\omega_n}+W_n
 +\mathcal N_{1,n}
 +\frac12(\mathcal N_{2,n}+\mathcal N_{2,n}^*),
\end{equation}
where
\[
A_{p,\omega}
=-\partial_r^2-\frac1r\partial_r+\frac{p^2}{r^2}+\omega
\]
and, using $V(t)=\kappa_nt^2$,
\begin{align*}
W_n={}&k_{\theta_n}C_n
 +\frac{a_n^2\cos(2\theta_n)-p_n^2}{r^2}\\
&+2\kappa_n\{S_n^2+(1-C_n)C_n\}
 +\omega_n(C_n-1).
\end{align*}
Equations \eqref{pn-dn-relation},
\eqref{uniform-ground-state-power-tail}, and
\eqref{uniform-hk-tail} give
\begin{equation}\label{local-remainder-tail}
|W_n(r)|\leq Cr^{-2-2\eta},
\qquad r\geq R_1.
\end{equation}

Let $0<\gamma<\eta$, let $R\geq R_1$, and abbreviate
\[
X_R(f)=\|r^\gamma f\|_{L^2(R,\infty;rdr)},
\qquad
Y_R(f)=\|f\|_{L^2(0,R;rdr)}.
\]
For $r\geq R$, Cauchy--Schwarz and
\eqref{uniform-ground-state-power-tail} give
\begin{equation}\label{Hnf-weighted-split}
|\mathcal H_{n,f}(r)|
\leq C_RY_R(f)+CR^{-\gamma-\eta}X_R(f).
\end{equation}
Using \eqref{N1-noescape}--\eqref{N2-noescape}, one obtains directly
\begin{align}
\|r^{\gamma+2}\mathcal N_{2,n}f\|_{L^2(R,\infty;rdr)}
&\leq CR^{-2\eta}X_R(f)+C_RY_R(f),                 \label{N2-weighted}\\
\|r^{\gamma+2}\mathcal N_{2,n}^*f\|_{L^2(R,\infty;rdr)}
&\leq CR^{-2\eta}X_R(f),                           \label{N2star-weighted}\\
\|r^{\gamma+2}\mathcal N_{1,n}f\|_{L^2(R,\infty;rdr)}
&\leq CR^{-4\eta}X_R(f)+C_RY_R(f).                 \label{N1-weighted}
\end{align}
For completeness, the first estimate follows from
$|S_n(r)|\leq Cr^{-1-\eta}$ and
\[
\|r^{\gamma-1-\eta}\|_{L^2(R,\infty;rdr)}
= C_{\gamma,\eta}R^{\gamma-\eta};
\]
the exterior part in \eqref{Hnf-weighted-split} therefore contributes
$CR^{-2\eta}X_R(f)$.  For the adjoint term we use
\[
(\mathcal N_{2,n}^*f)(r)
=4S_n(r)\int_r^\infty
 \frac{S_n(s)C_n(s)a_n(s)}s f(s)\,ds
\]
and weighted Cauchy--Schwarz.  Finally,
\[
\int_r^\infty\frac{S_n(\rho)^2}{\rho}\,d\rho
\leq Cr^{-2-2\eta}
\]
gives \eqref{N1-weighted}.  Combining these estimates with
\eqref{local-remainder-tail}, we obtain
\begin{equation}\label{weighted-remainder-estimate}
\|r^{\gamma+2}\mathcal R_nf\|_{L^2(R,\infty;rdr)}
\leq\varepsilon_R X_R(f)+C_RY_R(f),
\qquad
\varepsilon_R\longrightarrow0,
\end{equation}
where
\[
\mathcal R_n
=W_n+\mathcal N_{1,n}
 +\frac12(\mathcal N_{2,n}+\mathcal N_{2,n}^*).
\]
All constants in \eqref{weighted-remainder-estimate} are uniform in $n$.
For later use put $\varpi_L(r)=\min\{r,L\}^\gamma$, $L\ge R$, and
$X_{R,L}(f)=\|\varpi_L f\|_{L^2(R,\infty;rdr)}$.  The same remainder
estimate holds with $r^{\gamma+2}$ replaced by $r^2\varpi_L$ and
$X_R$ by $X_{R,L}$, uniformly in $L$.  Indeed, for $s\ge r\ge R$,
\[
 \varpi_L(R)=R^\gamma,\qquad
 \varpi_L(r)\le r^\gamma,\qquad
 1\le\frac{\varpi_L(s)}{\varpi_L(r)}\le(s/r)^\gamma.
\]
The first two properties give the same forward-integral bound as
\eqref{Hnf-weighted-split} and the same estimates for $\mathcal N_{1,n}$
and $\mathcal N_{2,n}$.  For the adjoint term, monotonicity gives
\[
 \left|\int_r^\infty
 \frac{S_n(s)C_n(s)a_n(s)}s f(s)\,ds\right|
 \le C\,\varpi_L(r)^{-1}r^{-2-\eta}X_{R,L}(f).
\]
Multiplication by $r^2\varpi_L(r)S_n(r)$ and integration again yield
$CR^{-2\eta}X_{R,L}(f)$.  The multiplier $W_n$ has the same bound
pointwise.  This proves the asserted truncated estimate without assuming
$X_R(f)<\infty$.

\smallskip
\noindent\emph{Step 4: a uniform exterior inverse-square estimate.}
We record an elementary estimate which is uniform as $\omega\downarrow0$.
Suppose $p\geq p_0>1$, $0<\gamma<p_0-1$, and
\[
A_{p,\omega}v=F
\quad\hbox{on }(R,\infty),
\qquad 0<\omega\leq1,
\]
with $v\in L^2(R,\infty;rdr)$.  Then
\begin{equation}\label{uniform-Bessel-estimate}
\|r^\gamma v\|_{L^2(R,\infty;rdr)}
\leq C_R|v(R)|
 +C\|r^{\gamma+2}F\|_{L^2(R,\infty;rdr)},
\end{equation}
where $C$ is uniform for $p\geq p_0$ in a fixed compact set and
$0<\omega\leq1$.

Here is a truncated-weight proof which does not assume the desired
weighted integrability in advance.  Let $\phi$ be the positive decaying
solution of $A_{p,\omega}\phi=0$ with $\phi(R)=1$.  The maximum principle
and $A_{p,\omega}(R/r)^p=\omega(R/r)^p\ge0$ give
\[
 0<\phi(r)\le(R/r)^p,
 \qquad \|r^\gamma\phi\|_2\le C_R,
\]
uniformly for the stated $p,\omega$, since $\gamma<p_0-1$.
Set $w=v-v(R)\phi$, so $w(R)=0$.  For $L\ge R$, define
\[
 \varpi_L(r)=\min\{r,L\}^\gamma,
 \qquad \mu_L(r)=r\varpi_L(r).
\]
This locally Lipschitz weight satisfies, almost everywhere,
\[
 0<\frac{r\mu_L'(r)}{\mu_L(r)}\le1+\gamma.
\]
Testing $A_{p,\omega}w=F$ with $\mu_L^2w$ and completing the square in
the derivative gives
\begin{align*}
 &\int_R^\infty |(\mu_Lw)'|^2r\,dr
 +\int_R^\infty
 \left\{\frac{p^2}{r^2}+\omega-
              \left(\frac{\mu_L'}{\mu_L}\right)^2\right\}
                |\mu_Lw|^2r\,dr\\
 &\hspace{25mm}=\operatorname{Re}\int_R^\infty
                F\overline w\,\mu_L^2r\,dr.
\end{align*}
The resulting coercive coefficient is bounded below by
\[
 c_0:=p_0^2-(1+\gamma)^2>0.
\]
Therefore
\begin{equation}\label{eq:truncated-inverse-square}
 c_0\|\varpi_L w\|_2
 \le \|r^2\varpi_L F\|_2.
\end{equation}

For full justification of the test, first insert a nonincreasing cutoff
$\chi_N$ equal to one on $[R,N]$, zero on $[2N,\infty)$, with
$|\chi_N'|\le C/N$, and take $N>\max\{L,R\}$.  The same square
identity uses $\chi_N\mu_L$.  Since $\mu_L'\ge0$ and $\chi_N'\le0$,
its cross term has the favorable sign; the only error to subtract is
bounded by
\[
 \int_N^{2N}|\chi_N'|^2\mu_L^2|w|^2r\,dr
 \le C L^{2\gamma}\int_N^{2N}|w|^2r\,dr\longrightarrow0
 \quad(N\to\infty)
\]
for each fixed $L$.  The nonnegative derivative and $\omega$ terms
may be discarded.  This proves \eqref{eq:truncated-inverse-square}
whenever the truncated right side is finite.  There is no boundary term
at $r=L$: the weight is locally Lipschitz and only its first weak
derivative is used.  Adding back $v(R)\phi$ proves
\[
 \|\varpi_L v\|_2\le C_R|v(R)|+c_0^{-1}\|r^2\varpi_L F\|_2
\]
with constants independent of $L$ and $\omega$.  Monotone convergence as
$L\to\infty$ gives \eqref{uniform-Bessel-estimate}.

By \eqref{pn-convergence-rigorous}, after increasing $R$ and discarding
finitely many indices, $p_n\geq p_0>1+\gamma$.  Apply
\eqref{uniform-Bessel-estimate} to the equation
\[
A_{p_n,\omega_n}f_n=-\mathcal R_nf_n.
\]
Using \eqref{eq:truncated-inverse-square} and the truncated remainder
estimate, we obtain, for every $L\ge R$,
\[
 X_{R,L}(f_n)\le C_R|f_n(R)|
       +C\varepsilon_R X_{R,L}(f_n)+C_RY_R(f_n).
\]
Choose $R$ so large that $C\varepsilon_R\le1/2$.  The trace theorem on
the fixed annulus $(R/2,2R)$ and \eqref{uniform-kernel-H-bound} give
$\sup_n|f_n(R)|<\infty$, while $Y_R(f_n)\le1$.  Absorption therefore
bounds $X_{R,L}(f_n)$ uniformly in both $n$ and $L$.  Sending $L\to\infty$
by monotone convergence gives
\begin{equation}\label{zero-mode-weighted-bound}
\sup_n\int_R^\infty r^{2\gamma}|f_n(r)|^2r\,dr<\infty.
\end{equation}

Finally, \eqref{uniform-kernel-H-bound} gives Rellich compactness on every
bounded disk, while \eqref{zero-mode-weighted-bound} gives the uniform
tightness
\[
\sup_n\int_A^\infty |f_n(r)|^2r\,dr
\leq CA^{-2\gamma}\longrightarrow0
\qquad(A\to\infty).
\]
A diagonal extraction therefore converges strongly in the full radial
$L^2$ space.  This proves the lemma.
\end{proof}

We can now state the near-BPS nondegeneracy result.

\begin{prop}[Near-BPS constrained nondegeneracy]
\label{prop-near-BPS-nondegeneracy}
Let $m\geq2$, and let $\mathcal K=[\lambda^{-},\lambda^{+}]$ be a compact
subinterval of $(0,\lambda_*)$, where $\lambda_*$ is given by
\eqref{lambda-star-definition}.  There exists $\delta_{\mathcal K}>0$ such that,
for
\[
0<\delta<\delta_{\mathcal K},
\qquad \kappa=1-\delta,
\qquad \omega=\delta\lambda,
\qquad \lambda\in \mathcal K,
\]
every radial ground state $\widetilde\theta_{\delta,\lambda}$ satisfies
\begin{equation}\label{near-BPS-ND}
\ker L_{\widetilde\theta_{\delta,\lambda}}^s
\cap\left\{f\in\mathcal H:
\langle f,\sin\widetilde\theta_{\delta,\lambda}\rangle_{L^2}=0
\right\}=\{0\}.
\end{equation}
Consequently the orbital stability theorem applies to these ground
states.
\end{prop}

\begin{proof}
Suppose otherwise.  Then there are sequences
$\delta_n\downarrow0$, $\lambda_n\in \mathcal K$, ground states $\theta_n$,
and functions $f_n$ such that
\begin{equation}\label{degenerate-kernel-sequence}
L_{\theta_n}^sf_n=0,\qquad
\|f_n\|_{L^2}=1,\qquad
\langle f_n,\sin\theta_n\rangle_{L^2}=0.
\end{equation}
Lemma \ref{lem-interior-BPS-selection} gives, after extraction,
\[
\theta_n\longrightarrow\Theta_\nu
\quad\hbox{strongly in }H^1_{\rm rad}(\mathbb R^2)
\]
for a parameter $\nu\geq\nu_{\mathcal K}>0$.  Lemma
\ref{lem-zero-mode-no-escape} therefore gives, after one more
extraction,
\[
f_n\longrightarrow f
\quad\hbox{strongly in }L^2,
\qquad \|f\|_{L^2}=1.
\]
The bound \eqref{uniform-kernel-H-bound} gives weak $H^1$
convergence.  The profile ODE gives local $C^2$ convergence away from the
origin.  Moreover, \eqref{ground-state-angular-strong-limit} implies
$S_n\to\sin\Theta_\nu$ both in $L^2(r\,dr)$ and in the angular space
$L^2(r^{-1}dr)$.  For the latter, write
\[
S_n-\sin\Theta_\nu
=\sin\Theta_\nu(\cos w_n-1)+\cos\Theta_\nu\sin w_n
\]
and use the two terms in $E_0(w_n)$ together with
$\|w_n\|_\infty\to0$.  Together with the strong $L^2$ convergence of $f_n$,
the explicit formulas \eqref{N1-noescape}--\eqref{N2-noescape} then allow
passage to the limit in the nonlocal bilinear forms by Cauchy--Schwarz.
To make the weak-limit step explicit, the first equation of
\eqref{degenerate-kernel-sequence} means
\begin{equation}\label{eq:weak-zero-mode-approximants}
 \mathfrak q_{\theta_n}[f_n,\varphi]=0
 \qquad\text{for every }\varphi\in C_c^\infty(0,\infty).
\end{equation}
The convergences above and the explicit nonlocal kernels yield
\begin{equation}\label{eq:weak-zero-mode-limit}
 \mathfrak q_{\Theta_\nu}[f,\varphi]=0
 \qquad\text{for every }\varphi\in C_c^\infty(0,\infty).
\end{equation}
As $f\in\mathcal H$ by the uniform form bound and weak convergence,
density extends this identity to all $\varphi\in\mathcal H$.
The operator definition \eqref{def-Ltheta-s} therefore implies
\[
 f\in D(L_{\Theta_\nu}^s),\qquad L_{\Theta_\nu}^sf=0.
\]
The strong convergence of both $f_n$ and $\sin\theta_n$ in $L^2$
also gives
\[
\langle f,\sin\Theta_\nu\rangle_{L^2}=0.
\]
Lemma \ref{lem-interior-BPS-kernel} now implies $f=0$, contradicting
$\|f\|_{L^2}=1$.  This proves \eqref{near-BPS-ND}.
\end{proof}

\begin{remark}\label{remark-critical-cusp-not-covered}
Proposition \ref{prop-near-BPS-nondegeneracy} assumes
$\lambda\geq\lambda^{-}>0$.  If
$\omega/(1-\kappa)\to0$, Lemma \ref{lem-BPS-family} shows that the
selected BPS parameter may satisfy $\nu\to0$.  To recall the relevant
indicial equation, the exterior inverse-square model
\[
 A_{p,0}:=-\partial_r^2-r^{-1}\partial_r+p^2r^{-2}
\]
satisfies
\begin{equation}\label{eq:indicial-equation-cusp}
 A_{p,0}(r^\sigma)=(p^2-\sigma^2)r^{\sigma-2},
 \qquad\sigma=\pm p\text{ for a homogeneous power solution}.
\end{equation}
Here $p=1+\nu$ for the limiting BPS tail.  Thus the indicial exponent
$1+\nu$ approaches the borderline value $1$, and the uniform weighted
estimate \eqref{zero-mode-weighted-bound} is lost.  The present argument
therefore no longer rules out threshold escape through the intermediate
region $1\ll r\ll\omega^{-1/2}$.  The full $L^2$ nondegeneracy of the
critical BPS operator proved in Proposition
\ref{prop-BPS-nondegeneracy} does not by itself exclude such an escape.
\end{remark}

\subsection{Completion of the main results}
\begin{proof}[Proof of Theorem~\ref{thm:global-stability}]
Let $u$ be the maximal $H^2$ solution from Theorem~\ref{thm:lwp}.  By a
fixed target rotation we may first arrange
$\|u_0-\widetilde u\|_{H^1}=d_{H^1}(u_0,\mathcal O(\widetilde u))$ up to an
arbitrarily small error.  Rotation covariance then reduces the argument to
this representative.  On the maximal existence interval,
Theorem~\ref{thm:conditional-stability}, energy conservation, and mass
conservation imply that sufficiently small initial $H^1$ distance to
$\mathcal O(\widetilde u)$ keeps the solution in the $H^1$ tube required by
Proposition~\ref{prop:continuation}.  If the
forward maximal time $T^{+}$ were finite, that proposition would give
\[
\sup_{0\le t<T^{+}}\|u(t)-\vec n\|_{H^2}<\infty,
\]
contradicting the blow-up alternative in Theorem~\ref{thm:lwp}.  Hence
$T^{+}=\infty$.  The same argument backward in time gives $T^{-}=\infty$.
Applying the conditional stability estimate on arbitrary finite intervals
then yields the stated uniform-in-time $H^1$ orbital bound.  The finite-time
$H^2$ estimate is exactly Proposition~\ref{prop:continuation} applied on
$[-T,T]$.
\end{proof}

\begin{proof}[Proof of Corollary~\ref{cor:near-BPS}]
The variational theory of \cite{JHS} gives a scalar noncritical ground
state in the stated parameter range.  By
Lemma~\ref{lem:great-circle-reduction}, the scalar and full vector-valued
infima coincide, so this profile is also a ground state in the full
$m$--equivariant class.  Conversely, every ground state in the full class
is a constant target rotation of a scalar ground state.  Proposition
\ref{prop-near-BPS-nondegeneracy} supplies \emph{(ND)} for every such scalar
profile, uniformly for $\lambda\in \mathcal K\Subset(0,\lambda_*)$ once $\delta$ is
sufficiently small.  Constant target rotation conjugates the Hessian and
preserves the mass constraint, so \emph{(ND)} is invariant under this
rotation.  Thus every full $m$--equivariant ground state satisfies
\emph{(ND)}, and the conclusion follows from
Theorem~\ref{thm:global-stability}.
\end{proof}

\section*{Data availability statement}
No datasets were generated or analysed during the current study.

\section*{Conflict of interest}
The authors declare that they have no conflict of interest. 
The authors also declare that this manuscript has not been previously published.

\section*{Acknowledgment}
The work of G. Jin was supported by NSFC under Grant Number 12201542 and the Jilin Provincial international Joint Research Center for Differential Equations and Their Innovative Applications under Grant number YDZJ202602CXJD091.
The work of J. Seok is supported by the National Research Foundation of Korea (NRF) grant funded by the Korea government (MSIT) (No. RS-2025-00562375).

\appendix
\section{Instantaneous geometric pairings for arbitrary profiles}
\label{app:instantaneous-geometric}

The purpose of this appendix is to prove the \emph{precisely stated}
fixed-profile estimates in the proof of
Lemma~\ref{lem:pure-geometric-difference}: the short-geodesic Jacobi
bound \eqref{eq:main-Jacobi-bound}, the pointwise bounds for all four
remainders \eqref{eq:main-geometric-remainder-bounds}, and their
instantaneous $L^2$ pairings \eqref{eq:four-pairings-instantaneous}.
We also prove the auxiliary endpoint-velocity Jacobi estimate
\eqref{eq:main-velocity-Jacobi-bound}, used in the forcing estimate of
Lemma~\ref{lem:forcing-uniqueness}.  The final subsection proves the
low-regularity statement \eqref{eq:main-geometric-low-regularity} and
thereby justifies the $H^2$ energy identity.

The four-term decomposition is the one used in the parallel-transport
method of \cite{M,ShimizuUnique}.  We derive the differentiated transport
bound from the short-geodesic Jacobi boundary problem, without assuming
that either profile solves the unforced equation or that their initial
data agree.  Proposition~2.2 of \cite{ShimizuUnique} is not an input.

\subsection{Fixed-time notation and norm bounds}\label{app:fixed-norms}
Fix smooth $u,\bar u\in\mathcal X_m^2$, $m\ge2$, with
$\|u-\vec n\|_{H^2}+\|\bar u-\vec n\|_{H^2}\le R$ and
$\|u-\bar u\|_\infty<1$.  Retain $w,\psi_i,\bar\psi_i,\Pi$,
$\mathcal Z_i$, and $G$ from Section~\ref{sec:lwp-proof}.  Define the
nonnegative spatial functions
\[
 \begin{split}
 L(x)&=d_{\mathbb S^2}(\bar u(x),u(x)),\\
 \mathsf a&=|du|+|d\bar u|,\\
 \mathsf b&=|D^2u|+|D^2\bar u|+\mathsf a^2,\\
 \mathsf z&=\left(\sum_{i=1}^2|\mathcal Z_i|^2\right)^{1/2}.
 \end{split}
\]
Here $D^2$ is the ordinary Cartesian Hessian.  Since
$|u-\bar u|=2\sin(L/2)<1$, one has $L<\pi/3$ and $L\le C|w|$.
Equation~\eqref{eq:equivariant-gradient-infty} and two-dimensional
Sobolev embedding imply
\begin{equation}\label{eq:app-fixed-norms}
 \|\mathsf a\|_\infty+\|\mathsf a\|_4+\|\mathsf b\|_2
       +\|\mathsf z\|_4\le C_R,
 \qquad \|\mathsf z\|_2\le (2G)^{1/2}.
\end{equation}
Moreover $|\Pi-I|\le C|w|$ on tangent vectors at $\bar u$ and
$\partial_iw=\mathcal Z_i+(\Pi-I)\bar\psi_i$.  Hence
\begin{equation}\label{eq:app-energy-equivalence}
 C_R^{-1}\|w\|_{H^1}^2\le G\le C_R\|w\|_{H^1}^2,
 \qquad 0\le G\le C_R.
\end{equation}
The uniform equivariant comparison preceding
Lemma~\ref{lem:pure-geometric-difference}, followed by
\eqref{eq:app-energy-equivalence}, also gives
\[
 \|L\|_\infty\le C\|w\|_\infty
       \le C_R\|w\|_{H^1}\le C_RG^{1/2}.
\]
In particular the pointwise comparison for $w$ is available before any
evolution estimate is used.  These statements are spatial estimates;
they involve no time history.

\subsection{The Jacobi bound}
We first prove \eqref{eq:main-Jacobi-bound}.  Let $\gamma(s,x)=\Phi_s(\bar u(x),u(x))$ and use the pullback
connection $\nabla^\gamma$ along it.  Smooth dependence of the shorter
geodesic on its endpoints, uniformly for $L\le\pi/3$, gives
\begin{equation}\label{eq:app-geodesic-basic}
 |\gamma_s|=L,\qquad
 |\partial_j\gamma|\le C\mathsf a,\qquad
 |\nabla_j^\gamma\partial_k\gamma|\le C\mathsf b.
\end{equation}
For example, the geodesic is
\[
 \Phi_s(a,b)=\frac{\sin((1-s)L)}{\sin L}a
               +\frac{\sin(sL)}{\sin L}b,
 \qquad L=d_{\mathbb S^2}(a,b).
\]
The quotients have smooth even extensions in $L$ at zero, so this
formula is smooth on the relevant set of endpoint pairs, including the
diagonal.  Differentiating it in the endpoints proves
\eqref{eq:app-geodesic-basic} by the chain rule.

The more precise estimate needed when a derivative hits $\gamma_s$ is
\begin{equation}\label{eq:app-mixed-Jacobi}
 |\nabla_j^\gamma\gamma_s|
 \le C\bigl(|\mathcal Z_j|+L^2\mathsf a\bigr).
\end{equation}
We prove it explicitly.  For $L>0$ use a parallel orthonormal basis
along $\gamma$ whose first vector is $\gamma_s/L$.  The Jacobi field
$\partial_j\gamma$ has components $y_\parallel,y_\perp$ satisfying
\[
 y_\parallel''=0,\qquad y_\perp''+L^2y_\perp=0.
\]
Its endpoint vectors, identified by parallel transport, are
$A_j=\bar\psi_j$ and $B_j=\Pi_{0,1}\psi_j$; thus
$|B_j-A_j|=|\mathcal Z_j|$.  For either scalar component write its
endpoint values as $\alpha_0,\alpha_1$, to distinguish them from the
spatial size functions $\mathsf a,\mathsf b$.  The respective formulas are
\[
 y_\parallel(s)=(1-s)\alpha_0+s\alpha_1,
 \qquad
 y_\perp(s)=\frac{\sin((1-s)L)}{\sin L}\alpha_0
             +\frac{\sin(sL)}{\sin L}\alpha_1.
\]
For the perpendicular component,
\[
 y_\perp'(s)
 =\frac{L\cos(sL)}{\sin L}(\alpha_1-\alpha_0)
  +\frac{L\{\cos(sL)-\cos((1-s)L)\}}{\sin L}\alpha_0.
\]
The first coefficient is bounded and the second is $O(L^2)$,
uniformly for $0<L\le\pi/3$ and $0\le s\le1$.
Also $y_\parallel'=\alpha_1-\alpha_0$.  Torsion freeness gives
$\nabla_j^\gamma\gamma_s=\nabla_s^\gamma\partial_j\gamma$,
which proves \eqref{eq:app-mixed-Jacobi} and thus the statement
\eqref{eq:main-Jacobi-bound} in the body.  At $L=0$ the same estimate
follows by continuity of the endpoint formula; no division by $L$ is
used in the resulting inequality.

The same Jacobi boundary-value calculation proves the auxiliary
endpoint-velocity bound \eqref{eq:main-velocity-Jacobi-bound}.
Indeed, for arbitrary tangent endpoint velocities $\bar V,V$, the field
$\mathcal J[\bar V,V]$ in \eqref{eq:geodesic-velocity-linear} obeys
\begin{equation}\label{eq:app-velocity-Jacobi}
 |\mathcal J[\bar V,V](s)|\le C(|\bar V|+|V|).
\end{equation}
This proves \eqref{eq:main-velocity-Jacobi-bound}: the Jacobi field
solves the same linear boundary problem, with these velocities as
boundary values.  In particular it splits exactly when the endpoint
velocities are written as the sum of the pure and source velocities.
Neither this fact nor \eqref{eq:app-mixed-Jacobi} uses initial data.

\subsection{The four pointwise and pairing bounds}
We next prove \eqref{eq:main-geometric-remainder-bounds} and
\eqref{eq:four-pairings-instantaneous}.  Write $F=u\times\Delta u$, $\bar F=\bar u\times\Delta\bar u$,
so $|F|+|\bar F|\le C\mathsf b$.
The transport identity \eqref{eq:transport-commutator-explicit}, also
on subsegments, implies
\begin{equation}\label{eq:app-transport-bounds}
 |\mathcal C_jY|\le CL\mathsf a|Y|,
 \qquad
 |\mathcal C[\bar F,F]Y|\le CL\mathsf b|Y|.
\end{equation}
Equations \eqref{eq:app-geodesic-basic} and
\eqref{eq:app-velocity-Jacobi} justify these estimates.  The first two
remainders in \eqref{eq:geometric-forced-difference} therefore satisfy
\begin{equation}\label{eq:app-R01-pointwise}
 |\mathcal R_{0,i}|+|\mathcal R_{1,i}|
       \le CL\mathsf a\mathsf b.
\end{equation}
Here $|\nabla_j^{\bar u}\bar\psi_i|\le |D^2\bar u|$.
Transport preserves curvature, and the curvature expression is cubic
in the first derivatives.  Subtracting its three factors in turn gives
\begin{equation}\label{eq:app-R3-pointwise}
 |\mathcal R_{3,i}|\le C\mathsf a^2\mathsf z.
\end{equation}

For completeness, the differentiated transport term is expanded before
being estimated.  Put
$A_s=\gamma_s$, $B_{j,s}=\partial_j\gamma$, and
$Y_{i,s}=\Pi_{s,0}\bar\psi_i$.  For a field $Y$ at parameter $a$,
define
\[
 \mathcal C_j^{b,a}Y
 =\nabla_j^\gamma(\Pi_{b,a}Y)
     -\Pi_{b,a}\nabla_j^\gamma Y.
\]
Thus $|\mathcal C_j^{b,a}Y|\le CL\mathsf a|Y|$ for every subsegment.
Since the curvature tensor of the unit sphere is parallel, differentiating
the integral for $\mathcal C_j\bar\psi_i$ gives exactly
\begin{align*}
 \nabla_j^u(\mathcal C_j\bar\psi_i)=\int_0^1\bigl\{&
 \mathcal C_j^{1,s}\bigl(R(A_s,B_{j,s})Y_{i,s}\bigr)\\
 &+\Pi_{1,s}R(\nabla_j^\gamma A_s,B_{j,s})Y_{i,s}\\
 &+\Pi_{1,s}R(A_s,\nabla_j^\gamma B_{j,s})Y_{i,s}\\
 &+\Pi_{1,s}R(A_s,B_{j,s})\mathcal C_j^{s,0}\bar\psi_i\\
 &+\Pi_{1,s}R(A_s,B_{j,s})
                    \Pi_{s,0}\nabla_j^{\bar u}\bar\psi_i
 \bigr\}\,ds.
\end{align*}
In this display $R$ denotes $R^{\mathbb S^2}$, not the Sobolev-ball
radius.  By \eqref{eq:app-geodesic-basic},
\eqref{eq:app-mixed-Jacobi}, and \eqref{eq:app-transport-bounds},
the five integrands are bounded, respectively, by
\[
 CL^2\mathsf a^3,\quad
 C\mathsf a^2(|\mathcal Z_j|+L^2\mathsf a),\quad
 CL\mathsf a\mathsf b,\quad
 CL^2\mathsf a^3,\quad CL\mathsf a\mathsf b.
\]
Thus
\begin{equation}\label{eq:app-R2-pointwise}
 |\mathcal R_{2,i}|
 \le C\bigl(\mathsf a^2\mathsf z
             +L\mathsf a\mathsf b+L^2\mathsf a^3\bigr).
\end{equation}
Together with \eqref{eq:app-R01-pointwise} and
\eqref{eq:app-R3-pointwise}, this completes the proof of all three lines
of \eqref{eq:main-geometric-remainder-bounds}.
In particular, no derivative of $\mathcal Z_i$ has to be bounded in
$L^2$ by the difference energy.  This is the purpose of the Jacobi
calculation, rather than merely assuming that bounded gradients control
a differentiated transport commutator.

We now perform the spatial pairings, still at a single fixed time.
Equations \eqref{eq:app-R01-pointwise}--\eqref{eq:app-R2-pointwise}
imply
\begin{equation}\label{eq:app-pairing-integrals}
 \sum_{i,\alpha}|\langle\mathcal R_{\alpha,i},\mathcal Z_i\rangle|
 \le C\int_{\mathbb R^2}
 \bigl(\mathsf a^2\mathsf z^2
       +L\mathsf a\mathsf b\mathsf z
       +L^2\mathsf a^3\mathsf z\bigr)\,dx.
\end{equation}
The first term is at most $C_RG$.  Since $L<\pi/3$,
$\|L\|_2\le C\|w\|_2$, and $\|\mathsf a\|_\infty\le C_R$,
the last term is also at most $C_R\|w\|_2\|\mathsf z\|_2\le C_RG$.
For the middle term use the equivariant bound
$\|L\|_\infty\le C_RG^{1/2}$ from Section~\ref{app:fixed-norms}.
H\"older's inequality now gives directly
\begin{equation}\label{eq:app-middle-pairing}
 \int L\mathsf a\mathsf b\mathsf z\,dx
 \le \|L\|_\infty\|\mathsf a\|_\infty
                   \|\mathsf b\|_2\|\mathsf z\|_2
 \le C_RG.
\end{equation}
All three integrals in \eqref{eq:app-pairing-integrals} are therefore
bounded by $C_RG$, proving the second inequality in
\eqref{eq:four-pairings-instantaneous}; the first inequality follows
from \eqref{eq:app-pairing-integrals}.
The Jacobi calculation \eqref{eq:app-mixed-Jacobi} and its consequence
\eqref{eq:app-R2-pointwise} are essential: only after eliminating
uncontrolled derivatives of $\mathcal Z_i$ can this $L^\infty$ pairing
be used.  No auxiliary exponent or logarithmic optimization is needed.
No time supremum, integration in time, equality of initial data, or
solution of an unforced equation has been used in this proof.

\subsection{The energy identity at \texorpdfstring{$H^2$}{H2} regularity}
We prove \eqref{eq:main-geometric-low-regularity} and the associated
energy identity.  Let now $u,\bar u$ be paths in the solution class of
Theorem~\ref{thm:lwp}, on a compact time interval on which the $H^2$
bound and short-geodesic condition above hold.  This paragraph uses
only the regularity $C_tH^2\cap C_t^1L^2$ and the distributional forced
equation, not uniqueness or continuous dependence of the solution map.

For $Y\cdot\bar u=0$, the transport is given by
\eqref{eq:geodesic-transport-first-use}; in the present notation this is
\begin{equation}\label{eq:app-transport-matrix}
 \Pi Y=Y-\frac{Y\cdot u}{1+u\cdot\bar u}(u+\bar u).
\end{equation}
The denominator is bounded below by $3/2$ under $|u-\bar u|<1$.
Consequently the matrix coefficients of $\Pi$ belong to
$L_t^\infty W_x^{1,\infty}$, and their time derivatives belong to
$L_t^\infty L_x^2$.  We used
\eqref{eq:equivariant-gradient-infty} for the spatial derivative and
$u_t,\bar u_t\in L_t^\infty L_x^2$ for the time derivative.
The chain and product rules then give
\[
 \mathcal Z_i\in C_tH_x^1,
 \qquad \partial_t\mathcal Z_i\in L_t^\infty H_x^{-1}.
\]
Indeed, in
\[
 \partial_t\mathcal Z_i
 =\partial_i u_t-\Pi\partial_i\bar u_t
                    -(\partial_t\Pi)\partial_i\bar u
\]
the first two terms are in $H^{-1}$, since multiplication by a
$W^{1,\infty}$ matrix is bounded there; the last is in $L^2$, since
$\partial_i\bar u\in L^\infty$.
This proves the first two statements of
\eqref{eq:main-geometric-low-regularity}.  The
$H^1$--$L^2$--$H^{-1}$ energy chain rule therefore makes
$\|\mathcal Z_i(t)\|_2^2$ absolutely continuous and gives
\[
 \frac12\frac d{dt}\|\mathcal Z_i\|_2^2
 =\langle\partial_t\mathcal Z_i,\mathcal Z_i\rangle_{H^{-1},H^1}
 =\langle\nabla_t^u\mathcal Z_i,\mathcal Z_i\rangle_{H^{-1},H^1}.
\]
The connection and transport identities follow by spatial approximation
from the smooth formulas.  They can also be differentiated directly in
\eqref{eq:app-transport-matrix}; the coefficients are $W^{1,\infty}$,
the second spatial derivatives are in $L^2$, and the endpoint velocities
are in $L^2$, so every product in the identities is defined in
$H^{-1}$.  In particular all four remainders are in $L^2$ by
\eqref{eq:app-R01-pointwise}--\eqref{eq:app-R2-pointwise} and
\eqref{eq:app-fixed-norms}.  Approximation passes their estimates to
$H^2$ profiles with the same constants.

The principal energy pairing is zero: metric compatibility and
$\nabla^u J=0$ give
\[
 \left\langle-J\sum_j(\nabla_j^u)^2\mathcal Z_i,
                            \mathcal Z_i\right\rangle
 =\sum_j\int
 \langle J\nabla_j^u\mathcal Z_i,\nabla_j^u\mathcal Z_i\rangle\,dx=0.
\]
The zeroth-order energy has the usual $L^2$ chain rule.  Hence $G$
is absolutely continuous in time, completing the justification of
\eqref{eq:main-geometric-low-regularity}.  For the
forced flow write $Z=\Gamma[u](\vec n\times u)$ and
$\bar Z=\Gamma[\bar u](\vec n\times\bar u)$.  Applying
\eqref{eq:geometric-forced-difference} now gives the exact
almost-everywhere decomposition
\begin{equation}\label{eq:app-exact-energy-split}
 \begin{split}
 G'_{\rm SM}&=-\sum_j\int w\cdot
       (\partial_j\bar u\times\partial_jw)\,dx
       +\sum_{i,\alpha}\langle\mathcal R_{\alpha,i},\mathcal Z_i\rangle,\\
 G'_{\rm CS}&=\int w\cdot(Z-\bar Z)\,dx
                     +\sum_i\langle\mathcal S_i,\mathcal Z_i\rangle,\\
 G'&=G'_{\rm SM}+G'_{\rm CS}.
 \end{split}
\end{equation}
The first line obeys \eqref{eq:pure-geometric-bound}; the second is
estimated in Lemma~\ref{lem:forcing-uniqueness}.  Integrating the
resulting almost-everywhere inequality is legitimate by the established
absolute continuity.  The initial value $G(t_0)$ remains arbitrary
throughout; it is set to zero only for the uniqueness consequence.

\raggedbottom
\bibliographystyle{amsplain}

\end{document}